\documentclass[12pt,oneside, double]{CUNY_PhD} 
 
\usepackage{amsfonts} 
\usepackage{amssymb} 
\usepackage{amsmath}
\usepackage{amsthm} 
\usepackage{amscd} 
\usepackage{verbatim} 
\usepackage[titletoc]{appendix}

\newtheorem{Thm}{Theorem}
\newtheorem{Cor}[Thm]{Corollary}
\newtheorem{Lemma}[Thm]{Lemma} 
\newtheorem{Sublemma}{Lemma}[Thm] 
\newtheorem{Def}[Thm]{Definition}

\newtheorem{Observation}[Thm]{Observation}
\newtheorem*{TestQuestion}{Test Questions}
\newtheorem{Question}[Thm]{Question}
\newcounter{saveenum} 

\def\<#1>{\langle#1\rangle} 
\def\(#1){\left(#1\right)} 
\def\bval#1{\|#1\|} 
\def\btimes{\cdotp} 
\def\forces{\Vdash} 
\def\proves{\vdash} 
\def\nforces{\nVdash} 
\def\iff{\leftrightarrow} 
\newcommand{\from}{\mathbin{\vbox{\baselineskip=2pt\lineskiplimit=0pt
                         \hbox{.}\hbox{.}\hbox{.}}}} 

\def\class#1{\mathbb #1} 
\def\classdot#1{\dot{\mathbb #1}} 
\def\corners#1{\ulcorner{#1}\urcorner} 

\newcommand{\cof}{\mathop{\rm cof}} 
\newcommand{\dom}{\mathop{\rm dom}} 
\newcommand{\ran}{\mathop{\rm ran}} 
\newcommand{\tc}{\mathop{\rm tc}} 
\newcommand{\rank}{\mathop{\rm rank}} 
\newcommand{\supt}{\mathop{\rm support}} 
\newcommand{\op}{\mathop{\rm op}} 
\newcommand{\Add}{\mathop{\rm Add}} 
\newcommand{\Coll}{\mathop{\rm Coll}} 
\newcommand{\Con}{\mathop{\rm Con}} 

\renewcommand{\P}{{\mathbb P}} 
\newcommand{\Q}{{\mathbb Q}} 
\newcommand{\B}{{\mathbb B}} 
\newcommand{\D}{{\mathbb D}} 
\newcommand{\I}{{\mathcal I}} 
\newcommand{\M}{{\mathcal M}} 
\newcommand{\V}{{\mathcal V}} 
\newcommand{\W}{{\mathcal W}} 

\newcommand{\GA}{\hbox{\sc ga}} 
\newcommand{\BA}{\hbox{\sc ba}} 
\newcommand{\GCH}{\hbox{\sc gch}} 
\newcommand{\AC}{\hbox{\sc ac}} 
\newcommand{\ORD}{\hbox{\sc ord}} 
\newcommand{\ZFC}{\hbox{\sc zfc}} 
\newcommand{\ZC}{\hbox{\sc zc}} 
\newcommand{\BGC}{\hbox{\sc bgc}} 
\newcommand{\BG}{\hbox{\sc bg}} 

\newcommand{\VHOD}{\hbox{\sc v=hod}} 
\newcommand{\nVHOD}{\hbox{\sc v$\neq$hod}} 
\newcommand{\CCA}{\hbox{\sc cca}} 
\newcommand{\PFA}{\hbox{\sc pfa}} 

\title{The Ground Axiom}
\author{Jonas Reitz}
\date{August, 2006}

\begin{document}

\frontmatter

\maketitle

\makecopyrightpage
\makeapprovalpage{Sergei Artemov}{Melvin Fitting}{Joel David Hamkins}{Attila Mate}{Joel David Hamkins}{Jozef Dodziuk}
\makeabstractpage{Joel David Hamkins}{A new axiom is proposed, the Ground Axiom, asserting that the universe is not a nontrivial set-forcing extension of any inner model.  The Ground Axiom is first-order expressible, and any model of $\ZFC$ has a class-forcing extension which satisfies it.  The Ground Axiom is independent of many well-known set-theoretic assertions including the Generalized Continuum Hypothesis, the assertion $\VHOD$ that every set is ordinal definable, and the existence of measurable and supercompact cardinals.  The related Bedrock Axiom, asserting that the universe is a set-forcing extension of a model satisfying the Ground Axiom, is also first-order expressible, and its negation is consistent.  As many of these results rely on forcing with proper classes, an appendix is provided giving an exposition of the underlying theory of proper class forcing.}

\chapter*{Acknowledgements}
First and foremost I would like to thank my advisor, Joel David Hamkins, an inspiration in mathematics and in life.  He sets high standards for his students and conforms to even higher standards himself.  An aspiring mathematician could not ask for a better mentor. Thank you also to my committee for their time, attention, and insightful comments.  My close friends and respected colleagues Victoria Gitman and Thomas Johnstone, who have shared five years of intense and wonderful study with me, deserve much credit for enabling me to produce this work.  Thank you for listening to my ideas and sharing your own.  Learning has never been such an adventure and such a pleasure as it has been with you.  Thank you to my family for instilling the crazy idea that I could accomplish whatever I wanted, and for teaching me about love.  Finally, to my wife and best friend Gwen I owe the gratitude of uncomplaining support and unconditional love through far too many trials.  You are amazing.

\tableofcontents

\mainmatter

\chapter*{Introduction}

The main contribution of this dissertation consists of Chapter \ref{C:TheGroundAxiom}.  Chapter \ref{C:GAandVnotHOD} represents further joint work, which is included because it relates very strongly to the results in Chapter \ref{C:TheGroundAxiom}.  The appendix is a largely supplemental review of material that is known, but which has rarely been written down in detail.  It is included for completeness.

\chapter{The Ground Axiom}\label{C:TheGroundAxiom}

Forty years of forcing has illustrated its efficacy and flexibility in producing models of set theory with a wide variety of properties. Each such example further enriches our knowledge of the collection of models obtainable by forcing.  This collection possesses such diversity of models and such intricate structure between them that its exploration will doubtless continue into the forseeable future.  

I am interested in exploring the boundaries and limitations of this collection, by considering under what circumstances is the universe \emph{not} a forcing extension.  This motivates the following axiom, formulated jointly with Joel Hamkins.  This axiom is proposed not in the sense of a ``self-evident truth,'' but rather as a proposition identifying an interesting and important feature that a model of set theory may exhibit.

\clearpage
\begin{Def}\rm\label{D:GAformal}The \emph{Ground Axiom} ($\GA$) is the assertion that the universe of sets $V$ is not a forcing extension of any inner model $W$ by nontrivial forcing $P\in W$.
\end{Def}

Some observations are in order.  As stated, the Ground Axiom is second order in nature, requiring quantification over classes $W$.  However,  I will show in Section \ref{S:GAExpressible} that the Ground Axiom has a first-order equivalent.  In addition, the Ground Axiom refers only to \emph{nontrivial} forcing, else any model is a forcing extension of itself, and only to \emph{set} forcing, a restriction that will be important in the first-order expressibility results.

There are a number of well known models of the Ground Axiom, including many of the `canonical models.'    

\begin{Observation}\label{T:CanonicalModelsofGA} The constructible universe $L$, the model $L[0^\#]$, and the canonical model of a measurable cardinal $L[\mu]$ all satisfy the Ground Axiom.
\end{Observation}

\begin{proof}In each case, the result follows from the uniqueness and minimality properties of the model.  If $L=W[g]$ were a forcing extension of an inner model $W$, then by the absoluteness of $L$ we would have $L\subset W$, and so $L=W$.  For $L[0^\#]$, suppose $L[0^\#]=W[g]$.  It is well known that $0^\#$ cannot be created by set forcing (for example, this appears in Jech \cite{Jech:SetTheory3rdEdition} as exercise 18.2, page 336), and so $0^\# \in W$.  Thus $L[0^\#]\subset W$, and so $W=L[0^\#]$.

In the case of $L[\mu]$, suppose $L[\mu]=W[g]$, where $g$ is $W$-generic for a poset $Q\in W$.  Fix a $Q$-name $\dot{\mu}$ such that $(\dot{\mu})_g = \mu$, and $q\in Q$ such that $ q \forces \dot{\mu}$ is a normal measure on $\check{\kappa}$.  I claim that for every $A \subset \kappa$ in $W$, either $q \forces \check{A}\in\dot{\mu}$ or $q \forces \check{A}\notin\dot{\mu}$.  If this is not the case, then there are $q_0,q_1 \leq q$ with $q_0 \Vdash \check{A}\in\dot{\mu}$ and $q_1 \Vdash \check{A}\notin\dot{\mu}$.  Let $g_0 \times g_1$ be $W$-generic for $Q\times Q$ such that $\<q_0,q_1>\in g_0 \times g_1$.  Then $W[g_0] \vDash A \in  (\dot{\mu})_{g_0}$ and $W[g_1] \vDash A \notin  (\dot{\mu})_{g_1}$.  Work in $W[g_0]$ and $W[g_1]$ to build $L[(\dot{\mu})_{g_0}]$ and $L[(\dot{\mu})_{g_1}]$ respectively.  The uniqueness of the model $L[\mu]$ implies that $L[(\dot{\mu})_{g_0}]=L[(\dot{\mu})_{g_1}]=L[\mu]$, and uniqueness of the normal measure in $L[\mu]$ implies that $(\dot{\mu})_{g_0}=(\dot{\mu})_{g_1}=\mu$. However, $A\in(\dot{\mu})_{g_0}$ but $A\notin (\dot{\mu})_{g_0}$, a contradiction. This shows that $W\cap \mu$ is definable in  $W$ as $\{A\subset \kappa \mid q \forces \check{A}\in\dot{\mu}\}$, and so $\kappa$ is measurable in $W$.  As $L[\mu]$ is the minimal model in which $\kappa$ is measurable, $L[\mu]\subset W$, and so $L[\mu] = W$.  Thus the forcing adding $g$ was trivial.
\end{proof}

In addition, in many cases the core model $K$ satisfies the Ground Axiom.  However, certain other canonical models do not. For example, Ralph Schindler observed that methods of Woodin show that the least model $M_1$ of one Woodin cardinal is in fact a forcing extension of an inner model.   The similarities between the known models of the Ground Axiom suggest many questions about the consequences of the axiom.  As a starting point, we can consider the relationship of the Ground Axiom to the Generalized Continuum Hypothesis ($\GCH$) and to the assertion that every set is ordinal definable ($\VHOD$).

\begin{TestQuestion}\mbox{ }
\begin{enumerate}
\item \label{Ques:GA->GCH} Does the $\GA$ imply the $\GCH$?
\item \label{Ques:GA->HOD} Does the $\GA$ imply $\VHOD$?
\end{enumerate}
\end{TestQuestion}

Both questions have negative answers.  The former is answered in Corollary \ref{C:GAand-GCH} below, and the latter is answered in Chapter \ref{C:GAandVnotHOD}.

Finally, the examples mentioned above leave open the question of compatibility of the Ground Axiom with various large cardinal hypotheses not covered by the Core Models, such as $\GA +$\emph{ there is a supercompact cardinal}.  While the Ground Axiom is not implied by any large cardinal axiom, as such axioms are generally preserved by small forcing, the relative consistency of the Ground Axiom with large cardinals remains to be considered.  In Section \ref{S:ForcingGA}, I will demonstrate  a method for building models of the Ground Axiom that accomodates many large cardinals, giving relative consistency of the Ground Axiom with supercompact cardinals among others.  This method will also show the consistency of $\GA + \neg \GCH$.  In Section \ref{S:GAandGCH}, I consider an adaptation of the method to produce diverse models of $\GA + \GCH$.  In Section \ref{S:BedrockAxiom}, I will turn my attention to a related notion, the \emph{Bedrock Axiom} ($\BA$), which asserts that either the universe is a model of the Ground Axiom (a `bedrock model') or is a forcing extension of a such a model.  I will show the consistency of the negation of the Bedrock Axiom.  In Chapter \ref{C:GAandVnotHOD}, I will discuss a very different method of producing models of the Ground Axiom which shows the consistency of $\GA+\nVHOD$.

Throughout the following I will use blackboard bold $\P$ for proper class partial orders and standard text $P$ for set partial orders.  In function definitions, three dots $\from$ indicates that a function may be a partial function, e.g. $f\from \ORD \rightarrow \ORD$.  The partial order $\Add(\gamma,\delta)$ consists of all functions $f\from \gamma\times\delta \rightarrow 2$ of size $<\gamma$, and $\Coll(\gamma,\delta)$ consists of functions $f\from \gamma \rightarrow \delta$ of size $<\gamma$. In describing initial segments of models other than $V$, such as forcing extensions, I will use $V[G]_\alpha$ to denote $\(V_\alpha)^{V[G]}$.

\section[The Ground Axiom is first-order expressible]{The Ground Axiom is first-order \\expressible}\label{S:GAExpressible}

I begin by showing that these notions are, indeed, first-order expressible.  

\begin{Thm}\label{T:GAIsFirstOrder}
There is a first-order formula which holds in a model of set theory exactly when that model is not a forcing extension of any inner model by nontrivial set forcing.
\end{Thm}

Theorem \ref{T:GAIsFirstOrder} is a consequence of Theorem \ref{T:GADefinition}, which gives a more detailed result.  A similar result is implicit in independent work of Woodin \cite{Woodin:GenericMultiverse}.

The formula asserting the Ground Axiom will be given explicitly below, but some definitions are required.  The first are the \emph{$\delta$ cover and $\delta$ approximation properties}, formulated by Hamkins \cite{Hamkins2003:ExtensionsWithApproximationAndCoverProperties}, which provide a framework for analyzing  extensions and inner models.

\begin{Def}\rm\label{D:deltacover&approx}\emph{(Hamkins)}.  Suppose that $W\subseteq V$ are transitive models of (some fragment of) $\ZFC$, and $\delta$ is a cardinal in $V$.
\begin{enumerate}
\item $\<W,V>$ has the \emph{$\delta$ cover property} if and only if for each $A\in V$ with $A\subset W$ and $|A|^{V} < \delta$ there is a covering set $B$ in $W$ with $A\subset B$ and $|B|^{W} < \delta$.
\item $\<W,V>$ has the \emph{$\delta$ approximation property} if and only if for each $A\in V$ with $A\subset W$, if $A\cap B\in W$ for every $B\in W$ with $|B|^{W} < \delta$, then $A\in W$.
\end{enumerate}
\end{Def}

As I will be working in initial segments of the universe, I will also need an appropriate variant of $\ZFC$.

\clearpage
\begin{Def}\rm Let $\ZFC_{\delta}$ be the theory consisting of Zermelo Set Theory, Choice,  and $\leq\!\delta$-Replacement (that is, Replacement holds for functions with domain $\delta$, a regular cardinal), together with the axiom $$(*) \quad \forall A \; \exists \alpha \in ORD \; \exists E\subseteq\alpha\times\alpha \thickspace \thickspace \<\alpha, E> \cong \<\tc(\{A\}), \in>$$ which asserts ``every set is coded by a set of ordinals.''
\end{Def}

Formally, $\ZFC_\delta$ is a theory in the language of set theory together with a symbol for $\delta$, and includes the assertion that $\delta$ is a regular cardinal.  In any model $V$ of $\ZFC$, if $\gamma$ is a $\beth$-fixed point of cofinality $>\!\delta$, then $V_\gamma \vDash \ZFC_\delta$.  The main result of this section is given in Theorem \ref{T:GADefinition}, providing an explicit first-order statement $\Phi(\delta,z,P,G)$ which holds if and only if the universe is a set-forcing extension of an inner model (that is, if and only if $V\vDash \neg \GA$).

\clearpage
\begin{Thm}\label{T:GADefinition} The Ground Axiom is first-order expressible.\\
Specifically, the Ground Axiom fails if and only if there exist $\delta, z, P,$ and $G$ satisfying the following statement.
\begin{description}
	\item[$\Phi(\delta,z,P,G):$] \mbox{ }\\$\delta$ is a regular cardinal, $P\in z$ is a poset of size $<\delta$, $G$ is $z$-generic for $P$, and for every $\beth$-fixed point $\gamma > \delta$ of cofinality $>\delta$, there exists a transitive structure $\M$ of height $\gamma$ such that:
		\begin{enumerate}
		\parskip=0pt
		\item \label{T:GADefinition:ZFCdelta}$\M$ is a model of $\ZFC_{\delta}$,
		\item \label{T:GADefinition:z=P(delta)}$z = \mathcal{P}(\delta)^\M$,
		\item \label{T:GADefinition:forcingext}$\M[G] = V_{\gamma}$, and
		\item \label{T:GADefinition:CoverandApproximation}$\M \subset V_{\gamma}$ satisfies the $\delta$ cover and $\delta$ approximation properties.
	\end{enumerate}
\end{description}
\end{Thm}

This result is closely related to a theorem of Laver \cite{Laver:CertainVeryLargeCardinalsNotCreated}, who showed that every model of set theory is a definable class in all of its set-forcing extensions (independently observed by Woodin \cite{Woodin:GenericMultiverse}).  The two directions will be stated and proven separately as Lemmas \ref{L:GADefinition1} and \ref{L:GADefinition2}.

\begin{Sublemma}\label{L:GADefinition1} If the Ground Axiom fails, then there exist $\delta, z, P,$ and $G$ satisfying $\Phi(\delta,z,P,G)$.
\end{Sublemma}

\begin{proof}
Suppose $V=W[G]$ is a forcing extension of $W$ by a poset $P\in W$. Let $\delta = (|P|^+)^V$ and let $z=\mathcal{P}(\delta)^W$.  I will argue that $\Phi(\delta,z,P,G)$ holds.  Fix $\gamma > \delta$ a $\beth$-fixed point of cofinality $>\!\delta$.  I will show that $W_\gamma$ witnesses the properties of the structure $\M$ for $\gamma$.  For property \ref{T:GADefinition:ZFCdelta}, note that $W_\gamma \vDash \ZC$ (Zermelo set theory and the Axiom of Choice) for any limit ordinal $\gamma$.  As $\gamma$ is a $\beth$-fixed point, every set in $W_\gamma$ has transitive closure of size $<\gamma$ and so is coded as a set of ordinals in $W_\gamma$.  That $W_\gamma$ satisfies $\leq\!\delta$-Replacement follows from the cofinality of $\gamma$, for any function $f: a\to W_\gamma$ with $|a|=\delta$ has its range contained in some $W_\beta$ for $\beta < \gamma$, and so $\ran(f)\in W_{\beta+1}$.  Thus $W_\gamma \vDash \ZFC_\delta$.  That $z=\mathcal{P}(\delta)^{W_\gamma}$ holds by definition, and so $W_\gamma$ witnesses property \ref{T:GADefinition:z=P(delta)}.  To see that $W_{\gamma}[G] = V_\gamma$, note that any $x\in W[G]$ of rank $<\!\gamma$ will have a name $\dot{x}\in W$ of rank $<\!\gamma$, provided $\gamma$ is a limit ordinal greater than $\rank(P)$.

It remains to consider the $\delta$ cover and $\delta$ approximation properties.  I first claim that $W \subset V$ satisfies these properties.  In fact, any forcing extension by forcing of size less than $\delta$ will satisfy the $\delta$ cover and $\delta$ approximation properties.  This is a special case of a more general result of Hamkins \cite{Hamkins2003:ExtensionsWithApproximationAndCoverProperties}, which is presented in detail in Chapter \ref{C:GAandVnotHOD}, Lemma \ref{L:closurepointforcing}.  That these properties are inherited by $W_\gamma \subset V_\gamma$ is straightforward.  For example, if a set $A\in V_\beta$ for some $\beta < \gamma$ is covered by a set $B\in W$, then taking $B\cap W_\beta$ gives a cover of $A$ in $W_\gamma$.  A similar argument shows that for a set in $V_\gamma$, if every $W_\gamma$ approximation of that set over $V_\gamma$ lies in $W_\gamma$, then every $W$ approximation of that set over $V$ lies in $W$.  Thus $W_\gamma \subset V_\gamma$ satisfies the $\delta$ cover and $\delta$ approximation properties.  
\end{proof}

For the reverse direction, a key step in the argument is that for each $\gamma$ the witnessing structure $\M$ is unique.  This is an application of a result of Laver \cite{Laver:CertainVeryLargeCardinalsNotCreated}, who showed that the power set of $\delta$ together with the $\delta$ cover and $\delta$ approximation properties uniquely determine an inner model.

\begin{Sublemma}\emph{(Laver)}\label{L:M=M'}
Suppose $\M, \M^\prime$ and $V$ are transitive models of $\ZFC_{\delta}$ for $\delta$  a regular cardinal of $V$, both $\M$ and $\M^\prime$ are submodels of $V$ for which $\M\subseteq V$ and $\M^\prime \subseteq V$ satisfy the $\delta$ cover and $\delta$ approximation properties, $\M$ and $\M^\prime$ have the same power set of $\delta$,  and $(\delta^{+})^{\M} = (\delta^{+})^{\M^\prime} = (\delta^{+})^{V}$.  Then $\M=\M^\prime$.
\end{Sublemma}

\begin{proof} (Hamkins)
Since $\M$ and $\M^\prime$ satisfy axiom $(*)$, in order to conclude $\M=\M^\prime$ it suffices to show that $\M$ and $\M^\prime$ have the same sets of ordinals.  Observe that for any $A\subset ORD$, the statement ``$|A|<\delta$'' is unambiguous between the models $\M, \M^\prime$, and $V$, for by the $\delta$ cover property if it is true in one, then it will be true in all others containing $A$.  Also, since $(\delta^{+})^{\M} = (\delta^{+})^{\M^\prime} = (\delta^{+})^{V}$ the statement $|A|=\delta$ is unambiguous as well.

I will first establish a certain simultaneous cover property, namely, that for every set of ordinals $A\in V$ of size $<\delta$ there is a set $B\in \M\cap \M^\prime$ of size $\leq\delta$ with $A\subseteq B$.  Let $A\subset \alpha$ be such a set and fix in $V$ a well-ordering $\preccurlyeq$ of $\mathcal{P}(\alpha)^{V}$.  Construct in $V$ a sequence $\<B_\xi \mid \xi < \delta>$ of subsets of $\alpha$, each of size $<\delta$.  Let $B_0 = A$.  If $B_\xi \in \M$, then let $B_{\xi+1} \in \M^\prime$ be the $\preccurlyeq$-least subset of $\alpha$ such that $B_{\xi+1}\supset B_\xi$ and $|B_{\xi+1}|<\delta$.  The $\delta$ cover property for $\<\M^\prime,V>$ guarantees the existence of such a set.  If $B_\xi \notin \M$, then let $B_{\xi+1} \in \M$ to be the  $\preccurlyeq$-least subset of $\alpha$ with $B_{\xi+1}\supset B_\xi$, and $|B_{\xi+1}|<\delta$.  For limit $\xi$, let $B_\xi = \bigcup_{\beta<\xi} B_\beta$ and note that $B_\xi \in V$ by $\leq\delta$-Replacement.  As $\delta$ is regular and each $B_\beta$ has size $<\delta$, I have $|B_\xi|<\delta$ as well.  Thus the construction can continue through every stage $\xi < \delta$.  Note that $B_\xi\in \M$ for cofinally many $\xi$, and this is true for $\M^\prime$ as well.  In $V$, let $B = \bigcup_{\xi<\delta} B_\xi$.  To show $B\in \M$, I will show that every $\delta$ approximation of $B$ is in $\M$.  If $C\in \M$ has size $<\delta$, then by regularity of $\delta$ the intersection $C\cap B$ is equal to $C\cap B_\xi$ for some $\xi<\delta$.  Without loss of generality $B_\xi\in \M$ as this occurs cofinally in the sequence, and so $C\cap B_\xi \in \M$.  Thus $C\cap B \in \M$ for every set $C$ of size $<\delta$ in $\M$, and so by the $\delta$ approximation property $B\in \M$.  The same proof shows that $B\in \M^\prime$.  This establishes the desired simulataneous cover property.

I next claim that $\M$ and $\M^\prime$ have the same sets of ordinals of size $<\delta$.  Fix $A$ a set of ordinals in $\M$ of size $<\delta$.  The simultaneous cover property gives $B\in \M\cap \M^\prime$ of size $\leq\delta$ with $A\subseteq B$.  Since $B$ has order type smaller than $\delta^{+}$ there is a well-ordering $w\in \M$ of $\delta$, or possibly of a subset of $\delta$, of order type $ot(B)$.  Since  $w\subset \delta\times\delta$, the assumption that $\M$ and $\M^\prime$ have the same subsets of $\delta$ implies that $w$ must be in $\M^\prime$ as well.  The ordering $w$ induces an enumeration $\{b_\alpha \mid \alpha <\delta^\prime\}$ of $B$ for some $\delta^\prime\leq\delta$, and this enumeration is in $\M^\prime$ by $\leq\delta$-Replacement.  The set $\{\alpha <\delta^\prime \mid b_\alpha\in A\}$ is definable in $\M$, and since $\mathcal{P}(\delta)^\M=\mathcal{P}(\delta)^{\M^\prime}$ it exists in $\M^\prime$ as well.  As $A$ is definable from this set, $B$, and $w$, it must be the case that $A\in \M^\prime$.  The same argment shows that every set of ordinals in $\M^\prime$ of size $<\!\delta$  is also in $\M$.

Finally, it remains to show that $\M$ and $\M^\prime$ have the same sets of ordinals.  Fix $A$ a set of ordinals in $\M$.  I will show that every $\delta$ approximation of $A$ is in $\M^\prime$.  Fix $B \in \M^\prime$ of size $<\delta$.  The claim above shows that $B\in \M$ and so $A\cap B \in \M$, and applying the claim once more shows that $A\cap B \in \M^\prime$.   Thus every $\delta$ approximation of $A$ is in $\M^\prime$, and so $A\in \M^\prime$.  This shows $\M\subset \M^\prime$, and the reverse inclusion follows by the same argument, yielding $\M=\M^\prime$. 
\end{proof}

\clearpage
\begin{Sublemma}\label{L:GADefinition2} If there exist $\delta, z, P,$ and $G$ satisfying $\Phi(\delta,z,P,G)$, then
\begin{enumerate}
\parskip=0pt
\item \label{L:GADefinition2:uniqueness}For each $\gamma$ as in $\Phi$ there is a unique $\M$ witnessing properties (1) through (4) of $\Phi$, which I will denote $\M_{\gamma}$.
\item \label{L:GADefinition2:coherence}The $\M_\gamma$ form a coherent sequence: $\gamma < \gamma^\prime \to \left(\M_{\gamma^\prime}\right)_\gamma = \M_\gamma$. 
\item \label{L:GADefinition2:ZFC}If $\M = \bigcup_{\gamma} \M_{\gamma}$, then $\M  \vDash \ZFC$.
\item \label{L:GADefinition2:forcingext}$V = \M[G]$ is forcing extension of $\M$ by $P$.
\item \label{L:GADefinition2:GAfails}Consequently, the Ground Axiom fails in $V$.
\end{enumerate}
\end{Sublemma}

\begin{proof}
Fix $\gamma$ as in $\Phi$ and suppose $\M$, $\M^\prime$ are transitive structures of height $\gamma$ witnessing the properties stated in $\Phi$.  Note that $(\delta^{+})^{\M} = (\delta^{+})^{\M^\prime} = (\delta^{+})^{V}$, with the first equality holding by the equality of $\mathcal{P}(\delta)$ in $\M$ and $\M^\prime$, and the second equality holding by the fact that $V_\gamma=\M^\prime[G]$ is a forcing extension of $\M^\prime$ by forcing of size $<\!\delta$.  Thus Lemma \ref{L:M=M'} applies, and so $\M=\M^\prime$.  Thus the structure of height $\gamma$ witnessing the properties of $\Phi$ is unique, and will be denoted $\M_\gamma$.

To show the coherence of the $\M_\gamma$ sequence, suppose $\gamma < \gamma^\prime$ and $\M_{\gamma^\prime}$ is a witnessing structure for $\gamma^\prime$.  If $\mathcal{W}=(\M_{\gamma^\prime})_\gamma$, I must show that $\mathcal{W}$ is a witnessing structure for $\gamma$.  That $\mathcal{W}\vDash \ZC$ is straightforward, as $\mathcal{W}$ is a transitive initial segment of a model of $\ZC$ of height a limit ordinal.  The argument used in Lemma \ref{L:GADefinition1} shows that $\mathcal{W}\vDash \leq\!\delta$-Replacement, the pertinent property being $\leq\delta$-Replacement in the larger model $\M_{\gamma^\prime}$.  To see that every set in $\mathcal{W}$ is coded by a set of ordinals, it suffices to show that the transitive closure of every $x\in\mathcal{W}$ has size $<\gamma$ in $\M_{\gamma^\prime}$.  Fix such $x$, and suppose that there is an injection in $\M_{\gamma^\prime}$ from $\gamma$ to $\tc(x)$.  Then this injection is in $V$ as well, contradicting the assumption that $\gamma$ is a $\beth$-fixed point of $V$.  Thus $\mathcal{W}\vDash \ZFC_{\delta}$.  The $\delta$ cover and $\delta$ approximation properties are inherited by $\mathcal{W}\subset V_\gamma$ from $\M_{\gamma^\prime} \subset V_{\gamma^\prime}$ by the same argument given in the proof of Lemma \ref{L:GADefinition1}.  That $z = \mathcal{P}(\delta)^{\M_{\gamma^\prime}} = \mathcal{P}(\delta)^\mathcal{W}$ follows trivially from the definition of $\mathcal{W}$.  Finally, every $x\in V_\gamma$ has a $P$-name in $W_{\gamma^\prime}$ of rank less than $\gamma$, and so $\mathcal{W}[G] = V_\gamma$.  Thus $\mathcal{W}$ is a witnessing structure of height $\gamma$ for $\Phi$, and by uniqueness $\mathcal{W} = (\M_{\gamma^\prime})_\gamma = \M_\gamma$.  This establishes coherency of the $\M_\gamma$ sequence.

Note that $\M = \bigcup_{\gamma} \M_{\gamma}$ is a definable transitive class in $V$.  That $\M$ satisfies Extensionality, Pairing, Union, Power Set, Infinity, Regularity, and Choice can be checked in each instance by simply going to a large enough $\M_\gamma$.  The only difficulty lies in showing that $\M$ satisfies the  Replacement and Separation Schemes.   For Replacement it suffices to show Collection.  Fix $A\in\M$ and $\psi(x,y)$.  As $\M$ is definable in $V$ it follows that satisfaction of $\psi$ in $\M$ can be computed in $V$, that is, for $a,b$ in $\M$ I have $\M \vDash \psi(a,b)\leftrightarrow V \vDash \psi^\M(a,b)$. Since Replacement holds in $V$, there must be an $\alpha$ such that $\forall x\in A \thickspace \left(\exists y\in\M \thickspace \psi^\M(x,y) \to \exists y\in\M_\alpha \thickspace \psi^\M(x,y)\right)$.  Thus $\M \vDash \forall x\in A \thickspace  (\exists y \, \psi(x,y) \to \exists y\in\M_\alpha \, \psi(x,y))$, and so $\M$ satisfies Collection.  Separation is established through the use of the reflection principle in $\M$.  To see that $\M$ satisfies the principle, suppose $\psi(x,b)$ is a first-order formula with parameter $b\in\M$.  I would like to show that there is a class club of $\alpha$ such that $\M_\alpha \vDash \psi(x,b) \leftrightarrow \M \vDash \psi(x,b)$ for every $x \in \M_\alpha$.  I apply reflection in $V$ to the formulas $x\in\M$ and $\psi^\M(x,b)$.  This yields a class club of $\alpha$ such that $V_\alpha$ calculates membership in $\M$ correctly and $V_\alpha \vDash \psi^\M(x,b) \leftrightarrow V \vDash \psi^\M(x,b)$ for every $x \in V_\alpha$. Since $\M^{V_\alpha}=\M_\alpha$, I have $V_\alpha \vDash \psi^\M(x,b) \leftrightarrow \M_\alpha \vDash \psi(x,b)$.  Thus $\M$ satisfies the reflection principle.  To see that $\M$ satisfies Separation, fix $A\in\M$ and $\psi(x)$.  For readability, the parameter $b$ will be suppressed.  Choose $\alpha$ sufficiently large that $A \in \M_\alpha$ and $\M_\alpha$ reflects $\psi(x)$.  Choose $\gamma > \alpha$ a $\beth$-fixed point of cofinality $>\delta$.  Separation holds in $\M_\gamma$, so $B=\{x\in A \mid \M_\gamma \vDash \psi^{\M_\alpha}(x)\}$ is in $\M_\gamma$ and hence in $\M$.  But for $x\in A$, reflection gives $\M_\gamma \vDash \psi^{\M_\alpha}(x) \leftrightarrow \M_\alpha \vDash \psi(x) \leftrightarrow \M \vDash \psi(x)$.  Thus $B=\{x\in A \mid \M\vDash\psi(x)\}$, and so Separation holds in $\M$.  

It remains to show $V=\M[G]$.  Note that $G$ is $\M$-generic for $P$.  Furthermore, for any $A\in V$ there is $\gamma$ above the rank of $A$ such that $A\in \M_\gamma[G]$, and so $A$ is the interpretation by $G$ of some $P$-name $\dot{A}\in\M_\gamma \subset \M$.  Similarly, for any $P$-name $\tau\in\M$ the valuation $\tau_G$ will be in $V$.  Thus $V=\M[G]$ is a forcing extension of $\M$, and so $V$ does not satisfy the Ground Axiom. 
\end{proof}

This concludes the proof of Theorem \ref{T:GADefinition}.  Note that if $V=\M[G]$ is a set-forcing extension, Lemma \ref{L:GADefinition2} gives a first-order definition of $\M$ as a class of $V$ based on the parameters $\delta, \mathcal{P}(\delta)^\M, P,$ and $G$.  By quantifying over $P$ and $G$ we obtain Laver's result \cite{Laver:CertainVeryLargeCardinalsNotCreated}.

\begin{Thm}\label{T:VDefinable}\emph{(Laver)} 
Suppose $V=\M[G]$ is a forcing extension of $\M$ by set forcing $P \in \M$.  Then $\M$ is definable in $V$ from parameters in $\M$.  
\end{Thm}

\section{Forcing the Ground Axiom}\label{S:ForcingGA}

In order to explore the relative consistency of the Ground Axiom with other set-theoretic assertions a general method for producing models of the Ground Axiom is required.  I will begin by identifying an axiom which implies the Ground Axiom.  Informally, this axiom asserts that every set in the universe is coded into the pattern of the $\GCH$ holding and failing at successor cardinals.  An exact statement of the assertion, which I will refer to as the Continuum Coding Axiom, is given below.  After demonstrating that the Ground Axiom is a consequence of the Continuum Coding Axiom, I will discuss a general method by which the property can be forced.   The apparent paradox -- using forcing to produce a model of the Ground Axiom -- is explained by noting that the Ground Axiom refers only to set forcing.  A proper class of forcing will be used to obtain the Continuum Coding Axiom.

\begin{Def}\rm\label{D:ContinuumCodingAxiom}The \emph{Continuum Coding Axiom} ($\CCA$) is the assertion that for every ordinal $\alpha$ and for every $a\subset \alpha$ there is an ordinal $\theta$ such that $\beta \in a \iff 2^{\aleph_{\theta+\beta+1}} = \(\aleph_{\theta+\beta+1})^+$ for every $\beta < \alpha$.
\end{Def}

Note that the $\CCA$ relies on coding only at the successor cardinals, over which we have greater control.  While the $\CCA$ refers only to sets of ordinals, this is equivalent in $\ZFC$ to every set being coded into the continuum function, since under $\AC$ every set is coded as a set of ordinals.  Also note that the $\CCA$ is essentially a strong form of $\VHOD$.

\begin{Thm}\label{T:CCAimpliesGA} 
The $\CCA$ implies the $\GA$.
\end{Thm}

\begin{proof}
Suppose $V$ satisfies the Continuum Coding Axiom.  Suppose further that $V$ is a set-forcing extension of an inner model $V=W[h]$, where $h$ is $W$-generic for some poset $Q\in W$. For  $\kappa > |Q|$, the models $W$ and $V$ will agree on the properties ``$\kappa$ is a cardinal'' and ``the $\GCH$ holds at $\kappa$.''  Every set of ordinals $a$ in $V$ is coded into the continuum function of $V$.  I claim that one such code for $a$ must appear above $|Q|$.  If $|Q|=\aleph_\delta$, consider the set of ordinals $a^\prime=\{\delta+\beta \mid \beta \in a\}$.  As $a^\prime$ is also coded into the continuum function, it is clear that the part of $a^\prime$ above $\delta$ must appear coded into the continuum function above $|Q|$.  Thus $a$ is coded into the continuum function of $V$ above $|Q|$, and so the code appears also in $W$.  Thus $a\in W$, and so every set of ordinals of $V$ is also in $W$.  This shows that $V=W$, and so the forcing $Q$ was trivial.  Thus $V \vDash \GA$. 
\end{proof}

To obtain a general method for forcing the Ground Axiom, I will present a proper class notion of forcing which gives $\CCA$ in the resulting extension.  This folkloric result is based on the work of Kenneth McAloon who presented the basic method in his paper on Ordinal Definability\cite{McAloon1971:ConsistencyResultsAboutOrdinalDefinability}.  See Appendix \ref{S:ClassForcing} for the basics of class forcing.

\begin{Thm}\label{T:CCA}
If $V$ satisfies $\ZFC$, then there is a forcing extension by class forcing which satisfies $\ZFC + \CCA$.
\end{Thm}

\begin{proof}The basic idea is to use forcing to code every set in the universe into the continuum function in the manner of the $\CCA$.  Coding a single set into the continuum function can be accomplished by applying Easton's celebrated result concerning powers of regular cardinals \cite{Easton70:PowersofRegularCardinals}.  Recall that a function $E\from \ORD \rightarrow \ORD$ is an \emph{Easton index function} if $E$ is nondecreasing, $\dom(E)$ is a set of regular cardinals,  and $E(\kappa)$ is a cardinal with $\cof\({E(\kappa)})>\kappa$ for each $\kappa \in \dom(E)$.  Associated with $E$ is a poset $Q(E)$, called Easton forcing, consisting of functions $q$ with domain $\subset \dom(E)$ and, for every regular $\lambda$, the domain of $q$ is bounded in $\lambda$, satisfying for each $\kappa$ in the domain, $q(\kappa)$ is a partial function from $E(\kappa)$ to $2$ of size $<\!\kappa$.  Conditions are ordered by extension on each coordinate and by enlarging the domain, with trivial extensions and enlargements allowed.  Easton's theorem tells us that, under mild $\GCH$ restrictions in the ground model, forcing with $Q(E)$ preserves cardinals and yields a model in which the power set of $\kappa$ has size $E(\kappa)$ for every $\kappa \in \dom(E)$.   What is more, for all $\lambda \notin \dom(E)$ the power set of $\lambda$ is `as small as possible,' i.e. the least cardinal $\alpha \geq \lambda^+$ such that $\cof(\alpha)>\lambda$ and $\alpha\geq E(\kappa)$ for all $\kappa \in \dom(E)\cap\lambda$.  The forcing relies only on a local $\GCH$ assumption, for if $\dom(E)\cup \ran(E)$ is contained in the closed interval $[\delta,\gamma]$ for regular cardinals $\delta$ and $\gamma$, then Easton's Theorem holds if $2^{<\!\delta}=\delta$ and  $2^\kappa = \kappa^+$ for every $\kappa \in [\delta,\gamma)$.  To code a set of ordinals $a\subset \alpha$ into the continuum function starting at  $\aleph_\theta$, define the Easton index function $E_a \from \aleph_{\theta+\alpha} \rightarrow \aleph_{\theta+\alpha}$ such that for all $\beta<\alpha$, the value $E_a \(\aleph_{\theta+\beta+1})$ is $\(\aleph_{\theta+\beta+1})^{+}$ for $\beta\in a$, and $\(\aleph_{\theta+\beta+1})^{++}$ otherwise.  Provided the $\GCH$ hypotheses described above hold on the interval $[\aleph_{\theta+1}, \aleph_{\theta+\alpha+1}]$, forcing with $Q(E_a)$ will code $a$ into the continuum function in exactly the manner of the $\CCA$.  Furthermore, the continuum function will be undisturbed outside this interval.

To prove the Theorem, I will encode every set into the continuum function. This will be accomplished by forcing with an $\ORD$-length iteration which, at stage $\xi$, uses Easton's result to encode the generics obtained at all previous stages.  Note that although the iteration only encodes generics, an easy density argument will show that in fact every set is encoded.  The coding at stage $\xi$ must take place on a ``clean'' interval, lying above the coding performed at all previous stages, to prevent later stages from destroying the coding at earlier stages.  This is accomplished by breaking the cardinals up into intervals, each of sufficient length to encode the generics added on all previous intervals.  The $\xi^{th}$ interval will be given by $[\aleph_{f(\xi)+1},\aleph_{f(\xi+1)})$, where $f(\xi)$ is defined recursively by $f(0)=0$, at successors $f(\xi+1)=\aleph_{f(\xi)}$, and for $\lambda$ a limit $f(\lambda)= \left(\sup_{\alpha<\lambda} f\(\alpha)\right) +1$.

The forcing is defined as follows.  Let $\P=\bigcup P_{\xi}$ be the $\ORD$-length iteration with Easton support, that is, bounded support at inaccessible $\xi$ and full support otherwise, satisfying the following property.  For each $\xi$, the poset $P_{\xi+1}=P_\xi*\dot{Q}_\xi$, where the forcing $\dot{Q}_\xi$ is (a $P_\xi$-name for) the Easton forcing encoding the generic $G_\xi \subset P_\xi$ into the interval $[\aleph_{f(\xi)+1},\aleph_{f(\xi+1)})$, that is, $P_{\xi} \forces \dot{Q}_\xi = Q(E_{\dot{G}_\xi})$.  Note that $\dot{G}_\xi$ is the canonical name for the generic subset of $P_\xi$.  Also observe that $G_\xi$ will not, strictly speaking, be a set of ordinals, but will rather be a sequence of sets $g_\eta$, each generic for the Easton forcing at stage $\eta$.  Each $g_\eta$ is itself a sequence  $g_\eta = \<g_\eta(\kappa)\mid \kappa \in \dom(g_\eta)>$, with each $g_\eta(\kappa)$ a binary sequence of length $\leq \kappa^{++}$.  By concatenating the $g_\eta(\kappa)$, I ``flatten'' $g_\eta$ into a single set $\corners{{g}_\eta} = \{\kappa^{++}+\beta \mid g_\eta(\kappa)(\beta)=1 \mbox{ for }\kappa \in \dom(g_\eta) \mbox{ and } \beta<\kappa\}$.  Note that $\corners{{g}_\eta} \subset \aleph_{f(\eta+1)}$.  The sequence of $\corners{{g}_\eta}$ are themselves concatenated to code $G_\xi$ as a single set of ordinals $\corners{{G}_\xi} = \{\aleph_{f(\eta+1)}+\beta \mid \beta \in \corners{{g}_\eta} \mbox{ for } \eta < \xi \mbox{ and } \beta < \aleph_{f(\eta+1)}\}$.  Here $\corners{{G}_\xi} \subset \aleph_{f(\xi)}$.

Now suppose $G$ is $V$-generic for $\P$.  I will assume that the $\GCH$ holds in $V$, as  if this is not the case, then it can be forced as an initial step of the proof.  That $V[G]$ satisfies $\ZFC$ follows from Easton support together with the increasing closure of the stage $\xi$ forcing as $\xi$ progresses through the ordinals.  For each $\xi$, let $G_\xi$ be the corresponding $V$-generic for $P_\xi$.  I claim that in $V[G]$, I can define $\corners{{G}_\xi}$ as $\{\beta<\aleph_{f(\xi)} \mid 2^{\aleph_{f(\xi)+\beta+1}}=(\aleph_{f(\xi)+\beta+1})^+\}$.  Factor $\P = P_\xi*\dot{Q}_\xi*\P_{tail}$. A standard argument establishes $|P_\xi|\leq \aleph_{f(\xi)+1}$, and it follows that the interval $[\aleph_{f(\xi)+1},\aleph_{f(\xi+1)+1}]$ satisfies the $\GCH$ requirements for Easton's Theorem in $V[G_\xi]$.  Therefore forcing with $(\dot{Q}_\xi)_{G_\xi}$ over $V[G_\xi]$ codes $\corners{{G}_\xi}$ into $[\aleph_{f(\xi)+1}, \aleph_{f(\xi+1)+1}]$ in the manner of the $\CCA$.  The tail forcing is $\leq \aleph_{f(\xi+1)}$-closed, so the coding is preserved in the full extension $V[G]$.  Thus in $V[G]$, every $\corners{{G}_\xi}$ is encoded in the manner of the $\CCA$.

It remains to show that every set of ordinals $a\in V[G]$ is encoded in the same fashion.  Fix $a\subset \alpha$ in $V[G]$ and $\xi$ such $\alpha<(\aleph_{f(\xi)})^+$.   Factor $\P = P_\xi*\dot{Q}_\xi*\P_{tail}$.  The closure of $\dot{Q}_\xi*\P_{tail}$ ensures that $a \in V[G_\xi]$.  Furthermore, the forcing $Q_\xi$ adds a subset to $(\aleph_{f(\xi)})^+$ which is Cohen generic over $V[G_\xi]$.  A standard density argument shows that every bounded subset of $(\aleph_{f(\xi)})^+$ appears as a block in the resulting Cohen generic.  That is, if $g\subset(\aleph_{f(\xi)})^+$ is $V[G_\xi]$-generic for $\Add((\aleph_{f(\xi)})^+,1)$, then there is an ordinal $\theta$ such that $\beta \in a \iff \theta+\beta+1 \in g$ for every $\beta < \alpha$.  The method used to ``flatten'' $g_{\xi}$ to $\corners{{g}_{\xi}}$ and $G_{\xi+1}$ to $\corners{{G}_{\xi+1}}$ preserves each of the generics added on each coordinate of $Q_\xi$ as a contiguous block, and so $a$ appears as a block in $\corners{{G}_{\xi+1}}$.  The forcing at stage $\xi+1$ codes $\corners{{G}_{\xi+1}}\subset \aleph_{f(\xi+1)}$ into the continuum function, and so $a$ is coded into the continuum function as well.  Thus in $V[G_{\xi+1}]$ there is an ordinal $\theta$ such that $\beta \in a \iff 2^{\aleph_{\theta+\beta+1}} = (\aleph_{\theta+\beta+1})^+$ for every $\beta < \alpha$.  Thus $V[G] \vDash \CCA$.
\end{proof}

As the $\GCH$ fails cofinally in $V[G]$, an immediate consequence is

\begin{Cor}\label{C:GAand-GCH} If $V$ satisfies $\ZFC$, then there is a forcing extension by class forcing which satisfies $\ZFC + \GA + \neg \GCH$.  Thus, the answer to Test Question \ref{Ques:GA->GCH} is negative.
\end{Cor}

Furthermore, the method for constructing $V[G]$ can be adapted to provide additional consistency results.  For example, the encoding does not have to begin at $\aleph_0$. The iteration can be easily modified to begin at $\delta$ for any regular $\delta$, which means that the entire iteration will be $<\!\delta$-closed and so no sets will be added to $H_\delta$.  Thus an arbitrary initial segment of the universe can be preserved while forcing the Ground Axiom.

\begin{Thm}\label{T:GAandPreserveVAlpha}If $V \vDash \ZFC$ and $\alpha$ is an ordinal, then there is a forcing extension $V[G]$ of $V$ satisfying the Ground Axiom and having the same initial segment of height $\alpha$, i.e. $V[G]_\alpha = V_\alpha$.
\end{Thm}

\begin{proof}  Working in $V$, fix $\delta > |V_\alpha|$.  If the $\GCH$ does not hold above $\delta$, then as a first step it can be forced to hold there with $<\!\delta$-closed forcing using standard methods.  Then carry out the coding forcing described above, modifying the iteration to begin at $\delta$.  The resulting extension will satisfy the $\CCA$ and thus the Ground Axiom, and as both the $\GCH$ forcing and the coding forcing are $<\delta$-closed no sets will have been added to $V_\alpha$.
\end{proof}

This allows, for example, a measure on $\kappa$ to be preserved by the forcing, because a measure on $\kappa$ is verified in $V_{\kappa+2}$.

\begin{Cor}\label{C:GAand-GCHandMsbl} If $V$ satisfies $\ZFC + \kappa$ is measurable, then there is a forcing extension by class forcing which satisfies $\ZFC + \kappa$ is measurable $+ \GA + \neg \GCH$. 
\end{Cor}

This fact holds not just for `existence of a measurable cardinal' but for any $\Sigma_2$ property.

\begin{Cor}\label{C:GAisSigma2Independent1} If $\phi$ is any $\Sigma_2$ assertion true in $V$, then there is a forcing extension of $V$ by class forcing satisfying $\phi + \GA + \neg \GCH$.
\end{Cor}

\begin{proof}
This is a consequence of the following characterization of $\Sigma_2$ properties. Given a formula $\phi$, the following are equivalent.
\begin{enumerate}
\parskip=0pt
\item \label{C:Sigma2Equivalence1}$\ZFC \proves \phi \iff \sigma$, for some $\Sigma_2$ formula $\sigma$.
\item \label{C:Sigma2Equivalence2}$\ZFC \proves \phi \iff \exists \kappa \; H_\kappa \vDash \psi$, for some first-order formula $\psi$.

\end{enumerate}
It is a straightforward exercise to show that the statement $\exists \kappa \; H_\kappa \vDash \psi$ is provably equivalent to a $\Sigma_2$ formula for any $\psi$, which gives (\ref{C:Sigma2Equivalence2})$\rightarrow$(\ref{C:Sigma2Equivalence1}).  For the other direction, assume (\ref{C:Sigma2Equivalence1}) holds and suppose $\ZFC$ proves $\phi$ is equivalent to $\exists x\forall y\: \phi_0(x,y)$ for some $\phi_0 \in \Delta_0$.  I will show that ``$\omega$ exists${}\wedge\exists x\forall y\: \phi_0(x,y)$'' is a witness for $\psi$ in (\ref{C:Sigma2Equivalence2}).  Fix a model $V$ of $\ZFC$.  If $V\vDash \phi$, then by the reflection principle there is a class club of cardinals $\kappa$ such that $H_\kappa \vDash \omega \mbox{ exists }\wedge\exists x\forall y\: \phi_0(x,y)$.  Now suppose that there is $\kappa_0 > \omega$ such that $H_{\kappa_0} \vDash \exists x\forall y\: \phi_0(x,y)$.  Fix $x_0 \in H_{\kappa_0}$ such that $H_{\kappa_0} \vDash \forall y \phi_0 (x_0,y)$.  Fix any $y_0 \in V$, and let $X$ be a structure containing $y_0$ as an element and $\tc(\{x_0\})$ as a subset such that $|X| < \kappa_0$ and $X$ reflects $\phi_0$.  If $\pi : X \rightarrow M$ is the Mostowski Collapse of $X$, then $\pi$ fixes $x_0$, so $X \vDash \phi_0(x_0,y_0) \iff M \vDash \phi_0(x_0,\pi(y_0))$.  But $H_{\kappa_0}\vDash \phi_0(x_0,\pi(y_0))$ and $M$ is a transitive subset of $H_{\kappa_0}$, so by absoluteness of $\Delta_0$ formulas $M\vDash \phi_0(x_0,\pi(y_0))$.  Thus $X \vDash \phi_0(x_0,y_0)$, and as $X$ reflects $\phi_0$ this shows that $V \vDash \phi_0(x_0,y_0)$.  Since the choice of $y_0$ was arbitrary, it follows that $V \vDash \exists x\forall y\: \phi_0(x,y)$.  This completes the proof of equivalence of (\ref{C:Sigma2Equivalence1}) and (\ref{C:Sigma2Equivalence2}).  The Corollary follows directly. Given a $\Sigma_2$ assertion $\exists x\forall y\: \phi_0(x,y)$ true in $V$, I choose $\kappa$ so that $H_\kappa \vDash \exists x\forall y\: \phi_0(x,y)$.  Then apply Theorem \ref{T:GAandPreserveVAlpha} to force the $\GA$ while preserving $H_\kappa$.  The resulting extension $V[G]$ will satisfy $\exists \kappa \, H_\kappa \vDash \exists x\forall y\: \phi_0(x,y)$, and so $V[G] \vDash  \exists x\forall y\: \phi_0(x,y)$.
\end{proof}

One suprising application of Corollary \ref{C:GAisSigma2Independent1} is the consistency of the Ground Axiom with Martin's Axiom.  Intuitively, Martin's Axiom seems to say ``a lot of forcing has been done,'' whereas the Ground Axiom seems to say the opposite.  In fact, an even stronger result along these lines holds, that the Ground Axiom is relatively consistent with the Proper Forcing Axiom.

\begin{Cor}\label{C:PFA+GA}If $V$ satisfies the Proper Forcing Axiom \emph{($\PFA$)}, then there is a class-forcing extension of $V$ satisfying $\PFA + \GA$.
\end{Cor}

\begin{proof}This follows from the fact that the $\PFA$ is indestructible by $<\aleph_2$-directed closed forcing.  By forcing the $\CCA$ as above but starting the iteration at $\aleph_2$, the iteration will be $<\aleph_2$-directed closed, preserving the $\PFA$ and forcing the $\GA$.
\end{proof}

Theorem \ref{T:GAandPreserveVAlpha} leaves open the question of consistency of the Ground Axiom with supercompact cardinals and other axioms not captured by $\Sigma_2$ formulas.  Just as in the previous corollary, however, the technology of indestructibility can be used to obtain further results.  Since we generally expect the Ground Axiom to hold in the canonical inner models of large cardinals, this theorem fits into the set theoretic program of obtaining large cardinal inner model properties by forcing.

\begin{Thm}\label{T:GAand-GCHandSupercompact} If $V$ satisfies $\ZFC + \kappa$ is supercompact, then there is a forcing extension by class forcing which satisfies $\ZFC + \kappa$ is supercompact $+ \GA + \neg \GCH$.
\end{Thm}

\begin{proof}
Laver's well-known result on indestructibility \cite{Laver78} shows that the supercompactness of $\kappa$ can be made indestructible by $<\!\kappa$-directed closed forcing.   For simplicity in the proof, I will assume that this condition holds in $V$.  Now force $\CCA$ as above with an iteration $\P$ encoding every set into the continuum function, but start the iteration at $\kappa$.  I will argue that $\kappa$ remains supercompact in the resulting extension $V[G]$.

Note that the entire iteration $\P$ is  $<\!\kappa$-directed closed.  However, Laver's result applies only to set forcing so an additional argument is needed.  Fix $\theta > \kappa$.  I will show that the $\theta$-supercompactness of $\kappa$ is preserved in $V[G]$.  Fix $\beta$ large enough so that $\P$ factors as $P_\beta * \P_{tail}$ where the second factor is $\leq \theta^{<\kappa}$-closed.  The first factor $P_\beta$ is $<\!\kappa$-directed closed set forcing, and so $\kappa$ remains $\theta$-supercompact in the partial extension $V[G_\beta]$.  The tail forcing adds no subsets to $\mathcal{P}_\kappa(\theta)$ and so $\kappa$ remains $\theta$-supercompact in the full extension $V[G]$.  Since this is true for every $\theta$, I have $\kappa$ supercompact in $V[G]$.  
\end{proof} 

This method generalizes to other large cardinals for which a comparable indestructibility theorem exists.

The methods for obtaining models of the Ground Axiom described above and in the following section work according to the same basic principle, forcing $\CCA$ or some similar coding axiom.  In particular, the models of the $\GA$ thus produced all satisfy strong versions of $\VHOD$.  It is natural to consider the relationship of $\VHOD$ to the $\GA$.  This relationship arose in Test Question \ref{Ques:GA->HOD}, which asked if the $\GA$ implies $\VHOD$.  The converse question, whether $\VHOD$ implies the $\GA$, is also natural.

In fact, neither implication holds.  The consistency of the $\GA$ with $\nVHOD$ is demonstrated in Chapter \ref{C:GAandVnotHOD}, in which it is shown that every model of $\ZFC$ has a forcing extension satisfying $\GA+\nVHOD$.  The consistency of $\neg\GA$ with $\VHOD$ was demonstrated by McAloon in 1970 \cite{McAloon1971:ConsistencyResultsAboutOrdinalDefinability}.  His result introduced the idea of coding information into the regular cardinals, and using this idea he produced a set-forcing extension $L[G]$ of $L$ in which the generic $G$ is ordinal definable.  Since every element in $L[G]$ is definable from $G$ together with a name from $L$, every element of $L[G]$ is ordinal definable.  As $L[G]$ is an extension by set forcing, clearly $L[G] \vDash \VHOD+\neg\GA$.

Note that the above idea will work for any for any set forcing over $L$ for which the generic is definable in the extension.  For example, Fuchs and Hamkins in \cite{FuchsHamkins:DegreesOfRigidity} describe Suslin trees with the unique branch property.  Forcing with such a tree $T$ adds a single branch to the tree, and so the generic is definable as the unique branch of $T$ in the extension.  Such trees exist in $L$, and forcing over $L$ with the $L$-least such Suslin tree yields an extension in which the generic is definable without parameters.  The extension will satify $\neg\GA + \VHOD$.

In fact, it is not necessary to work in $L$.  These techniques can be combined with the method for forcing the Ground Axiom to work over any model $V$.  

\begin{Thm}\label{T:notGA+V=HOD2}
If $V\vDash \ZFC$, then there is a class-forcing extension $V[G][H]$ satisfying  $\ZFC + \neg\GA + \VHOD$.
\end{Thm}

\begin{proof}
Using the methods of Theorem \ref{T:CCA}, go to a forcing extension $V[G]$ by class forcing in which every set is coded into the continuum function.  Working in $V[G]$, I follow McAloon's strategy described above, doing set forcing to obtain an extension $V[G][H]$ of $V[G]$ in which $H$ is coded into the continuum function.  Since $H$ was added by set forcing, $V[G]$ and $V[G][H]$ agree on the continuum function above the size of the forcing that added $H$.  Thus every set in $V[G]$ remains coded into the continuum function in $V[G][H]$, and so every set in $V[G]$ is ordinal definable in $V[G][H]$.  Since every set in $V[G][H]$ is definable from $H$ together with a name from $V[G]$, it follows that every set is ordinal definable in $V[G][H]$. This completes the proof.  
\end{proof}

\section{The Ground Axiom and the GCH}\label{S:GAandGCH}

The coding used above to force the Ground Axiom is quite flexible, but has the feature that in the resulting model the $\GCH$ fails quite strongly. In this section I will explore a different method of coding which preserves the $\GCH$.  This method relies on the same basic strategy of coding each set into the successor cardinals, differing only in the coding mechanism used.  Rather than coding according to the $\GCH$ holding or failing at the $\beta^{th}$ cardinal, the coding here will be accomplished by controlling whether the $\beta^{th}$ cardinal of the ground model is collapsed in the extension.  Unfortunately, this requires some method of computing in the extension the $\beta^{th}$ cardinal of the ground model.  This will limit the method to work only over models with a certain absoluteness property.

\clearpage
\begin{Def}\rm\label{D:AbsoluteInnerModel}Suppose $U\subset V$ is an inner model which is said to ``exist'' provided a certain first-order formula $\phi$ is satisfied.  Suppose further that $U$ is definable from ordinal parameters.  Then $U$ is \emph{absolutely definable} if for any model $N$ of $\ZFC$ having the same ordinals as $V$, if $N\vDash \phi$, then $U^N = U^V$.  Furthermore, $\phi$ is required to be upward absolute, so if $N\vDash \phi$ and $V\supset N$ has the same ordinals as $N$, then $V \vDash \phi$.  Such a model is \emph{forcing robust} if and only if for any model $N$ of set theory and any forcing extension $N[g]$ by set forcing, $N \vDash \phi \iff N[g] \vDash \phi$.
\end{Def}

Note that the canonical models $L$, $L[0^\#]$, and $L[\mu]$ are all absolutely definable, forcing robust models.

\begin{Thm}\label{T:GAandGCH1}
If $V=U[g]$ is a set-forcing extension of $U$, an absolutely definable, forcing robust model, then there is a forcing extension of $V$ by nontrivial class forcing which satisfies $\GA + \GCH$.
\end{Thm}

\begin{proof}To demonstrate the basic method, I will assume first that the forcing adding $g$ was trivial, that is, that $V=U$ is itself an absolutely definable, forcing robust model, and second that the $\GCH$ holds in $V$.  I will then describe the modifications to the argument necessary to avoid these assumptions.  Suppose that $V$ is an absolutely definable, forcing robust model satisfying the $\GCH$.  I will use an iteration to collapse cardinals of $V$, coding every set into the pattern of ``$(\aleph_\beta)^V$ is a cardinal'' holding or failing in the extension.  This coding differs from that described in the previous section in that the forcing at stage $\beta$ does not code all generics added at previous stages, but rather codes a single `bit' of a single generic added at a previous stage.  Furthermore, to avoid complications arising from collapsing many cardinals in a row, coding will take place not at every successor cardinal but at every other successor cardinal.  Let $\kappa_\beta$ be the $\beta^{th}$ odd cardinal of $V$, that is $\kappa_\beta$ is the $\beta^{th}$ cardinal of the form $\aleph_{\lambda+n}$ where $\lambda$ is a limit ordinal and $n<\omega$ is odd.  Let $\P=\bigcup P_\beta$ be the $\ORD$-length iteration with Easton support, such that for each $\beta$, the poset $P_{\beta+1}=P_\beta*\dot{Q}_\beta$, where the forcing $\dot{Q}_\beta$ at stage $\beta$ is trivial if $\beta \in \corners{{G}_\beta}$ and equals $\Coll(\kappa_\beta,\kappa_\beta^+)$ if $\beta \notin \corners{{G}_\beta}$. The set $\corners{{G}_\beta}$ is obtained by ``flattening'' the generic $G_\beta$ for $P_\beta$ into a single set of ordinals in the following way.  For each $\xi < \beta$, the generic $g_\xi$ added at stage $\xi$ is either empty or is a function $g_\xi \from \kappa_\beta \rightarrow \kappa_\beta^+$.  As such, $g_\xi$ is a collection of pairs of ordinals, and by use of a pairing function $g_\xi$ can be encoded as a single subset $\corners{{g}_\xi}\subset \kappa_\beta^+$ in a canonical way.  For concreteness, I will use the absolute pairing function in which pairs are ordered first by maximum, then by first coordinate, and then by second coordinate.  The entire generic $G_\beta$ for $P_\beta$ is then flattened into a single set of ordinals $\corners{{G}_\beta} = \left\{ \kappa_\xi + \alpha \mid \xi < \beta \mbox{ and } \alpha \in \corners{{g}_\xi}\right\}$.  Informally, a copy of $\corners{{g}_\xi}$ appears in $\corners{{G}_\beta}$ in the interval $[{ \kappa_\xi, \kappa_\xi^+ })$.  Note that for $\beta < \beta^\prime$, the larger generic $\corners{{G}_{\beta^\prime}}$ end-extends the smaller $\corners{{G}_\beta}$.   

Now suppose $G$ is $V$-generic for $\P$.  That $V[G] \vDash \ZFC$ follows from the fact that for any $\beta$ the partial order factors $\P=P_\beta*\P_{tail}$, where $P_\beta \forces \P_{tail}$ is $<\kappa_\beta$-closed.  I claim that the only cardinals collapsed by $\P$ are those of the form $\delta=\kappa_\beta^+$ such that $\beta \notin \corners{{G}_\beta}$.  I will first show that $\kappa_\beta$ remains a cardinal in $V[G]$ for every $\beta$.  Factor $\P=P_\beta * \P_{tail}$.  Inductively, $P_\beta$ has size at most $\kappa_\beta$ and the $\kappa_\beta$-c.c. so it does not collapse $\kappa_\beta$.  The tail forcing is $<\kappa_\beta$-closed and so does not collapse $\kappa_\beta$.  Thus $\kappa_\beta$ remains a cardinal in $V[G]$.  As the $\kappa_\beta$ are cofinal in every limit cardinal, this shows that limit cardinals are preserved in $V[G]$ as well.  Finally, fix $\delta=\kappa_\beta^+$ and factor $\P=P_\beta*\dot{Q}_\beta*\P_{tail}$.  The first factor has at worst the $\kappa_\beta$-c.c., so it cannot collapse $\delta$.  The forcing $\dot{Q}_\beta$ at stage $\beta$ will collapse $\delta$ if and only if $\beta \in \corners{{G}_\beta}$.  Closure of the tail forcing prevents it from collapsing $\delta$, and so $\delta$ is collapsed in the full extension $V[G]$ if and only if $\beta \notin \corners{{G}_\beta}$.  This shows that for any $\beta$, I can define $\corners{{G}_\beta}$ in $V[G]$ as $\{\xi < \kappa_\beta \mid (\kappa_{\xi}^+)^V \mbox{ is a cardinal} \}$.  The definition requires the extension to calculate $(\kappa_{\xi}^+)^V$, which  relies on the absolute definability of $V$.  Note that the code for the generic $\corners{{g}_\xi}$ of stage $\xi$ appears as a block in $\corners{{G}_\beta}$, and so it is coded in a similar fashion into an interval of the cardinals of $V[G]$.  I next show that every set of ordinals $a$ of $V[G]$ is definable from $g_\beta$ for $\beta$ sufficiently large.  Fix $a\subset \gamma$ in $V[G]$ and fix $\xi$ large enough that $\gamma < \kappa_\beta$, so $a \in V[G_\xi]$.  The next nontrivial stage of forcing $\beta > \xi$ will have the form $\Coll(\delta, \delta^+)$ for some $\delta$ and will be $\leq \(\kappa_\beta^+)^V$-closed.  A density argument shows that if $g$ is generic for $\Coll(\delta, \delta^+)$, then there is a $\theta < \delta$ such that $\alpha \in a \iff g(\theta+\alpha)=0$ for all $\alpha<\gamma$.  Thus for sufficiently large $\beta$ the set $a$ is definable from parameters $g_\beta$ and $\theta$.  That $V[G]$ satisfies the Ground Axiom follows directly.  Suppose to the contrary that $V[G] = W[h]$, a forcing extension of $W$ by a poset $Q$.  Forcing robustness combined with absolute definability of $V$ imply that $V\subset W$ and that $W$ computes $V$ correctly.  As $Q$ cannot collapse cardinals greater than $|Q|$, the models $V[G]$ and $W$ agree on the statement ``$(\kappa_\xi^+)^V$ is a cardinal'' for $\xi$ sufficiently large.  Thus for large enough $\beta$ the generic $g_\beta$ will be coded in $W$.  Since every set of ordinals in $V[G]$ is definable from $g_\beta$ for sufficiently large $\beta$, I have $V[G]\subset W$.  Thus $V[G]=W$ and the forcing $Q$ was trivial, so $V[G]$ satisfies the Ground Axiom.

It remains to shows that $V[G]$ satisfies the $\GCH$.  Fix $\beta$ an ordinal and working in $V$ factor $\P=P_\beta *\dot{Q}_\beta * \P_{tail}$.  The first factor has size at most $\kappa_\beta$ and the $\kappa_\beta^+$-c.c., so it cannot add $\kappa_\beta^{++}$ many subsets to $\kappa_\beta$.  The forcing $Q_\beta$ at stage $\beta$ is either trivial or collapses $\kappa_\beta^+$ to $\kappa_\beta$, in either case preserving $2^{\kappa_\beta}=\kappa_\beta^+$.  The tail forcing is $<\kappa_{\beta+1}$-closed, so $2^{\kappa_\beta}=\kappa_\beta^+$ in $V[G]$.  Essentially the same analysis shows that if $\kappa_\beta^+$ remains a cardinal in the extension, then $(2^{\kappa_\beta^+}=\kappa_\beta^{++})^{V[G]}$, noting that the forcing $\dot{Q}_\beta$ is trivial in this case.  The only other cardinals in $V[G]$ are limit cardinals $\delta$ of the form $\delta = \bigcup_{\beta<\lambda} \kappa_\beta$ for some limit ordinal $\lambda$.  In this case, factoring $\P=P_\lambda * \P_{tail}$, the first factor has size $\delta^+$ and the $\delta^+$-c.c., and so it cannot add  $\delta^{++}$ many subsets to $\delta$.  The closure of the tail forcing prevents it from adding subsets to $\delta$, so $2^\delta = \delta^+$ in $V[G]$.  Thus $V[G]$ satisfies the $\GCH$. 

This completes the proof under the additional simplifying assumptions on $V$.  Now suppose $V$ is an absolutely definable, forcing robust model but $V$ does not satisfy the $\GCH$.  It is possible that the coding as described above may fail in this case.  In particular, if $\gamma < \delta$ has power set of size $>\delta^{+}$, then $\Coll(\delta,\delta^+)$ will collapse all cardinals in the interval $[\delta^+, 2^\gamma]$ to $\delta$.  To avoid this issue, I will begin by forcing the $\GCH$ over $V$ in the canonical way, forcing with a proper class iteration that adds a single Cohen subset to each regular cardinal.  The resulting model $\overline{V}$ satisfies the $\GCH$, but of course $\overline{V}$ may no longer be an absolutely definable, forcing robust model.  However, it is completely determined in $V$ whether a given cardinal is collapsed in $\overline{V}$.  That is, if $\Q$ is the canonical forcing of the $\GCH$, then either $1_\Q \forces \kappa$ is a cardinal or $1_\Q \forces \kappa$ is collapsed, for every cardinal $\kappa$.  Thus any model that can compute $V$ can compute ``the $\beta^{th}$ odd cardinal of $\overline{V}$.''  Thus the coding described above can be carried out over $\overline{V}$ to obtain $\overline{V}[G]$, and the proof that $\overline{V}[G]$ satisfies $\GA + \GCH$ goes through as before.

Finally, suppose $V=U[g]$ is a forcing extension of an absolutely definable, forcing robust model $U$ by set forcing $Q \in U$.  As set forcing cannot collapse cardinals above the size of the forcing, it follows that $U$ and $V$ have the same cardinals above $|Q|$.  Begin by forcing the $GCH$ to obtain $\overline{V}$.  Note that $Q$ has no effect on whether the canonical forcing of the $\GCH$ collapses cardinals, at least for cardinals larger than $|Q|$.  Thus, for $\kappa > |Q|$ a regular cardinal, $U$ and $U[g]$ agree on the statement ``the canonical forcing of the $\GCH$ collapses $\kappa$.'' Next, perform the collapsing coding as described above, but begin the iteration after $|Q|$.  The resulting model $\overline{V}[G]$ will be able to correctly calculate the $\beta^{th}$ odd cardinal of $\overline{V}$, at least for cardinals above $|Q|$.  The remainder of the proof follows as above.
\end{proof}

Note that this result applies to any set-forcing extension of $L$, $L[0^\#]$, $L[\mu]$, and many instances of the core model $K$.  The method can be adapted to yield slightly more general results.  Both the coding iteration and the canonical forcing of the $\GCH$ can begin at any regular $\delta$, allowing the preservation of an arbitrary initial segment of the universe.

\begin{Cor}\label{C:GAandGCHoverL} Suppose $U$ is an absolutely definable, forcing robust model.  If $\phi$ is any $\Sigma_2$ assertion forceable over $U$ by set forcing, then there is a forcing extension of $U$ satisfying $\phi + \GA + \GCH$ holds beyond some cardinal $\delta$.
\end{Cor}

\begin{proof} Suppose $\phi$ is a $\Sigma_2$ assertion that holds in $U[g]$, a set-forcing extension of $U$ by $Q$.  Fix $\delta$ so that $H_\delta \vDash \phi$ and $\delta > |Q|$, and force the Ground Axiom using the above coding but beginning both the canonical forcing of the $\GCH$ and the coding iteration at $\delta^+$.  The resulting model $V[G]$ will satisfy $\GA+\GCH$ holds above $\delta$, and the argument given in the proof of Corollary \ref{C:GAisSigma2Independent1} shows that preservation of $H_\delta$ implies that $V[G]$ will also satisfy $\phi$.
\end{proof}

\section{The Bedrock Axiom}\label{S:BedrockAxiom}

What are the models of $\ZFC$ familiar to the working set theorist?  On the one hand, there are canonical models, including $L$, $L[0^\#]$, $L[\mu]$, the core model $K$, and many more, which are  generally characterized by some notion of `minimality,' and in many cases satisfy the Ground Axiom.  On the other hand, many consistency results can be established by performing set forcing over the minimal models, such as the  consistency of $\neg${\sc CH}, of Martin's Axiom, and so on.  These models clearly do not satisfy the Ground Axiom, but they `sit above' a model of the Ground Axiom with only set forcing separating them.  Indeed, a common property of many models of set theory is that they are either models of the Ground Axiom or set-forcing extensions of such models.  Does this hold in general? To investigate this question, I define the \emph{Bedrock Axiom}.  In the usual forcing paradigm one starts in a ground model $V$ and does forcing to obtain an extension $V[G]$.  I would like to shift perspective to the extension. If a model $V$ is a set-forcing extension of some inner model $W$, then I will refer to $W$ as a ground model of $V$.  If a ground model of $V$ satisfies the Ground Axiom, then I will call it a \emph{bedrock model} of $V$, for if one descends through the ground models of $V$ one reaches bedrock when one can descend no further.  The Bedrock Axiom asserts the existence of a bedrock model.

\begin{Def}\rm\label{D:BA}
The \emph{Bedrock Axiom} ($\BA$) asserts there is an inner model $W$ such that $V$ is a set-forcing extension of $W$ and $W\vDash \GA$.
\end{Def}

In particular, any model of the Ground Axiom is a model of the Bedrock Axiom, since every model is trivially a forcing extension of itself.  Once again this apparently second-order statement has a first-order equivalent.

\begin{Thm}\label{T:BAExpressible}
The Bedrock Axiom is first-order expressible.
\end{Thm}

\begin{proof}
The Bedrock Axiom is expressed by the statement ``either $V$ satisfies the Ground Axiom or there exist $\delta, z, P$, and $G$ that satisfy the formula $\Phi$ of Theorem \ref{T:GADefinition} and, in the resulting inner model $W$, there are no $\delta^\prime, z^\prime, P^\prime$, and $G^\prime$ that satisfy $\Phi$ relativized to $W$.''
\end{proof}

Despite the many natural examples of models satisfying the Bedrock Axiom, it is consistent that the axiom fails.

\begin{Thm}\label{T:Con-BA}
There is a forcing extension of $L$ by class forcing which satisfies $\ZFC + \neg \BA$.
\end{Thm}

\begin{proof}
I will start in $L$ and build a class-forcing extension, this time using a product rather than an iteration.  It is the commutative property of products that will be key in showing that the resulting model $L[\class{G}]$ satisfies $\neg \BA$.

For each regular cardinal $\lambda$ of $L$, let $P_\lambda = \Add(\lambda,1)$.  Let $\P=\prod_{\lambda}P_\lambda$ be the canonical Easton product adding a single subset to each regular $\lambda$.  Easton's forcing is explored in detail in the standard texts (such as Jech \cite{Jech:SetTheory3rdEdition} and Kunen \cite{Kunen:Independence}), and I will use the basic results about the forcing without proof.  If $\class{G}$ is $V$-generic for $\P$, then I claim that $L[\class{G}]\vDash \neg \BA$.  I must show that if $L[\class{G}]$ is a forcing extension of some model $W$, then $W$ does not satisfy the Ground Axiom.  As a warm-up, I observe that $L[\class{G}]$ itself does not satisfy the Ground Axiom.  Note that for any regular $\lambda$, I can factor $\P$ as $\P\cong \P^{>\!\lambda}\times P^{\leq\!\lambda}$ where $P^{\leq\lambda}$ has the $\lambda^+$-c.c. and $\P^{>\!\lambda}$ is $\leq\!\lambda$-closed.  Thus $L[\class{G}]=L[\class{G}^{>\lambda}][G^{\leq\lambda}]$, a nontrivial set-forcing extension of $L[\class{G}^{>\lambda}]$, and so $L[\class{G}]$ does not satisfy the Ground Axiom. 

Now suppose $L[\class{G}] = W[h]$ where $h$ is $W$-generic for some poset $Q\in W$.  Let $\delta = \left(|Q|^{+}\right)^W$.  Then $\delta$ is a regular cardinal of $L$, and $\P\cong \P^{>\delta}\times P^{\leq\delta}$.  I claim that $L[\class{G}^{>\delta}]\subset W$.  Clearly $L\subset W$, so it suffices to show that $g_\lambda\in W$ for regular $\lambda>\delta$.  Fix such $\lambda$ and consider $g_\lambda\subset \lambda$.  Since $g_\lambda$ is $L$-generic for $\Add(\lambda, 1)$, every initial segment of $g_\lambda$ is in $L$ and thus in $W$.  As shown in the proof of Lemma \ref{L:GADefinition1}, if $\lambda$ is a regular cardinal larger than $|Q|$, then $W \subset W[h]$ satisfies the $\lambda$ cover and $\lambda$ approximation properties.  I will show that $W$ contains every $\lambda$ approximation of $g_\lambda$, and thus contains $g_\lambda$.  Fix $C\in W$ of size $<\lambda$.  Then $C\cap\lambda$ is bounded by some $\beta<\lambda$, and so $C\cap g_\lambda = C\cap(\beta\cap g_\lambda)$.  Since every initial segment of $g_\lambda$ is already in $L\subset W$, it follows that $C\cap g_\lambda\in W$ as well.  Thus $g_\lambda\in W$ by the $\lambda$ approximation property, so $L[\class{G}^{>\delta}]\subset W$.  

Thus $L[\class{G}^{>\delta}]\subset W\subset L[\class{G}]$, so $W$ is an intermediate model between  $L[\class{G}^{>\delta}]$ and a forcing extension $L[\class{G}^{>\delta}][G^{\leq\delta}]=L[\class{G}]$ by set forcing $P^{\leq\delta}$.  It follows that $W$ is a possibly trivial forcing extension of $L[\class{G}^{>\delta}]$, a result which appears in many places in the literature. For example in Jech \cite{Jech:SetTheory3rdEdition} p.265 it is stated as follows:

\begin{Sublemma}\label{L:FactoringBAs}
Let $G$ be $V$-generic for a complete Boolean algebra $B$.  If $N$ is a model of $\ZFC$ such that $V\subseteq N\subseteq V[G]$, then there exists a complete subalgebra $D\subseteq B$ such that $N=V[D\cap G]$. 
\end{Sublemma}

Continuing with the proof of Theorem \ref{T:Con-BA}, if $W\ne L[\class{G}^{>\delta}]$, then $W$ is a set-forcing extension of $L[\class{G}^{>\delta}]$, and so does not satisfy the Ground Axiom.  If $W = L[\class{G}^{>\delta}]$, then observe that for any regular $\gamma>\delta$, I can factor $\P^{>\delta}\cong\P^{>\gamma}\times P^{(\delta,\gamma]}$, and $L[\class{G}^{>\delta}]=L[\class{G}^{>\gamma}][G^{(\delta,\gamma]}]$.  In this case, $W$ is a set-forcing extension of $L[\class{G}^{>\gamma}]$ by set forcing $P^{(\delta,\gamma]}$ and so does not satisfy the Ground Axiom.  This shows that no ground model $W$ of $L[\class{G}]$ satisfies the Ground Axiom, and so $L[\class{G}] \vDash \neg \BA$. 
\end{proof}

Note that in the above construction forcing need not occur at every regular $\lambda$ but may be restricted to a definable proper subclass.  The same result may be obtained by starting the forcing above any fixed $\kappa$, or by forcing at every other regular cardinal, etc.  This allows the technique to be combined with that used in Section \ref{S:ForcingGA} to provide results of much greater generality.

\begin{Thm}\label{T:-BA1}
If $V$ satisfies $\ZFC$, then there is a forcing extension by class forcing which satisfies $\ZFC + \neg \BA$.
\end{Thm}

\begin{proof}
Begin by following the strategy described in Section \ref{S:ForcingGA} to obtain a model $V[\class{G}]\vDash \CCA$ in which every set is definable from the continuum function.  Now proceed as in the proof above, forcing over $V[\class{G}]$ with a class product to add a single subset to regular cardinals $\lambda$.  However, in order to ensure preservation of cardinals and of the continuum function I must force only at those regular $\lambda$ for which $2^{<\lambda}=\lambda$.  If $\class{H}$ is $V[\class{G}]$-generic for this product, then the resulting model $V[\class{G}][\class{H}]$ has the same cardinals and the same continuum function as $V[\class{G}]$. Thus in $V[\class{G}][\class{H}]$ every set in $V[\class{G}]$ is definable from the continuum function.  Now suppose $V[\class{G}][\class{H}]=W[h]$ where $h$ is $W$-generic for some poset $Q\in W$.  Once again, $W$ and $V[\class{G}][\class{H}]$ must agree on the value of $2^{\alpha}$ for $\alpha$ sufficiently large, as well as on the statement ``$\alpha$ is the $\beta^{th}$ cardinal.''  Thus every set in $V[\class{G}]$ is definable from the continuum function in $W$, and so $V[\class{G}] \subset W$.  Now set $\delta = (|Q|^{+})^W$.  For $\lambda > \delta$, I will show $h_\lambda \in W$.  As $V[\class{G}] \subset W$, every initial segment and hence every $\lambda$ approximation of $h_\lambda$ is in $W$.  Since $\<W,W[h]>$ satisfies the $\lambda$ approximation property, $h_\lambda \in W$.  Thus $V[\class{G}][\class{H}^{>\delta}]\subset W$.  I now have $V[\class{G}][\class{H}^{>\delta}] \subseteq W \subset V[\class{G}][\class{H}]$ and so, by factoring $\class{H}^{>\delta}$ if necessary, $W$ is a forcing extension of an inner model.  Thus $V[\class{G}][\class{H}] \vDash \neg \BA$.   
\end{proof}

Various modifications provide further results.  An arbitrary initial segment of $V$ can be preserved when constructing $V[\class{G}][\class{H}]$ by restricting both the class iteration and the class product to stages above some suitably chosen $\kappa$.  This allows the preservation of a measurable cardinal, as in Corollary \ref{C:GAand-GCHandMsbl}, and in combination with indestructibility allows preservation of a supercompact cardinal, as in Theorem \ref{T:GAand-GCHandSupercompact}, or of any other large cardinal for which a comparable indestructibility result exists.

\clearpage
\begin{Cor}\label{C:-BAandMsbl} Suppose $V$ satisfies $\ZFC$.
\begin{enumerate}
\parskip=0pt
\item If $V$ satisfies $\kappa$ is measurable, then there is a forcing extension by class forcing which satisfies $\ZFC + \kappa$ is measurable  $+ \, \neg\BA$.
\item If $V$ satisfies $\kappa$ is supercompact, then there is a forcing extension by class forcing which satisfies $\ZFC + \kappa$ is supercompact $+ \, \neg\BA$. 
\end{enumerate}
\end{Cor}

In contrast to the models of the Ground Axiom, none of the models of $\neg\BA$ produced above satisfy $\VHOD$.   This is because the Easton product forcing above is almost homogeneous, and such extensions, if nontrivial, always produce models of $\nVHOD$.  A proof of this fact appears in Chapter \ref{C:GAandVnotHOD}, in the proof of Theorem \ref{Theorem.GA+VneqHODoverL}.  However, by combining the product forcing above with a set version of the coding used in Theorem \ref{T:CCA} to force $\CCA$,  a model of $\VHOD + \neg\BA$ can be obtained.

\begin{Thm}\label{T:-BA+VHOD}There is a forcing extension of $L$ by class forcing satisfying $\VHOD + \neg\BA$.
\end{Thm}

\begin{proof}I will force with a class product $\P=\prod_{\lambda} P_\lambda$ with Easton support, that is, bounded support at inaccessibles and full support otherwise.  Rather than each factor simply adding a set to a cardinal, the factor $P_\lambda$ will be a short iteration coding its own generic into the continuum function on a particular interval of cardinals.  To keep track of the intervals, I will once again use the function $f$ defined in Theorem \ref{T:CCA}, with $f(0)=0$, at successors $f(\xi+1)=\aleph_{f(\xi)}$, and for $\lambda$ a limit $f(\lambda)= \left(\sup_{\alpha<\lambda} f\(\alpha)\right) +1$.  Forcing will occur only at limit ordinals $\lambda$, and for such a $\lambda$ the forcing $P_\lambda$ is defined as $P_\lambda = \bigcup_{n\in\omega} P_{\lambda,n}$ a forcing iteration of length $\omega$ with full support.  For each $n<\omega$, the forcing $P_{\lambda,n+1}=P_{\lambda,n}*\dot{Q}_{\lambda,n}$ such that $\dot{Q}_{\lambda,n}$ is (a $P_{\lambda,n}$-name for) the Easton forcing coding the generic $G_{\lambda,n} \subset P_{\lambda,n}$ into the interval $[\aleph_{f(\lambda+n)+1},\aleph_{f(\lambda+n+1)})$.  Note that each generic $G_{\lambda,n}$ is ``flattened'' into a single set of ordinals $\corners{{G}_{\lambda,n}}$ before being encoded, just as in Theorem \ref{T:CCA}.  The generic $G_\lambda$ for $P_\lambda$ is completely determined by the collection of $\{G_{\lambda,n} \mid n\in \omega\}$, and forcing with $P_\lambda$ codes all $G_{\lambda,n}$ into the continuum function in the interval $[\aleph_{f(\lambda)+1},\aleph_{f(\lambda+\omega)})$ in the manner of the $\CCA$.  Furthermore, a standard analysis shows that $P_\lambda$ preserves cardinals and leaves the continuum function undisturbed outside of this interval.  The size and chain condition of $P_\lambda$ are both $\aleph_{f(\lambda+\omega)}$. For any limit $\lambda$, the product $\P$ factors as $\P \cong P^{\leq\lambda}\times \P^{>\lambda}$, where $P^{\leq\lambda}$ consists of all conditions $p\in\P$ such that $\dom(p)\subset\lambda+1$, and $\P^{>\lambda}$ consists of those conditions for which $\dom(p)\cap\lambda+1 = \emptyset$.  A straightforward calculation shows $P^{\leq\lambda}$ has size $\aleph_{f(\lambda+\omega)}$ and the $\aleph_{f(\lambda+\omega)}$-c.c., and $P^{>\lambda}$ is $\leq \aleph_{f(\lambda+\omega)}$-closed.  It is a standard fact from the theory of products that for $\kappa$ a regular cardinal, if $\P$ factors at $\kappa$, that is, $\P\cong P^1\times \P^2$ where $P^1$ has the $\kappa^+$-c.c. and $\P^2$ is $\leq\kappa$-closed, then every function $f\from \kappa \rightarrow L$ in the extension by $\P$ is already in the extension by $P^1$ (Jech \cite{Jech:SetTheory3rdEdition} p.234).  Recalling that $\aleph_{f(\lambda+\omega)}$ is a successor cardinal and therefore regular, this gives arbitrarily large regular $\kappa$ such that $\P$ factors at $\kappa$.  This is sufficient to ensure that forcing with $\P$ preserves $\ZFC$.

Suppose $\class{G}$ is $L$-generic for $\P$.  I will first show that $\P$ preserves cardinals.  Fix $\kappa$ a regular cardinal of $L$.  If there exists a limit ordinal $\lambda$ such that $\aleph_{f(\lambda)} \leq \kappa < \aleph_{f(\lambda+\omega)}$, then factor $\P \cong P^{\leq\lambda}\times \P^{>\lambda}$.  If $\kappa$ is collapsed in $L[\class{G}]$, then it must already be collapsed in $L[G^{\leq\lambda}]$.  I now factor $P^{\leq\lambda} \cong P_\lambda \times P^{<\lambda}$, and so $L[G^{\leq\lambda}] = L[G_\lambda][G^{<\lambda}]$.  The forcing $P^{<\lambda}$ has size $\aleph_{f(\lambda)}$ and the $\aleph_{f(\lambda)}$-c.c. in $L$, and the closure of $P_\lambda$ means that this remains true in $L[G_\lambda]$.  Thus $P^{<\lambda}$ will not collapse $\kappa$ over $L[G_\lambda]$, so if $\kappa$ is collapsed it must already be collapsed in $L[G_\lambda]$.  This is impossible as $P_\lambda$ preserves cardinals, so $\kappa$ must remain a cardinal in $L[\class{G}]$.  If there is no limit $\lambda$ as described above, then $\kappa$ must be a limit of $\aleph_{f(\beta)}$ for $\beta<\lambda$.  If $\kappa$ is collapsed, then some successor $\delta<\kappa$ is also collapsed, which is impossible.  Thus $\P$ preserves all cardinals.

I claim that for every limit ordinal $\lambda$ and $n<\omega$, the code for the generic $\corners{{G}_{\lambda,n}}$ is definable in $L[\class{G}]$ as $\corners{{G}_{\lambda,n}} = \{ \beta < \aleph_{f(\lambda+n)} \mid  2^{\aleph_{f(\lambda+n)+\beta+1}} = ( \aleph_{f(\lambda+n)+\beta+1} )^+ \}$.  Fix $\lambda$ a limit and factor $\P \cong P^{\leq\lambda}\times \P^{>\lambda}$.  The second factor $\P^{>\lambda}$ is $<\aleph_{f(\lambda+\omega)}$-closed and so it suffices to show that the claim holds in $L[G^{\leq\lambda}]$.  I once again factor $P^{\leq\lambda} \cong P_\lambda \times P^{<\lambda}$, and as $P^{<\lambda}$ is too small to affect the $\GCH$ above $\aleph_{f(\lambda)+1}$ it suffices to show the claim holds in $L[G_\lambda]$.  However, the forcing $P_\lambda$ codes $\corners{{G}_{\lambda,n}}$ into the continuum function in exactly the way described in the claim.  Thus the claim holds in $L[\class{G}]$.  This shows that in $L[\class{G}]$ every $\corners{{G}_{\lambda,n}}$, and therefore every $G_{\lambda,n}$, is definable from parameters $\lambda$ and $n$.  As the definition is uniform, it follows that $G^{\leq\lambda}$ is ordinal definable in $L[\class{G}]$ for every $\lambda$.  Any set $a\in L[\class{G}]$ is in $L[G^{\leq\lambda}]$ for some $\lambda$, and so $a$ is definable from $G^{\leq\lambda}$ together with a name $\dot{a}\in L$ for $a$.  As every member of $L$ is ordinal definable in $L[\class{G}]$, $a$ is ordinal definable there as well.  Thus $L[\class{G}] \vDash \VHOD$.

It remains to show that $L[\class{G}]$ satisfies $\neg\BA$.  That $L[\class{G}]$ itself does not satisfy the Ground Axiom follows from the observation that for any $\lambda$ the model $L[\class{G}]$ can be written $L[\class{G}]=L[\class{G}^{>\lambda}][G^{\leq\lambda}]$, where the latter extension is by the set forcing $P^{\leq\lambda}$.  Now suppose $L[\class{G}]=W[h]$, a forcing extension of $W$ by a poset $Q$.  As $Q$ cannot affect the continuum function above $|Q|$, the models $W$ and $L[\class{G}]$ agree on the assertion $2^\delta = \delta^+$ for $\delta>|Q|$.  Thus for $\lambda$ such that $\aleph_{f(\lambda)} > |Q|$, the generic $G_\lambda$ is definable in $W$.  Fixing such a $\lambda$, I have $L[\class{G}^{>\lambda}]\subset W$.  The argument given in the proof of Theorem \ref{T:CCA} shows that $W$ is therefore a set-forcing extension of $L[\class{G}^{>\lambda^\prime}]$ for $\lambda^\prime>\lambda$.  Thus $W$ does not satisfy the Ground Axiom, and so $L[\class{G}]\vDash \neg\BA$.
\end{proof}

\chapter[The Ground Axiom and $\nVHOD$]{The Ground Axiom is consistent with $\nVHOD$}\label{C:GAandVnotHOD}

In this chapter, I will describe a method for forcing $\GA + \nVHOD$ over $L$.  I will then show how to adapt the method to work over any $ZFC$ model.  Furthermore, the method is flexible and simple variations provide consistency of $\GA + \nVHOD$ with various large cardinals.  The material in this chapter was developed jointly by Hamkins, Woodin, and myself \cite{HamkinsReitzWoodin:TheGroundAxiomandVnotHOD}.

\section{Forcing the Ground Axiom and $\nVHOD$}

Prior methods of forcing \GA\ use an $\ORD$-length iteration adding subsets to regular cardinals.  Whether or not sets are added to a given cardinal is determined according to a carefully defined coding scheme.  However, over $L$ this coding is in fact unnecessary.  A simple iteration adding a single subset to each regular cardinal will also force \GA.

\begin{Thm} There is a proper-class-forcing extension of $L$ which satisfies $\GA + \nVHOD$. \label{Theorem.GA+VneqHODoverL}\end{Thm}

A preliminary result will be necessary before proceeding.  A key step in the proof of the Theorem will rely on the result of Laver given in Chapter \ref{C:TheGroundAxiom}, Lemma \ref{L:M=M'}, which states that the $\delta$ cover and $\delta$ approximation properties (Definition \ref{D:deltacover&approx}) together with $\mathcal{P}(\delta)$ uniquely determine an inner model.  However, in order to the apply Laver's result I will need the following Lemma of Hamkins \cite{Hamkins2003:ExtensionsWithApproximationAndCoverProperties}, which gives certain conditions under which a forcing extension will satisfy the $\delta$ cover and $\delta$ approximation properties.

\begin{Sublemma}\label{L:closurepointforcing}\emph{ (Hamkins)}  
Fix $W$ a model of $\ZFC$ and $\delta$ a regular cardinal of $W$.  Suppose $G*H$ is $W$-generic for a partial order $P*\dot{Q} \in W$ such that the extension $W[G][H]\vDash \ZFC$, and satisfying $P$ is nontrivial, $|P|<\delta$ and $P \Vdash \dot{Q}$ is $<\!\delta$-strategically closed.  Then $\<W, W[G][H]>$ satisfies the $\delta$ cover and $\delta$ approximation properties.
\end{Sublemma}

\begin{proof}  Note first that the forcing $\dot{Q}$ is allowed to be trivial, which shows that any set forcing $P$ will satisfy the $\delta$ cover and $\delta$ approximation properties for $\delta > |P|$.  On the other hand, $\dot{Q}$ is allowed to be a proper class, definable in $W$ with parameters, which is the case that will be used in the proof of Theorem \ref{Theorem.GA+VneqHODoverL}. 

That $\<W, W[G]>$ satisfies the $\delta$ cover property is a well-known result for any forcing with the $\gamma$ chain condition for some  $\gamma<\delta$.  The strategic closure of $\dot{Q}$ ensures that it adds no new subsets of $W$ of size $<\!\delta$, and so the full extension $\<W, W[G][H]>$ also satisfies the $\delta$ cover property.  This proof of the $\delta$ approximation property follows an idea due to Mitchell \cite{Mitchell2003:ANoteOnHamkinsApproximationLemma}.  Note that it suffices to show that the property holds for sets of ordinals.  Suppose that $A\subset \theta$ is in $W[G] \backslash W$ and every $\delta$ approximation of $A$ is in $W$.  Formally, the latter statement is the assertion $\forall C \in W \, \({|C|^W < \delta \rightarrow C\cap A \in W})$.  Fix a name $\dot{A}$ for $A$ such that $\Vdash_{P*\dot{Q}} ``\dot{A}\subset \theta \wedge{}$ every $\delta$ approximation of $\dot{A}$ is in $W$.''  To see that the latter property is first-order expressible, note that it is equivalent to consider only those sets $C\in\mathcal{P}(\theta)^W$.   Now suppose there is no condition $\<p,q>$ forcing $\dot{A} \in W$.  Observe that for every $\<p,q>$ there must be an $\alpha \in \theta$ and strengthenings $\<p^0, q^0>$ and $\<p^1, q^1>$ below $\<p,q>$  such that $\<p^0, q^0> \Vdash \alpha\in \dot{A}$ and $\<p^1, q^1> \Vdash \alpha\notin \dot{A}$.  I can assume that $p^0$ and $p^1$ are incompatible, and by mixing I can obtain a single condition $q^\prime$ such that $\<p^0, q^\prime> \Vdash \alpha\in \dot{A}$ and $\<p^1, q^\prime> \Vdash \alpha\notin \dot{A}$.  Without loss of generality I can assume that $\Vdash_P q^\prime < q$, for if necessary I extend $\{p^0,p^1\}$ to a maximal antichain of $P$ such that for every $p^* \in P$ there is $q^*$ such that $p^* \forces q^* < q$, and then apply mixing to obtain $q^\prime$.  As $|P|=\gamma<\delta$, there is an enumeration $\{p_\beta \,|\, \beta < \gamma\}$ of $P$.  Working in $W$, fix a $P$-name $\dot{\sigma}$ for a strategy for $\dot{Q}$.  Use $\dot{\sigma}$ to construct $\gamma$ sequences $\<p^0_\beta>$ and $\<p^1_\beta>$ in $P$,  a sequence $\<q_\beta>$ in $\dot{Q}$, and a sequence of ordinals $\<\alpha_\beta>$ such that
\begin{enumerate}
\parskip = 0pt
\item $\<p^0_\beta, q_\beta> \Vdash \check{\alpha}_\beta \in \dot{A}$ and $\<p^1_\beta, q_\beta> \Vdash \check{\alpha}_\beta \notin \dot{A}$ for every $\beta < \gamma$, and
\item $\Vdash_P \<q_\beta \,|\, \beta < \check{\gamma}>$ is strictly descending in $\dot{Q}$.
\end{enumerate}
Suppose the sequence has been defined below $\beta$.  I will use $\dot{\sigma}$ to obtain a condition $q\in \dot{Q}$ that is forced to lie below the sequence $\<q_{\beta^\prime}  \,|\, \beta^\prime < \beta >$.  Then apply the observation to the condition $\<p_\beta, q>$ to obtain $p^0_\beta$ and $p^1_\beta$ below $p_\beta$ with associated $q_\beta$ and $\alpha_\beta$ satisfying the two conditions above.  Thus the construction can be carried out for every $\beta < \gamma$.  Furthermore, I can apply $\dot{\sigma}$ one more time to obtain $q_\gamma$ which is forced to lie below the entire sequence $\<q_\beta>$.  Now consider the set $C = \{\alpha_\beta \,|\, \beta < \gamma\}$.  This is a set in $W$ of size $<\!\delta$, and so by assumption its intersection with $\dot{A}$ is forced to be in $W$.  That is, there is a set $B$ in $W$ and a condition $\<p,q> < \<1,q_\gamma>$ such that $\<p,q> \Vdash \dot{A} \cap \check{C} = \check{B}$.  In particular, for every $\alpha \in C$, the condition $\<p,q>$ decides the statement $\check{\alpha} \in \dot{A}$.  However, $p=p_\beta$ for some $\beta$, and so $\<p,q> < \<p_\beta, q_\beta>$.  But by extending the first coordinate I obtain $\<p^0_\beta, q> < \<p^0_\beta, q_\beta>$ that forces $\alpha_\beta \in \dot{A}$ and $\<p^1_\beta, q> < \<p^1_\beta, q_\beta>$ that forces $\alpha_\beta \notin \dot{A}$.  Thus $\<p,q>$ does not decide the value of $\check{\alpha_\beta} \in \dot{A}$, a contradiction.  Therefore there must be a condition forcing $A\in W$.  The argument is readily adapted to show that such conditions are dense in $P*\dot{Q}$, and so $\<W,W[G]>$ satisfies the $\delta$ approximation property. \end{proof}

I now show that there is a class-forcing extension of $L$ satisfying $\GA + \nVHOD$.

\begin{proof}[Proof of Theorem \ref{Theorem.GA+VneqHODoverL}] Working in $L$, let $\P$ be the Easton support iteration forcing with $\Add(\gamma,1)$ at every regular cardinal stage $\gamma$.  Suppose $\class{G}$ is $L$-generic for $\P$.  I will first establish that $L[\class{G}] \vDash \nVHOD$.  It is known that the forcing to add a subset to a cardinal is almost homogeneous, and it follows that the full iteration is as well.  A poset is \emph{almost homogeneous} if for any  $p$ and $q$ there is a automorphism of the poset sending $p$ to a condition compatible with $q$.  Consider a set $x$ of ordinals added by the forcing.  There must be an $\alpha$ and incompatible conditions $p$ and $q$, one forcing that $\alpha \in x$ and the other forcing the opposite.  If the forcing is almost homogeneous, then there is an automorphism sending $p$ to something compatible with $q$.  If $x$ is defined by a formula with ordinal parameters in the extension, then applying the automorphism will not affect the value of the formula, as automorphisms fix check names.  However, the automorphism will change the value of the assertion that $\alpha$ is in the set defined by the formula, a contradiction.  Thus $L[\class{G}]$ satisfies $\nVHOD$.  

I will next show that $L[\class{G}]$ satisfies \GA.  Suppose to the contrary that $\GA$ fails.  Then $L[\class{G}]=W[h]$ is a forcing extension of an inner model $W$, where $h$ is $W$-generic for some $Q \in W$.  Fix $\kappa$ a singular strong limit cardinal of $W[h]$ with $\kappa > |Q|$.  I can factor $\P$ as $\P = P_\kappa * \P_{tail}$, where $P_\kappa$ is the forcing up to stage $\kappa$ and $P_\kappa \forces \P_{tail}$ is $\leq\kappa$-closed.  The extension $L[\class{G}]$ can be rewritten $L[\class{G}] = L[G_\kappa][\class{G}_{tail}]$.  The generic $G_\kappa$ can be considered a subset of $\kappa$.  As $G_\kappa \in W[h]$, it has a $Q$-name $\dot{G_\kappa}$ in $W$ such that $(\dot{G_\kappa})_h = G_\kappa$.  Since $|Q|<\kappa$, I can ensure that $\dot{G_\kappa}$ is coded as a subset of $\kappa$ in $W$.  Furthermore, since $\kappa$ is a strong limit I can code all bounded subsets of $\kappa$ from $W$ as a single subset of $\kappa$, for there are $\kappa$ many such subsets.  Thus there is a set $A \subset \kappa$ in $W$, such that $A$ codes both $\dot{G_\kappa}$ and $(2^{<\kappa})^W$.  Now consider the model $L[A]$.  First note that $(2^{<\kappa})^{L[A]} = (2^{<\kappa})^W$.  The forward inclusion holds since $L \subset W$ and $A\in W$.  For the reverse inclusion, observe that every bounded subset of $\kappa$ from $W$ appears coded in $A$ and thus is in $L[A]$.  This means that $Q \in L[A]$ and that $h$ is $L[A]$-generic for $Q$, so I can consider the forcing extension $L[A][h]$.  I next show that $L[A][h] = L[G_\kappa]$.  For the forward inclusion, note that $A$ and $h$ are subsets of $\kappa$ and both are in full extension $L[G_\kappa][\class{G}_{tail}]$.  As the forcing $\P_{tail}$ is $\leq\kappa$-closed, $A$ and $h$ must already be in $L[G_\kappa]$.  For the reverse inclusion, recall that $\dot{G_\kappa}$ appears in $A$, and so $(\dot{G_\kappa})_h = G_\kappa$ is in $L[A][h]$.  This means that $\class{G}_{tail}$ is $L[A][h]$-generic for $\P_{tail}$, and so $L[\class{G}]$ is a forcing extension of $L[A]$ by $Q*\P_{tail}$, that is, $L[\class{G}] = L[A][h][\class{G}_{tail}]$.

It remains to show that $W=L[A]$.  I will argue that $W\subset L[G]$ and $L[A] \subset L[G]$ satisfy the hypotheses of Lemma \ref{L:M=M'} and thus $W=L[A]$.  Fix any successor $\delta$ larger than $|Q|^+$ and smaller than $\kappa$.  Consider $L[A][h][\class{G}_{tail}]$ as a forcing extension of $L[A]$ by $Q*\P_{tail}$.  Lemma \ref{L:closurepointforcing} tells us that the extension $L[A] \subset L[A][h][\class{G}_{tail}]$ satisfies the $\delta$ cover and approximation properties.  The extension $W \subset W[h]$ also satisfies $\delta$ cover and approximation by the same argument, noting that in Lemma \ref{L:closurepointforcing} the second stage of forcing may be trivial.  In both cases the larger model is the same, $L[A][h][\class{G}_{tail}]=W[h]=L[\class{G}]$.  As $W$ and $L[A]$ have the same bounded subsets of $\kappa$ I have $\mathcal{P}({\delta})^W=\mathcal{P}({\delta})^{L[A]}$.  Finally, the three models $W$, $L[A]$, and $L[\class{G}]$ agree on $\delta^+$.  Thus by Lemma \ref{L:M=M'}, I have $W=L[A]$.  This implies that $W[h]=L[A][h]=L[G_\kappa]$, contradicting the assumption that $W[h]=L[G_\kappa][\class{G}_{tail}]$.\end{proof}

The above method has some flexibility.  In particular, it is not necessary to add a subset to every regular cardinal.  The same proof will work when forcing in the iteration is restricted to any absolutely definable subclass of the cardinals provided it is cofinal in $\ORD$, for example the odd cardinals or the successors of limit cardinals.  This also means that the first nontrivial stage of forcing can be as large as I like, allowing the preservation of an arbitrary initial segment of the universe.  A consequence is that the theorem holds not just for $L$ but for any model of the form $L[g]$, a set-forcing extension of $L$.

Extending the theorem to all models of $\ZFC$ requires a little more work.  The property of $L$ essential to the proof is that, given $L[G]=W[h]$, I can always conclude $L\subset W$.  In order to obtain the analogous property for $V$, I prepare the universe by forcing the $\CCA$.  After this preparatory forcing, a similar reverse Easton iteration will once again force $\GA + \nVHOD$.

\begin{Thm} If $V$ is a model of $ZFC$, then there is a proper-class-forcing extension of $V$ which satisfies $\GA + \nVHOD$. \label{Theorem.GA+VneqHODoverV}\end{Thm}

\begin{proof}  I prepare the universe by forcing the $\CCA$ as in the proof of Theorem \ref{T:CCA}.  For simplicity of notation, I assume that this preparatory forcing has taken place and that $V$ therefore satisfies the $\CCA$.  Let $\P$ be the reverse Easton iteration adding a single subset to each regular cardinal $\kappa$ satisfying $2^{<\kappa} = \kappa$.    It follows from the definition of the preparatory forcing that the collection of such $\kappa$ is cofinal in $\ORD$, so forcing will occur at arbitrarily large stages.  Suppose $\class{G}$ is $V$-generic for $\P$.  Once again $\P$ is almost homogeneous and so $V[\class{G}]$ satisfies $\nVHOD$.  It follows from Easton support together with the fact that forcing only occurs at $\kappa$ satisfying $2^{<\kappa} = \kappa$ that $\P$ does not collapse cardinals, and so $V[\class{G}]$ has the same cardinals and the same continuum function as $V$.  Thus every set in $V$ is coded into the continuum function in $V[\class{G}]$.  I now show that $V[\class{G}]$ satisfies $\GA$.  Suppose to the contrary that $V[\class{G}]=W[h]$, where $h$ is $W$-generic for $Q$.  I first show that $V\subset W$.  Since $Q$ is a set, $W$ and $V[\class{G}]$ must agree on both cardinals and the continuum function above $|Q|$.  Fix any $x$ in $V$.  Certainly $x$ is coded into the continuum function of $V[\class{G}]$, and since the coding is duplicated arbitrarily high up a copy of the code must appear starting above $|Q|$.  Thus the code appears also in $W$, and so $x\in W$.  Thus $V \subset W$.

Now choose $\kappa>|Q|$ a strong limit, and factor $\P = P_\kappa * \P_{tail}$ as before, giving $V[\class{G}] = V[G_\kappa][\class{G}_{tail}]$.  Once again I can find $A\subset \kappa$ in $W$ such that $A$ codes both $(2^{<\kappa})^W$ and $\dot{G_\kappa}$, a $Q$-name for $G_\kappa$.

The remainder of the proof follows as above.  I get $(2^{<\kappa})^{V[A]} = (2^{<\kappa})^W$ by the same argument, relying on $V\subset W$ for the forward inclusion.  Thus $h$ is $V[A]$-generic for $Q$, and $V[A][h]=V[G_\kappa]$.  For any successor $\delta$ between $|Q|^+$ and $\kappa$, both $V[A]\subset V[A][h][\class{G}_{tail}]$ and $W\subset W[h]$ satisfy the $\delta$ cover and approximation properties, and $V[A]$ and $W$ have the same $\mathcal{P}(\delta)$.  Thus $V[A]=W$.  Finally, $W[h]=V[A][h]=V[G_\kappa]$, contradicting $W[h]=V[G_\kappa][\class{G}_{tail}]$.  \end{proof}

\section{Consistency with large cardinals}

The argument given above is flexible enough to allow the preparatory coding to begin arbitrarily high in the ordinals.  This allows preservation of many large cardinal notions.

\begin{Cor}  If $V$ satisfies $ZFC + \kappa$ is supercompact, then there is a proper-class-forcing extension of $V$ satisfying $\kappa$ is supercompact $+\GA+\nVHOD$.\end{Cor}

\begin{proof}  Start in $V$ with $\kappa$ supercompact.  Begin by forcing the supercompactness of $\kappa$ to be indestructible by ${<\kappa}$-directed closed forcing.  This can be accomplished by means of the Laver Preparation \cite{Laver78}.  I then follow the proof of Theorem \ref{Theorem.GA+VneqHODoverV}, starting both the preparatory forcing and the iteration $\P$ at $\kappa$.   Both the preparatory forcing and $\P$ will be $<\kappa$-directed closed, and the argument used in the proof of Theorem \ref{T:GAand-GCHandSupercompact} shows that $\kappa$ remains supercompact in the extension.  The argument that the extension satisfies $\GA + \nVHOD$ goes through as before. \end{proof}

\section{Open Questions}

A number of questions regarding these axioms remain.  Of particular interest is the question of uniqueness of bedrock models.

\begin{Question}If the Bedrock Axiom holds, is the bedrock model unique?
\end{Question}

A negative answer would settle another natural question:  Given two ground models of $V$ can we always find a third ground model contained in their intersection?  If the answer is yes, it would indicate that forcing can be used to amalgamate only models that are in some sense ``close together.''  Another natural structure to consider is the class obtained by intersecting all ground models.  Is it a model of $\ZFC$?  In many examples the answer is yes, even in cases as in the previous section where the Bedrock Axiom fails.  The models of $\neg\BA$ described above also have the property that there are a large number of ground models, proper class many.  Is it possible to have $\neg\BA$ with only set many ground models?  If so, can we reduce the number of ground models to countable?

Another area of interest is analyzing the restriction of these axioms to various classes of forcing such as {\sc ccc} or proper forcing.  For example, $\GA_{ccc}$ is the assertion that the universe is not a set-forcing extension of an inner model by {\sc ccc} forcing.  Can we separate notions of the Ground Axiom by obtaining models of, for example, $\GA_{ccc} + \neg\GA_{proper}$?  How about consistency of restricted notions of $\GA$ and $\neg\BA$, such as $\GA_{ccc} + \neg\BA$?

Also of interest is the extension of these axioms to include class forcing.  Of course the class versions of the Ground and Bedrock Axioms can be easily stated in second-order terms, but first-order equivalents have not yet been formulated.  While this may prove impossible in the general case, there is some hope that by restricting attention to a particular class of class forcing, for example, forcing with a closure point at $\delta$, a first-order expression may be possible.

\clearpage

\begin{appendices}
\chapter{Forcing with proper classes}\label{S:ClassForcing}

The main contribution of this dissertation is now complete.  The material of this appendix should be considered merely supplementary to the main part of the dissertation, consisting of Chapters \ref{C:TheGroundAxiom} and \ref{C:GAandVnotHOD}.  Many of the theorems in those chapters depend on forcing with proper classes.  In standard texts, forcing is usually presented first for sets, with complete proofs of all the relevant theorems.  The generalization to class forcing is then sketched, with descriptions given of the appropriate modifications to the proofs.  My goal in this appendix is to develop the theory of forcing with classes from the ground up, presenting the `classes' versions of the theorems, e.g. the Forcing Theorem and the Generic Model Theorem, in their entirety.   However, it should be noted that the use of class forcing in Chapters \ref{C:TheGroundAxiom} and \ref{C:GAandVnotHOD} is of a largely unproblematic variety, namely, progressively closed Easton support products and iterations, and my presentation will be confined to these and similar cases.  For a precise statement of the restrictions placed on the class partial orders considered in this appendix see Definition \ref{D:CCP}, with additional requirements in the cases of iterations and products given in Definitions \ref{D:IterationRestrictions} and \ref{D:ProductRestrictions}, respectively.  While definitions of basic terms will be provided, some familiarity with forcing and its accompanying machinery of partial orders, boolean algebras, $P$-names, etc. will be helpful.  For example, this exposition might be appropriate for the reader who is familiar with forcing in the set case and wishes to see the theory of forcing with classes presented in full.  

A rigorous treatment of forcing with classes suggests the use of some axiomatic theory beyond the first-order set theory $\ZFC$.  Statements of the form ``for every partially ordered class $\P$, the following assertion holds'' require quantification over classes $\P$, beyond the scope of $\ZFC$.  One alternative is to think of such statements as schemes, including for each formula $\phi(x)$ the statement ``if the class defined by $\phi(x)$ is a partial order, then the following assertion holds.''  However, a more attractive alternative is to use a true second-order axiomatization of set theory.  I will use that of Bernays-G\"odel, or $\BGC$, and a statement of the axioms of $\BGC$ will appear in Section \ref{SS:BGC}.  This theory is attractive for several reasons.  Any model $V \vDash \ZFC$ is also a model of $\BGC$, whose classes are simply the definable classes of $V$.  Furthermore, $\BGC$ is a conservative extension of $\ZFC$, so theorems provable in $\BGC$ \emph{using only set variables} are provable in $\ZFC$.  

The material in this appendix is a synthesis of Jech \cite{Jech:SetTheory3rdEdition}, Chapters 7 and 14, and Kunen \cite{Kunen:Independence}, Chapters 7 and 8, both excellent texts.  Kunen's treatment provides a very lucid and friendly introduction to the theory of forcing with partially ordered sets, and a detailed look at products and iterations.  Jech takes the Boolean-valued model approach to forcing and provides great background on the theory of Boolean algebras, as well as an introduction to Bernays-G\"odel set theory.

\section{Bernays-G\"odel set theory}\label{SS:BGC}

In Bernays-G\"odel set theory, there are two types of objects, \emph{sets}, for which I will use standard letters $x, A, \alpha$, etc., and \emph{classes}, which will be denoted by blackboard bold letters $\class{X}, \class{Y}$, etc..  An exception to this rule will be the class of all sets, which in keeping with standard notation I will denote $V$. As in the usual first-order set theory, there is a single binary relation $\in$.  The axioms of Bernays-G\"odel set theory, or $\BGC$, are as follows.

\clearpage
\begin{Def}\rm\label{D:BGC}\emph{Bernays-G\"odel set theory} ($\BG$) consists of axioms \ref{D:BGC:Extensionality}--\ref{D:BGC:Regularity} below.  \emph{Bernays-G\"odel set theory with Choice} ($\BGC$) also includes axiom \ref{D:BGC:Choice}.
\begin{enumerate}
\parskip=0pt

\item \label{D:BGC:Extensionality}$\forall u \({u\in \class{X} \iff u\in\class{Y}}) \rightarrow \class{X}=\class{Y}$ \emph{ (extensionality)}. 

\item \label{D:BGC:SetsAreClasses}Every set is a class.

\item \label{D:BGC:OnlySetsAreMembers}If $\class{X}\in\class{Y}$, then $\class{X}$ is a set.

\item \label{D:BGC:Pairing}For any sets $x$ and $y$ there is a set $\{x,y\}$ \emph{ (pairing)}. 

\item \label{D:BGC:Comprehension}For any formula $\phi$ in which only set variables are quantified, $$\forall \class{X}_1, ..., \class{X}_n \exists \class{Y} \; \({\class{Y}=\{x \mid \phi(x, \class{X}_1, ...,\class{X}_n)\} })\mbox{\emph{ (comprehension)}.}$$

\item \label{D:BGC:Infinity}There is an inductive set \emph{ (infinity)}.

\item \label{D:BGC:Union}For every set $x$ the set $\bigcup x$ exists \emph{ (union)}.

\item \label{D:BGC:PowerSet}For every set $x$ the power set $\mathcal{P}(x) = \{\class{Y} \mid \class{Y} \subset x\}$ of $x$ exists \emph{ (power set)}.

\item \label{D:BGC:Replacement}If a class $\class{J}$ is a function and $x$ is a set, then $\{\class{J}(z) \mid z\in x\}$ is a set \emph{ (replacement)}.

\item \label{D:BGC:Regularity}Every nonempty class has an $\in$-minimal element \emph{ (regularity)}.

\item \label{D:BGC:Choice}There is a function $\class{F}$ such that $\class{F}(x)\in x$ for every nonempty set $x$ \emph{ (choice)}.

\end{enumerate}
\end{Def}

A model $\<V,\V,\in>$ of $\BGC$ consists of a collection $\V$ of classes together with a subcollection $V\subset \V$ of sets, and a relation ${\in} \subset \(V\times\V)$.  Depending on the context, the model may be referred to as $\<V,\V>$ or simply $\V$.  Note that $\V$ and $\in$ completely determine the model, for the sets $V$ are definable as $V=\left\{\class{X}\mid \exists \class{Y} \; \({\class{X}\in\class{Y}})\right\}$.  From a model of $\BGC$ one can produce a model of $\ZFC$ by simply `throwing away $\V$,' that is, by taking $V$ and ${\in}\cap \(V\times V)$.  Conversely, if $V$ is a model of $\ZFC$, then taking $\V$ to be the collection of all classes definable in $V$ with set parameters, and taking $\in$ to be the obvious extension of $V$'s epsilon relation to $\V$, then $\<V,\V,\in>$ is a model of $\BG$.  While the model may not satisfy the uniform version of the Axiom of Choice given in Axiom \ref{D:BGC:Choice}, one can add a uniform choice function through a class forcing that is $\kappa$-closed for every $\kappa$ and so adds no sets to the universe.  The resulting model $\<V,\V^\prime,\in^\prime>$ of $\BGC$ extends $\<V,\V,\in>$ and has the same sets and the same  restriction of $\in$ to sets.  A nice consequence of this, and the essential property of $\BGC$ for the purposes of this exposition, is the following.

\begin{Thm}\label{T:BGCandZFC}If $\phi$ is a formula in the language of Bernays-G\"odel set theory in which no class variables appear, then $(\BGC \proves \phi) \iff (\ZFC \proves \phi)$.
\end{Thm}

\section{Partially ordered classes}\label{SS:POclasses}

Note that, as every set is a class, the definitions below apply equally to both sets and proper classes.

\begin{Def}\rm \label{D:POclass} A \emph{partially ordered class}, or \emph{partial order}, is a class $\P\in \V$ consisting of ordered pairs $\<p,q>$ defining a relation $p \leq q$ satisfying
\begin{enumerate}
\parskip=0pt
\item \label{D:POclass:reflexive}$\forall p\in\P \, \(p\leq p)$ \emph{ (reflexive)}.
\item \label{D:POclass:transitive}$\forall p, q, r \in \P \, \(p\leq q \wedge q \leq r \rightarrow p \leq r)$ \emph{ (transitive)}.

\setcounter{saveenum}{\value{enumi}} 
\end{enumerate} 
$\P$ is a  partial order \emph{in the strict sense} if and only if
\begin{enumerate}
\setcounter{enumi}{\value{saveenum}} 
\parskip=0pt
\item \label{D:POclass:strict}$\forall p, q \in \P \, \(p \leq q \wedge q \leq p \rightarrow p=q)$ \emph{ (antisymmetric)}.
\end{enumerate}
If $\P$ is a set then $\P$ is sometimes referred to as a \emph{poset}.
\end{Def}

There is disagreement in the literature, with some authors taking conditions \ref{D:POclass:reflexive} and \ref{D:POclass:transitive} to define a `partial pre-order,' and the addition of condition \ref{D:POclass:strict} defining a partial order.  As iterations rarely satisfy \ref{D:POclass:strict}, however, I have opted to exclude this condition from the definition of partial order. Note that the notation is abused by using $\P$ to represent the collection of underlying elements, allowing statement such as $p\in\P$.  The usual terminology of forcing and partial orders applies equally in the class context.

\begin{Def}\rm\label{D:POclassProperties} For $\P$ a partial order and $p, q \in \P$:
\begin{enumerate}
\parskip=0pt
\item $p$ is \emph{stronger than} $q$, or $p$ \emph{extends} $q$, if and only if $p\leq q$.
\item $p$ and $q$ are \emph{compatible}, written $p \parallel q$, if and only if there is $r\in\P$ extending both $p$ and $q$.
\item If $p$ and $q$ are not compatible they are \emph{incompatible}, $p\perp q$.
\item A subclass $\class{A}\subset \P$  is an \emph{antichain} if and only if the elements of $\class{A}$ are pairwise incompatible.
\item An antichain $\class{A}$ is \emph{maximal} if and only if there is no $p\in \P$ such that $\{p\} \cup \class{A}$ is an antichain.
\item A subclass $\class{A}\subset \P$ is \emph{open} if and only if it is closed downwards, that is, $\forall p\in \class{A} \, \forall q\in \P \, \(q \leq p \rightarrow q\in \class{A})$.
\item A subclass $\class{D}\subset \P$ is \emph{dense} if and only if it intersects every open set $\P$,  that is, $\forall p\in \P \, \exists q\in \class{D} \, \(q\leq p)$.
\item For $\kappa$ a regular cardinal, $\P$ satisfies the \emph{$\kappa$ chain condition}, or $\kappa$-c.c., if and only if every antichain of $\P$ has cardinality $<\kappa$.
\item For $\kappa$ a regular cardinal, $\P$ is \emph{$<\!\kappa$-closed} if and only if for every descending sequence $\<p_\beta \mid \beta < \alpha>$ in $\P$ of length $\alpha < \kappa$ there is $p\in \P$ lying below every $p_\beta$.
\end{enumerate}
\end{Def}

Forcing with $\P$ will adjoin to the universe a \emph{generic filter} $\class{G}\subset \P$.

\begin{Def}\rm\label{D:POClassFilters} For $\P$ a partial order, a subclass $\class{G}\subset \P$ is a \emph{filter} if and only if the following conditions hold.
\begin{enumerate}
\parskip=0pt
\item $\forall p, q \in \class{G} \, \exists r\in \class{G}$ extending $p$ and $q$ \emph{ (directed)}.
\item $\forall p \in \class{G} \, \forall q \in \P \, \(p \leq q \rightarrow q \in \class{G})$ \emph{ (upwards closed)}.
\end{enumerate}
\end{Def}

\begin{Def}\rm\label{D:POGenericFilters}
A filter $\class{G}\subset \P$ is \emph{$\V$-generic} if and only if any of the following equivalent conditions holds.
\begin{enumerate}
\parskip=0pt
\item For every dense class $\class{D}\subset \P \, \({\class{D} \cap \class{G} \ne \emptyset})$.
\item For every open dense class $\class{D}\subset \P \, \({\class{D} \cap \class{G} \ne \emptyset})$.
\item For every maximal antichain $\class{A} \subset \P \, \({\class{A} \cap \class{G} \ne \emptyset})$.
\end{enumerate}
\end{Def}

Equivalence of the three definitions of generic filter is a standard exercise and is left to the reader.

Consistency of the existence of $\V$-generic filters can be established by applying essentially the same metamathematical arguments as in the set forcing case (see Section \ref{SS:ConsistencyOfGenericFilters}).  However, unlike set forcing, the resulting model $\V[\class{G}]$ will not necessarily satisfy $\BGC$.  Later in the text I will consider two different restrictions on $\P$ that guarantee $\BGC$ in the extension $\V[\class{G}]$, each useful in a different type of forcing argument.  However, there is a very basic underlying assumption on $\P$ which will be required even for the preliminary development.  This assumption, that $\P$ is a \emph{chain of complete subposets}, appears below as Definition \ref{D:CCP}.  It is based on the notion of \emph{ complete embeddings}.

\begin{Def}\rm\label{D:POEmbeddings} If $\P$ and $\Q$ are partial orders and $i: \P \rightarrow \Q$, then $i$ is an \emph{embedding} if and only if
\begin{enumerate}
\parskip=0pt
\item $\forall p, q \in \P \, \left( p \leq q \rightarrow i(p) \leq i(q) \right)$ \emph{ ($i$ preserves $\leq$)}, and
\item $\forall p, q \in \P \, \left( p \perp_\P q \rightarrow i(p) \perp_\Q i(q) \right)$ \emph{ ($i$ preserves incompatibility)}.
\end{enumerate}
The basic notion of embedding can be strengthened in various ways.
\begin{enumerate}
\parskip=0pt
\item An embedding $i$ is \emph{dense} if the image of $i$ is a dense subclass of $\Q$.
\item An embedding $i$ is \emph{complete} if and only if for every maximal antichain  $\class{A} \subset \P$, the image $i^{\prime\prime} \class{A}$ is a maximal antichain of $\Q$.  Note that every dense embedding is complete.
\item If $\P \subset \Q$ and $i$ a complete embedding is the identity map, then $\P$ is a \emph{complete suborder} of $\Q$, denoted $\P \subset_c \Q$.
\end{enumerate}
\end{Def}

The primary interest in dense and complete embeddings lies in the relationship to generic filters.  A dense embedding implies that the two posets are equivalent in terms of forcing, as a generic for one always gives a generic for the other.

\begin{Lemma}\label{L:PODenseEmbeddingsAndGenericFilters} Suppose $i: \P \rightarrow \Q$ is a dense embedding.
\begin{enumerate}
\parskip=0pt
\item If $\class{G}\subset \Q$ is $\V$-generic for $\Q$, then $i^{-1} \class{G}$ is $\V$-generic for $\P$.
\item If $\class{G}\subset \P$ is $\V$-generic for $\P$, then $\class{G}^\prime = \{q \in \Q \mid \exists p\in \class{G} \left( i(p) \leq q \right) \}$ is $\V$-generic for $\Q$.  Note that $\class{G}^\prime$ is the `upwards closure' of the image of $\class{G}$ in $\Q$.
\end{enumerate}
\end{Lemma}

The proof is left as an exercise.

Complete embeddings are similar, but the implication goes only one direction.  A generic for $\Q$ always gives a generic for $\P$.

\begin{Lemma}\label{L:POCompleteEmbeddingsAndGenericFilters} If $i:\P \rightarrow \Q$ is a complete embedding and $\class{G} \subset \Q$ is a $\V$-generic filter for $\Q$, then $i^{-1} \class{G}$ is a $\V$-generic filter for $\P$.
\end{Lemma}
\begin{proof}  That $i^{-1} \class{G}$ is a filter straightforward and follows from the properties of embeddings given above.  It remains to show that $i^{-1} \class{G}$ meets every maximal antichain $\class{A}\subset \P$.  Fix such an $\class{A}$, and it follows by completeness of $i$ that $i^{\prime\prime} \class{A}$ is maximal in $\Q$.  By genericity of $\class{G}$ for $\Q$ there is $q\in i^{\prime\prime} \class{A} \cap \class{G}$.  Clearly $i^{-1}(q) \in i^{-1} \class{G}  \cap  \class{A}$. \end{proof}

The underlying assumption made on partially ordered classes is that they can be written as the union of a sequence of partially ordered sets, each a complete suborder of those that follow.  

\begin{Def}\rm \label{D:CCP} A partial order $\P$ is a \emph{chain of complete subposets} if and only if $\, \P = \bigcup_{\alpha\in\ORD} P_\alpha$, a union of  partially ordered sets $P_\alpha$, such that $\alpha \leq \beta \rightarrow P_\alpha \subset_c P_\beta$.  Furthermore, the sequence $\{\<\alpha,P_\alpha>\mid \alpha \in \ORD\}$ is a class. 
\end{Def}

Note there is no \emph{continuity} assumption made at limits, i.e.\ that $P_\lambda = \bigcup_{\alpha < \lambda} P_\alpha$.  In fact, this assumption is quite strong, and neither products nor iterations using Easton support, such as those employed throughout this paper and in many common class forcing arguments, will satisfy it.  Easton support is used in order to preserve closure conditions, a requirement for many forcing proofs.

It is not difficult to see that for a partial order $\P$ that is a chain of complete subposets, each $P_\alpha$ is a complete subposet of $\P$.

\begin{Lemma}\label{L:EachSubposetIsCompleteInP} If $\P = \bigcup_{\alpha\in\ORD} P_\alpha$ is a chain of complete subposets, then $P_\alpha \subset_c \P$ for all $\alpha$.
\end{Lemma}

\begin{proof}  I will show that maximal antichains of $P_\alpha$ remain maximal in $\P$.  The other properties, such as preservation of ordering and incompatibility, are proved in a similar fashion.  Suppose $A\subset P_\alpha$ is a maximal antichain, but $A$ is not maximal in $\P$.  Then there is $p \in \P$ such that $\{p\} \cup A$ is an antichain in $\P$.  The set $\{p\} \cup A$ is contained in $P_\beta$ for some $\beta \geq \alpha$, and it remains an antichain there.  However, $P_\alpha$ is a complete subposet of $P_\beta$ and so $A$ is a maximal antichain in $P_\beta$.  Thus $\{p\} \cup A$ cannot be an antichain in $P_\beta$, a contradiction. \end{proof}

\section{Boolean algebras}\label{SS:BooleanAlgebras}

The theory of forcing for sets is often developed using Boolean algebras in the place of partial orders.  The two approaches are equivalent, as will be shown below, and the Boolean algebra approach gives rise to the pleasingly intuitive notion of Boolean-valued models.  However, to take full advantage of Boolean-valued models requires the use of \emph{complete} Boolean algebras.  While every set Boolean algebra has a completion, this is not true of many of the proper class Boolean algebras considered here.  However, the Boolean algebra approach still provides simplified proofs of many of the basic forcing theorems.  I will therefore use a hybrid approach, employing Boolean algebras where possible and relying on partial orders otherwise.  I will begin by outlining the basics of Boolean algebras.  Notice that many of the properties defined above for partial orders have corresponding notions in the Boolean algebra context.

\clearpage
\begin{Def}\rm\label{D:BooleanAlgebra} A \emph{Boolean algebra} is a class $\B$ with at least two elements $0$ and $1$, two binary operations $+$ and $\btimes$, and one unary operation $-$, satisfying for all $u,v,w \in \B$
\begin{enumerate}
\parskip=0pt
\item $u+v=v+u$ and $u\btimes v = v\btimes u$\emph{ (commutativity)}.
\item $u+(v+w) = (u+v)+w$ and $u\btimes(v\btimes w) = (u\btimes v)\btimes w$  \emph{ (associativity)}.
\item $u\btimes (v+w) = u\btimes v + u \btimes w$ and $u + (v\btimes w) = (u+v)\btimes(u+w)$ \emph{ (distributivity)}.
\item $u \btimes (u+v) = u$ and $u + (u\btimes v)=u$ \emph{ (absorption)}.
\item $u+(-u)=1$ and $u\btimes (-u)=0$ \emph{ (complementation)}.
\end{enumerate}
\end{Def}

From these properties, many other useful facts can be derived.  The following will be useful in later proofs.

\begin{Lemma}\label{L:BAproperties} If $\B$ is a Boolean algebra, then for all $u,v \in \B$
\begin{enumerate}
\parskip=0pt
\item $u\btimes u = u + u = u$.
\item $u+0 = u \btimes 1 = u$.
\item $u\btimes 0 = 0$ and $u + 1 = 1$.
\item $-(u+v) = -u \btimes -v$ and $-(u\btimes v) = -u + -v$ \emph{ (De Morgan's laws)}.
\end{enumerate}
\end{Lemma}

Boolean algebras are often equipped with additional functions and relations useful in various contexts  The second definition below puts a partial order structure on every Boolean algebra.

\begin{Def}\rm\label{D:BAordering} Suppose $\B$ is a Boolean algebra. 
\begin{enumerate}
\parskip=0pt
\item The operation $u-v$ is defined by $u-v = u\btimes (-v)$.
\item The relation $u\leq v$ is defined by $u\leq v$  if and only if $u-v=0$.
\item The binary operation $u \rightarrow v$ is defined by $u\rightarrow v = -u+v$.
\end{enumerate}
\end{Def}

It can be useful to think of the properties of Boolean algebras in terms of the ordering, and of the partial order properties in terms of the Boolean operations.  For forcing purposes, it generally useful to consider the partial order formed by a Boolean algebra with the $0$ element removed, $\B \setminus \{0\}$.

\begin{Lemma}\label{L:PropertiesofBAsfromOrdering} Suppose $\B$ is a Boolean algebra with ordering $\leq$ as defined above.
\begin{enumerate}
\parskip=0pt
\item $1$ is the greatest element of $\B$.
\item $0$ is the least element of $\B$.
\item $u+v$ is the least upper bound of $\{u,v\}$, and $u\btimes v$ is the greatest lower bound of $\{u,v\}$.
\item $-u$ is the unique element $v$ such that $u+v=1$ and $u\btimes v=0$.
\item $u$ and $v$ are incompatible in the partial order $\B \setminus \{0\}$ if and only if $u\btimes v =0$.
\end{enumerate}
\end{Lemma}

Embeddings of Boolean algebras are more restrictive than embeddings of partial orders.  In particular, it not sufficient that an embedding preserve the ordering, it must also preserve the Boolean operations.

\begin{Def}\rm\label{D:BAEmbeddings} If $\B$ and $\D$ are Boolean algebras and $i: \B \rightarrow \D$, then $i$ is an \emph{embedding} if and only if 
\begin{enumerate}
\item $i(0_\B)=0_\D$ and $i(1_\B)=1_\D$,
\setcounter{saveenum}{\value{enumi}}
\end{enumerate}
and for and for all $u,v \in \B$,
\begin{enumerate}
\setcounter{enumi}{\value{saveenum}}
\item $i(u)+i(v)=i(u+v)$, and
\item $i(u)\btimes i(v)=i(u\btimes v)$, 
\item $-i(u)=i(-u)$.
\end{enumerate}
\end{Def}

Complete embeddings of Boolean algebras require an additional restriction.  They are based on a generalization of the operations $+$ and $\btimes$ to sums and products of many elements, which in turn are defined in terms of the ordering on $\B$.

\begin{Def}\rm\label{D:BASumsandProducts} For $\B$ a Boolean algebra and $\class{X}\subset \B$, 
\begin{enumerate}
\parskip=0pt
\item $\sum \class{X} =$ the least upper bound of $\class{X}$, provided it exists, and $\sum \emptyset = 0$.
\item $\prod \class{X} =$ the greatest lower bound of $\class{X}$, provided it exists, and $\prod \emptyset = 1$.
\end{enumerate}
\end{Def}

Generalized sums and products also have certain nice properties with respect to the Boolean operations.

\begin{Lemma}\label{L:BASumsandProductsProperties} Suppose $\B$ is a Boolean algebra, $\class{X}\subset \B$ and the sum $\sum \class{X}$ and product $\prod \class{X}$ exist.  Then for all $u\in \B$,
\begin{enumerate}
\parskip=0pt
\item $u + \sum \class{X} = \sum \{u+v \mid v\in \class{X}\}$,
\item $u \btimes \prod \class{X} = \prod \{u\btimes v \mid v\in \class{X}\}$,
\item $u \btimes \sum \class{X} = \sum \{u\btimes v \mid v\in \class{X}\}$, and
\item $u + \prod \class{X} = \prod \{u + v \mid v\in \class{X}\}$.
\setcounter{saveenum}{\value{enumi}}
\end{enumerate}
Furthermore, if $\class{Y}\subset \B$ and $\sum \class{Y}$ exists, then 
\begin{enumerate}
\parskip=0pt
\setcounter{enumi}{\value{saveenum}}
\item $\({\sum \class{X}}) \btimes \({\sum \class{Y}}) = \sum \{v\btimes w \mid \<v,w>\in \class{X} \times \class{Y}\}$. 
\end{enumerate}
\end{Lemma}

Note that many of the above properties are generalizations of the distributivity laws, and hold in all Boolean algebras.  However, the most general distributive laws are quite strong and do not necessarily hold true even in complete Boolean algebras.  Distributivity of a Boolean algebra has important consequences in terms of its forcing extensions, but these are not directly related to development of class forcing presented here.

The notions of \emph{dense subset, open subset, antichain, maximal antichain, filter} and \emph{generic filter} in the Boolean algebra context can be obtained by applying the definitions given in Section \ref{SS:POclasses} to the partial order $\B \setminus \{0\}$.  In the theory of Boolean algebras, alternative terminology is used for many of these concepts and definitions are given in terms of the Boolean operations rather than in terms of the ordering, for example, incompatible elements are \emph{disjoint} and maximal antichains is are \emph{partitions}.  While it not strictly necessary to consider the Boolean algebra versions of all of these structures, certain algebraic properties will be extremely useful in later proofs.  I will present the most important definitions below.

\begin{Def}\rm\label{D:BAfilters}Suppose $\B$ is a Boolean algebra and $\class{F}\subset \B$. 
\begin{enumerate}
\parskip=0pt
\item $\class{F}$ is a \emph{filter} if and only if 
	\begin{enumerate}
	\parskip=0pt
	\item $1\in \class{F}$, 
	\item $0 \notin \class{F}$, 
	\item $\class{F}$ is closed upwards, that is, $u\in \class{F}$ and $v \geq u \rightarrow v\in \class{F}$, and 
	\item $\class{F}$ is closed under finite products, that is, $u,v \in \class{F} \rightarrow u\btimes v \in \class{F}$.
	\end{enumerate}
\item $\class{F}$ is an \emph{ultrafilter} if and only if $\class{F}$ is a filter and for all $u\in \B$, either $u\in \class{F}$ or $-u \in \class{F}$.
\item $\class{F}$ is a \emph{$\V$-generic ultrafilter} if and only if $\class{F}$ is an ultrafilter and for all $\class{X}\in \V$, if $\class{X}\subset \class{F}$, then $\prod \class{X} \in \class{F}$, provided the product exists.
\end{enumerate}
\end{Def}

The correspondence between the notions of generic ultrafilter on a Boolean algebra and generic filter on a partial order is stated below.  The proof is left as an exercise. 

\begin{Lemma}\label{L:GenericUltrafilter=GenericFilter}If $\B$ is a Boolean algebra and $\class{F}\subset \B$, then $\class{F}$ is a $\V$-generic ultrafilter on $\B$ if and only if $\class{F}$ is a $\V$-generic filter on the partial order $\B\setminus \{0\}$. 
\end{Lemma}

The existence of arbitrary sums and products in a Boolean algebra is captured in the notion of \emph{completeness} of Boolean algebras, not to be confused with completeness of embeddings of partial orders, defined above, or completeness of embeddings of Boolean algebras, defined below.

\begin{Def}\rm\label{D:BACompleteness} Suppose $\B$ is a  Boolean algebra. 
\begin{enumerate}
\parskip=0pt
\item $\B$ is \emph{complete} if $\sum \class{X}$ exists for all $\class{X}\subset \B$.
\item $\B$ is \emph{$\kappa$-complete} for $\kappa$ a regular cardinal if $\sum X$ exists for all $X\subset \B$ of size $<\!\kappa$.
\end{enumerate}
\end{Def}

A complete embedding of Boolean algebras is an embedding that preserves arbitrary sums and products.

\clearpage
\begin{Def}\rm\label{D:BACompleteEmbeddings} An embedding $i: \B \rightarrow \D$ of Boolean algebras is \emph{complete} if and only if for every $\class{X}\subset \B$,
\begin{enumerate}
\parskip=0pt
\item $\sum \({i^{\prime\prime}\class{X}}) = i\({ \sum \class{X}})$, if either exists.
\item $\prod \({i^{\prime\prime}\class{X}}) = i\({ \prod \class{X} })$, if either exists.
\end{enumerate}
If the map $i$ is the identity, then $\B$ is a \emph{complete subalgebra} of $\D$, denoted $\B \subset_c \D$.
\end{Def}

Note that in the usual exposition complete embeddings are defined only for complete Boolean algebras, and the existence requirement is dropped in the two conditions of Definition \ref{D:BACompleteEmbeddings}.

It is an exercise that a complete embedding $i: \B \rightarrow \D$ of Boolean algebras is also a complete embedding of the associated partial orders  $i: \B \setminus \{0\} \rightarrow \D\setminus \{0\}$ in the sense of Definition \ref{D:POEmbeddings}.  The analogues of Lemmas \ref{L:PODenseEmbeddingsAndGenericFilters} and \ref{L:POCompleteEmbeddingsAndGenericFilters} also hold for Boolean algebras, e.g.

\begin{Lemma}\label{L:BACompleteEmbeddingsAndGenericFilters} If $i:\B \rightarrow \D$ is a complete embedding of Boolean algebras and $\class{G} \subset \D$ is a $\V$-generic ultrafilter for $\D$, then $i^{-1} \class{G}$ is a $\V$-generic ultrafilter for $\P$.
\end{Lemma}

\section[Equivalence of partial orders and Boolean algebras]{Equivalence of partial orders and \\Boolean algebras}\label{SS:POsandBAs} 

Definition \ref{D:BAordering} provides a natural way of viewing any Boolean algebra as a partial order.  While this observation is not reversible, for it is not the case that every partial order is a Boolean algebra, it is close to true in the case of sets.  Lemma \ref{L:ro(P)} will show that every set partial order embeds densely into a complete Boolean algebra.  The construction of a complete Boolean algebra from a partially ordered set involves a notion called \emph{regular open sets}.

\begin{Def}\rm Suppose $P$ is a partially ordered set. 
\begin{enumerate}
\parskip=0pt
\item A subset $A\subset P$ is \emph{open} if and only if it is closed downwards, i.e. $\forall p\in A \, \forall q\in P \, \(q\leq p \rightarrow q\in A)$.
\item For each $p\in P$, the \emph{basic open set} $P_p$ is defined $P_p = \{q\leq p\}$.
\item An open set $A$ is \emph{regular} if and only if $\forall q \notin A \, \exists p\in A \, (q \perp p)$.
\end{enumerate}
\end{Def}

The construction of a Boolean algebra relies on certain nice properties of regular open sets.  

\clearpage
\begin{Lemma}\label{L:RegularOpenSets}  Suppose $P$ is a partially ordered set.
\begin{enumerate}
\parskip=0pt
\item If $\mathcal{A}$ is a collection of regular open sets of $P$, then $\bigcap \mathcal{A}$ is a regular open set.
\item Every subset $A\subset P$ is contained in a least regular open set $\overline{A}^\circ$, the intersection of all regular open sets containing $A$.  A subscript is sometimes used to indicate the partial order in which this operation takes place, e.g. $\overline{A}_P^\circ$.
\item Equivalently, $\overline{A}^\circ$ can be defined $\overline{A}^\circ = \{p\in P \mid \forall q\leq p \( P_q \cap A) \ne \emptyset\}$.
\item If $C$ is a \emph{maximal antichain of A}, that is, $C\subset A$ is an antichain of $P$ and no element of $A$ can be added to $C$ while preserving pairwise incompatibility, then $\overline{A}^\circ = \overline{C}^\circ$.
\end{enumerate}
\end{Lemma}

The collection of all regular open sets of a partially ordered set forms a complete Boolean algebra.  For partially ordered classes, an open set may be a proper class, which prevents us from completing the construction in the class case.

\clearpage
\begin{Lemma}\label{L:ro(P)} If $P$ is a partially ordered set and $ro(P)$ is the collection of all regular open sets of $P$, called the \emph{regular open algebra of $P$}, then $ro(P)$ together with the operations
\begin{enumerate}
\parskip=0pt
\item $u\btimes v = u \cap v$
\item $u + v = \overline{u \cup v}^\circ$
\item $-u = \{p \mid P_p \cap u = \emptyset\}$
\end{enumerate}
forms a complete Boolean algebra with $1=P$ and $0=\emptyset$.  Furthermore, the map $e:P \rightarrow ro(P)$ given by $e(p)=\overline{P}^\circ_p$ is a dense and therefore complete embedding of the partial order $P$ into the partial order $ro(P) \setminus \{0\}$.
\end{Lemma}

Lemma \ref{L:ro(P)} shows that every partially ordered set can be densely embedded into a complete Boolean algebra, and it follows from Lemma \ref{L:PODenseEmbeddingsAndGenericFilters} that $P$ and $ro(P)$ are forcing equivalent.  It is for this reason that the development of set forcing can be restricted to the consideration of complete Boolean algebras.  Unfortunately, this is not true in the case of proper classes.  However, for a partially ordered class that is a chain of complete subposets, complete Boolean algebras can still be of some use.  This is because many of the arguments about $\P$ can be restricted to arguments about some $P_\alpha$, in which case it is sufficient to consider the complete Boolean algebra $ro(P_\alpha)$.  A nice property of the $ro(P)$ construction is that it preserves complete embeddings.

\begin{Lemma}\label{L:POEmbeddingsInduceBAEmbeddings} If $i:P \rightarrow Q$ is a complete embedding of posets and the maps $e_P:P \rightarrow ro(P)$ and $e_Q:Q \rightarrow ro(Q)$ are the usual embeddings, then $j:ro(P) \rightarrow ro(Q)$ defined by $j(u) = \overline{i^{\prime\prime} u}^\circ_Q$ is a complete embedding of Boolean algebras and the following diagram commutes:

 \[
 \begin{CD}
P @>i>> Q \\
@V{e_P}VV @V{e_Q}VV\\
{ro(P)} @> {j} >> {ro(Q)} \\
\\
 \end{CD}
 \]

\end{Lemma}

It would be convenient in the proofs that follow to have $P \subset_c Q$ implies $ro(P) \subset_c ro(Q)$.  Unfortunately this is not the case, as a regular open set of $P$ may not even be open in $Q$.  However, it worthwile to consider an alternative $ro^*(Q)$ for which the embedding $j$ is the identity.

\begin{Def}\rm\label{D:ro*(Q)} Suppose $P$ and $Q$ are posets and $P \subset_c Q$.  For each $C\in ro(Q)$, let $C^* = A$ such that $A\in ro(P)$ and $\overline{A}^\circ_Q = C$, if such $A$ exists.  Otherwise let $C^* = C$.  Let $ro^*(Q)$ be the collection of all $C^*$.  A natural Boolean algebra structure is induced on $ro^*(Q)$ by the bijection $C \mapsto C^*$, giving $ro(Q) \cong ro^*(Q)$.
\end{Def}

Note that in the Definition, if $A$ exists in $ro(P)$ with $\overline{A}^\circ_Q = C$, then $A$ is unique.  This follows from the fact about maximal antichains given in Lemma \ref{L:RegularOpenSets} together with completeness of the embedding $P \subset_c Q$.

\begin{Lemma}\label{L:SubposetsGiveSubalgebras} Suppose $P$ and $Q$ are posets with $P \subset_c Q$.  Suppose further that the maps $e_P$ and $e_Q$ and the embedding of Boolean algebras $j$ are defined as in Lemma \ref{L:POEmbeddingsInduceBAEmbeddings} and $ro^*(Q)$ is given by Definition \ref{D:ro*(Q)}.  Then $j(A)^*=A$ for all $A\in ro(P)$ and so $ro(P) \subset_c ro(Q)$.  Furthermore, the following diagram commutes in all squares, the maps forming the bottom square are complete embeddings of complete Boolean algebras, and the remaining maps are complete embeddings of posets:

 \[
 \begin{CD}
P @>\subset_c>> Q \\
@V{e_P}VV @V{e_Q}VV\\
{ro(P)} @> j>> {ro(Q)} \\
@V{=}VV @V{*}VV\\
{ro(P)} @> \subset_c >> {ro^*(Q)} \\
\\
 \end{CD}
 \]
\end{Lemma}

In the case of a class partial order $\P$ that is a chain of complete subposets, the idea above can be extended to build a corresponding chain of complete Boolean algebras.

\clearpage
\begin{Def}\rm\label{D:ChainOfCompleteSubalgebras} Suppose $\P = \bigcup_{\alpha\in\ORD} P_\alpha$ is a chain of complete subposets.  Then the associated \emph{chain of complete subalgebras} $\<ro^*(P_\alpha) \mid \alpha\in \ORD>$ is defined as follows.  For each $P_\alpha$, let $ro^*(P_\alpha)$ be $\{C^* \mid C\in ro(P_\alpha)\}$, where $C^*$ is defined $C^*=A$ such that 
\begin{enumerate}
\parskip=0pt
\item $\exists \beta \leq \alpha$, $A\in ro(P_\beta)$,
\item $\overline{A}^\circ_{P_\alpha} = C$, and
\item $\beta$ is least for which such an $A$ exists.
\end{enumerate}
For simplicity of notation, $ro^*(P_\alpha)$ will henceforth be denoted $B_\alpha$.
\end{Def}

\begin{Lemma}\label{L:ChainOfCompleteSubalgebras} Suppose $\P = \bigcup_{\alpha\in\ORD} P_\alpha$ is a chain of complete subposets and $\<B_\alpha \mid \alpha\in \ORD>$ is the associated chain of complete subalgebras. 
\begin{enumerate}
\parskip=0pt

\item For all $\alpha$ the map $C \mapsto C^*$ is an isomorphism of complete Boolean algebras, $ro(P_\alpha) \cong B_\alpha$.

\item For all $\alpha$ and $\beta$, if $\alpha < \beta$, then $B_\alpha$ is a complete subalgebra of  $B_\beta$, that is, $B_\alpha\subset_c B_\beta$.

\item The union $\B = \bigcup_{\alpha\in\ORD} B_\alpha$ is a Boolean algebra.

\item  \label{L:ChainOfCompleteSubalgebras:SetComplete} While $\B$ is not complete it is \emph{set complete} or $\ORD$\emph{-complete}, that is, $\sum A$ and $\prod A$ exists for any set $A\subset\B$.

\item If $\alpha < \beta$, then the diagram below commutes in all squares and every map is a complete embedding, where the partial order embedding $e_\beta$ and $e_\alpha$ and the Boolean algebra embedding $j$ are defined as in Lemma \ref{L:POEmbeddingsInduceBAEmbeddings}:

 \[
 \begin{CD}
P_\alpha @>\subset_c>> P_\beta \\
@V{e_\alpha}VV @V{e_\beta}VV\\
{ro(P_\alpha)} @> j >> {ro(P_\beta)} \\
@V{*}VV @V{*}VV\\
{B_\alpha} @> \subset_c >> {B_\beta} \\
\\
 \end{CD}
 \]
 
\item For all $\alpha$, the algebra $B_\alpha$ is a complete subalgebra of $\B$.

\item The map $e:\P \rightarrow \B$ given by the $e^*_\alpha$ is a dense embedding of partial orders.

\end{enumerate}
\end{Lemma}

The proof, while tedious, is essentially a matter of checking definitions.  Set completeness of $\B$ follows from the fact that any set $A\subset \B$ is contained in some $B_\alpha$, and we can apply completeness in $B_\alpha$ to obtain the sum or product.  The penultimate assertion is the analogue of Lemma \ref{L:EachSubposetIsCompleteInP}.  Note that the sequence $\<B_\alpha \mid \alpha\in\ORD>$ together with the associated Boolean operations on each $B_\alpha$ forms a class.

\section{Boolean-valued models}\label{SS:BVMs}

A Boolean algebra $\B$ can be thought of as extending our usual notion of true and false, represented by the two-element Boolean algebra $\{0,1\}$, to a more general setting, with `truth values' lying in $\B$.  A Boolean-valued model is a natural application of this idea in the context of sets.  Boolean-valued models are defined only in terms of complete Boolean algebras, and thus will not be applicable in the case of class forcing.  However, it will still be possible to construct a limited Boolean-valued model for class forcing, and to define the Boolean values of certain formulas, including all $\Delta_0$ formulas involving only set variables. This will help to provide the base case for the induction in the proof of the Forcing Theorem (Theorem \ref{T:ForcingTheorem}).

\begin{Def}\rm\label{D:BooleanValuedModels} Suppose $\B$ is a complete Boolean algebra.  A \emph{Boolean-valued model} $\mathfrak{W}=\<W,\W,\bval{{=}},\bval{{\in}}>$  consists of a universe $\W$ and a subcollection $W\subset \W$, and binary functions $\bval{=}$ and $\bval{\in}$ on $\W\times \W$ taking values in $\B$ (the \emph{Boolean values} of $=$ and $\in$), satisfying:
\begin{enumerate}
\parskip=0pt
\item $\bval{\class{X}=\class{X}}=1$ (reflexive),
\item $\bval{\class{X}=\class{Y}}=\bval{\class{Y}=\class{X}}$ (symmetric),
\item $\bval{\class{X}=\class{Y}}\btimes\bval{\class{Y}=\class{Z}}\leq\bval{\class{X}=\class{Z}}$ (transitive),
\item $\bval{\class{X}\in\class{Y}}\btimes\bval{\class{X}=\class{Z}}\btimes\bval{\class{Y}=\class{W}}\leq \bval{\class{Z}\in\class{W}}$.
\end{enumerate}
\end{Def}

Informally, $\W$ and $W$ are the Boolean analogues of the classes and sets $\V$ and $V$, respectively.  I will take the underlying universe $\W$ itself to be a class, rather than a collection of classes, as this will ease certain metamathematical difficulties in the exposition.  Note that in the case of class forcing and the corresponding model $\<V^\B, \V^\B>$ (Section \ref{SS:TheBooleanValuedModelV^B}), the functions $\bval{=}$ and ${\bval{\in}}$ will be defined only on sets, although the four conditions given will be satisfied wherever they are defined.  I would next like to define the Boolean value of any set-theoretic assertion about members of $\mathfrak{W}$.  I will describe the recursion below with special attention to the case where the definition fails if $\B$ is not complete.

\begin{Def}\rm\label{D:BooleanValuesOfAllFlas}For any formula $\phi(\class{X}_1, ... \class{X}_n)$ and any members $\class{A}_1, ... \class{A}_n \in \W$, the \emph{Boolean value} $\bval{\phi(\class{A}_1, ... \class{A}_n)}$ is defined recursively on the complexity of $\phi$ as follows.
\begin{enumerate}
\parskip=0pt
\item For atomic formulas the Boolean values are provided by Definition \ref{D:BooleanValuedModels}.
\item $\bval{\neg\phi(\class{A}_1, ... \class{A}_n)} = -\bval{\phi(\class{A}_1, ... \class{A}_n)}$.
\item $\bval{\(\phi\wedge\psi)(\class{A}_1, ... \class{A}_n)} = \bval{\phi(\class{A}_1, .. \class{A}_n)}\btimes\bval{\psi(\class{A}_1, ... \class{A}_n)}$.
\item \label{D:BValExistential}$\bval{\exists x \phi(x, \class{A}_1, ... \class{A}_n)} = \sum_{\class{A}\in \W} \bval{\phi(\class{A},\class{A}_1, ... \class{A}_n)}$.
\end{enumerate}
\end{Def}

If $\B$ is not complete, then the sum in step \ref{D:BValExistential} may not exist, and the induction will halt.  Thus for complete Boolean algebras we obtain the Boolean value of every statement, but for incomplete Boolean algebras we are guaranteed a Boolean value for $\Delta_0$ formulas only.

Note that Boolean values of other logical connectives are defined from those above in the usual fashion, e.g. $\bval{\phi\rightarrow\psi}=\bval{\neg\psi\vee\phi}=-\bval{\phi}+\bval{\psi}$. A formula $\phi$ is $\emph{valid}$ in $\mathfrak{W}$ if $\bval{\phi}=1$.  The first three conditions of Definition \ref{D:BooleanValuedModels} guarantee that the axioms for equality are valid in $\mathfrak{W}$.  Note that if $\bval{\phi}$ and $\bval{\psi}$ exist, then the implication $\phi\rightarrow\psi$ is valid if and only if $\bval{\phi}\leq\bval{\psi}$.  For example, this holds for every implication if the Boolean algebra is complete.  Furthermore, it is a matter of carefully checking definitions to verify that the axioms of predicate calculus are valid in $\mathfrak{W}$, provided their Boolean values exist, and that validity is preserved under the rules of inference.

This last fact provides a new approach to consistency proofs, and it is this approach that gives Boolean algebras their appeal in forcing arguments.  If $\mathfrak{W}$ a Boolean-valued model, and if the $\BGC$ axioms are valid in $\mathfrak{W}$, then for any assertion $\phi$ with $\bval{\phi}>0$ we have $\phi$ is consistent with $\BGC$, since $\bval{\phi}\neq 0$ implies $\bval{\neg\phi}\neq 1$ and so $\neg\phi$ cannot be proved from $\BGC$.

The remaining material in this section describes the method for obtaining a standard $2$-valued model from an arbitrary Boolean-valued model.  The main work here is in finding a way to reduce the Boolean algebra $\B$ to a $2$-element Boolean algebra, which is accomplished through the use of an ultrafilter (Definition \ref{D:BAfilters}).  Once again, completeness of the Boolean algebra together with additional restrictions on $\mathfrak{W}$ (see \emph{fullness}, Definition \ref{D:BVMfull}) will be required to give a `nice' characterization of satisfaction in the resulting model.  However, a  $2$-valued model can be obtained for any Boolean algebra, with satisfaction characterized for $\Delta_0$ formulas only.

Given a Boolean-valued model $\mathfrak{W}=\<W,\W,\bval{=},\bval{\in}>$ for a Boolean algebra $\B$ and an ultrafilter $\class{G}\subset \B$, I define an equivalence relation $\sim$ on $\W$ by $$\class{X}\sim \class{Y} \mbox{ if and only if } \bval{\class{X}=\class{Y}}\in \class{G}.$$
I now consider the quotient $\W/\sim$.  Unfortunately, it may be the case that a $\sim$ equivalence class may a proper class, which leads to difficulties in the definition of $\W/\sim$ as classes can contain only sets as members.  However, we can use a standard trick (Scott's Trick) from the construction of ultrapowers to reduce each equivalence class to a set.  For $\class{A}\in \W$, let $\alpha$ be least such there is $\class{C} \in \W$ of rank $\alpha$ with $\class{A} \sim \class{C}$.  Let $[\class{A}] = \{\class{C}\in \W \mid \class{A}\sim \class{C} $ and $\class{C}$ has rank $\alpha \}$.  While it is no longer necessarily true that $\class{A}\in [\class{A}]$, the definition of $\W/\sim$ can be carried out without difficulty.  Next define a binary relation $\class{E}$ on $\W/\sim$ by $$[\class{X}] \, \class{E} \, [\class{Y}] \mbox{ if and only if } \bval{\class{X}=\class{Y}}\in \class{G}.$$
It is left to the reader to verify that $\sim$ is an equivalence relation on $\W$ and that the value of $\class{E}$ does not depend on the choice of representatives $\class{X}$ and $\class{Y}$, both consequences of conditions 1-4 in Definition \ref{D:BooleanValuedModels}.  The relation $\sim$ defines the equality relation for $\mathfrak{W}$, and $\class{E}$ defines the epsilon relation of $\W/\sim$.  Satisfaction is characterized in $\(W/\sim,\class{E})$ as follows:

\begin{Lemma}\label{L:BVMSatisfactionDelta0} Suppose $\mathfrak{W}$ is a Boolean-valued model for a complete Boolean algebra $\B$, and $\class{G}\subset \B$ is an ultrafilter.  For any $\phi(\class{X}_1, ... \class{X}_n)\in \Delta_0$ and any $\class{A}_1, ... \class{A}_n\in \W$, $$\(\W/\sim,\class{E}) \vDash \phi\([\class{A}_1],...,[\class{A}_n]) \mbox{ if and only if } \bval{\phi\(\class{A}_1,...,\class{A}_n)}\in \class{G}$$
\end{Lemma}

\begin{proof}   By induction on complexity of $\phi$.  The proof for atomic $\phi$ was given in the preceding discussion.  For negations, recall that $\bval{\neg\phi}=-\bval{\phi}$ (the parameters $\class{A}_1, ... \class{A}_n$ are suppressed, as they do not affect the proof).  Thus 
\begin{eqnarray*}
	\bval{\neg\phi}\in \class{G} & \iff & -\bval{\phi} \notin \class{G} \mbox{ ($\class{G}$ is an ultrafilter)}\\
	& \iff &  \(W/\sim,\class{E}) \nvDash \phi   \mbox{ (inductive assumption)}\\
	& \iff &  \(W/\sim,\class{E}) \vDash \neg\phi \mbox{ (definition of satisfaction)}
\end{eqnarray*}
 For conjunctions, 
\begin{eqnarray*}
	\lefteqn{\bval{\phi\vee\psi}\in \class{G}}\\
	& \iff & \bval{\phi}\btimes\bval{\psi} \in \class{G} \mbox{ (definition of $\bval{}$ for conjunctions)}\\
	& \iff &  \bval{\phi} \in \class{G} \mbox{ and } \bval{\psi} \in \class{G}   \mbox{ (since $\class{G}$ is a filter)}\\
	& \iff &  \(W/\sim,\class{E}) \vDash \phi  \mbox{ and } \(W/\sim,\class{E}) \vDash \psi \mbox{ (inductive assumption)}\\
	& \iff &  \(W/\sim,\class{E}) \vDash \phi\wedge\psi \mbox{ (definition of satisfaction)}
\end{eqnarray*}

The second line above may require some explanation.  If $\class{G}$ is a filter on $\B$ and $u\btimes v\in \class{G}$, then it is clear that $u\in \class{G}$ and $v\in \class{G}$, since $u$ and $v$ are $\geq u\btimes v$ and $\class{G}$ is closed upwards.  For the reverse implication, suppose $u$ and $v$ are in $\class{G}$.  As $\class{G}$ is a filter, there is $r\in \class{G}$ such that $r$ lies below both $u$ and $v$.  However, $u\btimes v$ is the greatest lower bound of $\{u,v\}$, so $r\leq u\btimes v$.  Thus $u\btimes v \in \class{G}$, once again by upwards closure of the filter. \end{proof}

I would like to extend Lemma \ref{L:BVMSatisfactionDelta0} to all formulas in $\(\W/\sim, \class{E})$ and thus obtain $2$-valued model with satisfaction entirely characterized by $\bval{\phi}\in \class{G}$.  However, this requires an additional condition, that the Boolean-valued model be $\emph{full}$.  While this condition will not hold in the case of forcing with proper class partial orders, the methods demonstrated can still be of use.  

\clearpage
\begin{Def}\rm\label{D:BVMfull} If $\mathfrak{W}$ is a Boolean-valued model for a complete Boolean algebra $\B$, then $\mathfrak{W}$ is \emph{full} if and only if for every formula $\phi(\class{X}, \class{X}_1, ... \class{X}_n)$ and every $\class{A}_1, ... \class{A}_n \in \W$ there is $\class{A} \in \W$ such that 
$$\bval{\exists x \phi(x,\class{A}_1, ... \class{A}_n)}=\bval{\phi(\class{A},\class{A}_1, ... \class{A}_n)}$$
\end{Def}

For a full Boolean-valued model, Lemma \ref{L:BVMSatisfactionDelta0} holds generally.

\begin{Lemma}\label{L:BVMSatisfaction} Suppose $\mathfrak{W}$ is a full Boolean-valued model for a complete Boolean algebra $\B$, and $\class{G}\subset \B$ is an ultrafilter.  For any formula $\phi(\class{X}_1, ... \class{X}_n)$ and any $\class{A}_1, ... \class{A}_n\in \W$, $$\(\W/\sim,\class{E}) \vDash \phi\([\class{A}_1],...,[\class{A}_n]) \mbox{ if and only if } \bval{\phi\(\class{A}_1,...,\class{A}_n)}\in \class{G}$$
\end{Lemma}

\begin{proof} By induction on complexity of $\phi$.  Lemma \ref{L:BVMSatisfactionDelta0} provides the proof for atomic formulas, negations, and conjunctions, and so it suffices to consider existential formulas.  Consider $\exists \class{X}\phi(\class{X}, \class{A}_1, ... \class{A}_n)$. 
\begin{eqnarray*}
	\lefteqn{\bval{\exists \class{X}\phi(\class{X}, \class{A}_1, ... \class{A}_n)}\in \class{G} }\\
	& \iff & \exists \class{A}\in \W \, (\bval{\phi(\class{A}, \class{A}_1, ... \class{A}_n)} \in \class{G}) \mbox{ (fullness)}\\
	& \iff &  \exists \class{A}\in \W \, \({ \(W/\sim,\class{E}) \vDash \phi([\class{A}], [\class{A}_1],...,[\class{A}_n]) })   \mbox{ (inductive assumption)}\\
	& \iff &  \(W/\sim,\class{E}) \vDash \exists \class{X} \phi(\class{X}, [\class{A}_1],...,[\class{A}_n]) \mbox{ (definition of satisfaction)}
\end{eqnarray*}

Thus satisfaction in the model $(\W/\sim,\class{E})$ holds exactly for formulas with Boolean value in $\class{G}$. \end{proof}

\section{The  model $\V^\B$}\label{SS:TheBooleanValuedModelV^B}

Given a model $\<V,\V>$ of $\BGC$, I will now define the model $\<V^\B,\V^\B>$ associated with a Boolean algebra $\B$.  While the model will not be a Boolean-valued model as in Definition \ref{D:BooleanValuedModels}, it will satisfy a sufficient portion of the definition to be of use in later forcing arguments.  I begin with the case where where $B$ is a set.  The generalization to proper class $\B$ will follow, with certain restrictions on $\B$.  The definition is intended to mirror that of the powerset sequence $V_\alpha$.  However, instead of taking $2$-valued subsets at each stage I consider $B$-valued subsets.  A $B$-valued subset of a set $A$ is a function $\tau:A\rightarrow B$, where $\tau(a)$ is understood to be something like the Boolean value of `$a\in \tau$.'  The relationship between $\tau(a)$ and $\bval{a\in \tau}$ is in fact slightly more subtle, and is stated precisely in Lemma \ref{L:V^BisBVM} assertion \ref{L:V^BisBVM:tau(x)}.

\begin{Def}\rm\label{D:V^Bsets}If $\<V,\V>$ is a model of $\BGC$ and the set $B \in V$ is a Boolean algebra, then $V^B$ is defined $V^B = \bigcup_{\alpha\in\ORD} V^B_\alpha$, where
\begin{enumerate}
\parskip=0pt
\item $V^B_0=\emptyset$,
\item $V^B_{\alpha+1} = \{\tau \in V \mid \tau$ is a partial function from $V^B_\alpha$ to $B\}$, and
\item $V^B_\lambda = \bigcup_{\alpha<\lambda} V^B_\alpha$ for $\lambda$ a limit,
\end{enumerate}
and $\V^B$ is defined to be the collection of classes $\class{A}\in \V$ such that $\class{A} \subset V^B \times B$.  
\end{Def}

Members of $\V^B$ and $V^B$ are called \emph{$B$-names}.  It might seem more natural to define $\V^B$ to be all \emph{functions} $\class{A}\from V^B \rightarrow B$, and in the case of $B$ a complete Boolean algebra this definition is equivalent.  In such a case, for $\class{A} \subset V^B \times B$ one can associate with each $a\in \dom(\class{A})$ the single Boolean value $\sum_{\<a,u>\in \class{A}} u$, thus producing a function $\class{A}^\prime\from V^B \rightarrow B$ which will be equivalent to $\class{A}$ in all forcing arguments.  However, for incomplete Boolean algebras the indicated sum may not exist and so the more general definition is necessary.  Observe that $V^B \subset \V^B$.  Also note that if $B$ is a set, then every $V^B_\alpha$ is a set.  For proper class $\B$, however, $V^\B_2$ is already a proper class and, since induction is defined only for sets, Definition \ref{D:V^Bsets} fails.  However, a direct definition is possible.

\begin{Def}\rm\label{D:V^Bclasses}If $\<V,\V>$ is a model of $\BGC$ and $\B$ a Boolean algebra, then $V^\B$ is the collection of all $x\in V$ such that there is a sequence $\<x_\alpha \mid \alpha\leq\delta> \in V$ satisfying 
\begin{enumerate}
\parskip=0pt
\item $x_0=\emptyset$,
\item for $\alpha<\delta$, $\forall y\in x_{\alpha+1}$, $y$ is a partial function from $x_\alpha$ to $\B$,
\item for $\lambda \leq \delta$ a limit, $x_\lambda \subset \bigcup_{\alpha<\lambda} x_\alpha$, and
\item $x_\delta=x$,
\end{enumerate}
and $\V^\B$ is defined to be the collection of all $\class{A} \subset V^\B \times \B$.
\end{Def}

Observe once again that $V^\B \subset \V^\B$.  That Definition \ref{D:V^Bclasses} agrees with Definition \ref{D:V^Bsets} in the case of $\B$ a set is left as an exercise.  However, a more useful definition can be obtained in the case where $\B=\bigcup B_\alpha$ is a chain of complete subalgebras associated with $\P=\bigcup P_\alpha$ a chain of complete subposets.  In this case, $V^\B$ is simply the union of the $V^{B_\alpha}$.  However, it will be useful for induction purposes to write $V^\B$ as an increasing union of sets as follows.

\begin{Lemma}\label{L:V^Bclasses}If $\<V,\V>$ is a model of $\BGC$ and $\P=\bigcup P_\alpha$ is a chain of complete subposets with associated chain of complete subalgebras $\B=\bigcup B_\alpha$, then $x\in V^\B$ if and only if $x \in V^{B_\alpha}$ for some $\alpha$.  Thus $$V^\B = \bigcup_{\gamma\in\ORD} \({\bigcup_{\alpha,\beta \leq \gamma} V^{B_\alpha}_\beta}),$$
where $V^{B_\alpha}_\beta$ is the $\beta^{th}$ stage of the construction given in Definition \ref{D:V^Bsets} for the set Boolean algebra $B_\alpha$.
\end{Lemma}

\begin{proof}   Note that for an set $A\subset \B$, if $\alpha$ is sufficiently large, then $A$ will be contained in $B_\alpha$.  Thus for any $x\in V^\B$, the range of $x$ is contained in some $B_\alpha$.  Of course, this does not guarantee that $x\in V^{B_\alpha}$.  For example, there may be $y\in \dom(x)$ which mentions conditions in some later $B_\beta$.  Thus we must consider not only the range of $x$, but the range of $y$ for each $y\in \dom(x)$, and the range of $z$ for $z\in \dom(y)$, etc..  A moments thought will verify that the transitive closure $\tc(x)$ contains all the pertinent information, and as $\tc(x)$ is a set, $\tc(x)\cap \B \subset B_\alpha$ for some $\alpha$. Thus $x\in V^{B_\alpha}$. \end{proof}

Definitions \ref{D:V^Bsets} and \ref{D:V^Bclasses} give us the underlying universe $\<V^\B,\V^\B>$ of our Boolean-valued model.  It remains to equip this universe with Boolean-valued analogues of the relations $=$ and $\in$.  These relations, denoted $\bval{x=y}$ and $\bval{x\in y}$, will be functions on $V^\B \times V^\B$ taking values in $\B$.  Unfortunately, extending these functions to $\V^\B \times \V^\B$, the `classes' of our model, requires a complete Boolean algebra.  However, for the purposes of forcing it will be sufficient to have a definition in the set case.  The definition is recursive, and some care is required to carry out the recursion rigorously.  The recursion is based on two notions, the \emph{canonical well-ordering of pairs of ordinals} and the \emph{$\B$-rank} of members of $V^\B$.

\begin{Def}\rm\label{D:CanonicalWellOrderingofPairs}Define an ordering on pairs of ordinals by $\<\alpha,\beta> \prec \<\xi,\eta>$ if and only if
\begin{enumerate}
\parskip=0pt
\item $\max\{\alpha,\beta\} < \max\{\xi,\eta\}$, or
\item $\max\{\alpha,\beta\} = \max\{\xi,\eta\}$ and $\alpha < \xi$, or
\item $\max\{\alpha,\beta\} = \max\{\xi,\eta\}$ and $\alpha = \xi$ and $\beta < \eta$.
\end{enumerate}
\end{Def}

The ordering $\prec$ is a well ordering, the \emph{canonical well ordering of pairs of ordinals}.  If $\<\alpha,\beta>$ is the $\gamma^{th}$ pair under $\prec$, then the map taking $\<\alpha,\beta>$ to $\gamma$ is taken to be `understood' and $\<\alpha,\beta>$ is simply identified with $\gamma$, that is, $\<\alpha,\beta>=\gamma$.  Furthermore, in practice the $\prec$ is replaced with the standard $<$, e.g. $\<\alpha,\beta> < \<\xi,\eta>$.

The \emph{canonical well ordering of triples} of ordinals is a natural extension of the ordering on pairs, in which triples $\<\alpha,\beta,\gamma>$ are ordered first according maximum, then by first coordinate, then second coordinate, and then third coordinate.

Note that if $\<\alpha,\beta>=\gamma$, then $\alpha$ and $\beta$ are both $\leq \gamma$.

\begin{Def}\rm\label{D:V^BRecursion}For $\B$ and $V^\B$ as in Lemma \ref{L:V^Bclasses}, define for each $x\in V^\B$ the \emph{$\B$-rank} $\rho(x)$ as the least $\<\alpha,\beta>$ under the canonical well ordering such that $x \in V^{B_\alpha}_\beta$.  In the case of Definition \ref{D:V^Bsets}, where $B$ is a set, $\rho(x)$ is taken to the be least $\alpha$ such that $x\in V^B_\alpha$.
\end{Def}

Recursion according to $\rho(x)$ is similar to $\in$-recursion, and depends on the following basic fact.

\begin{Lemma}\label{L:BRankBasicFact} If $x\in V^\B$ and $y\in \dom(x)$, then $\rho(y) < \rho(x)$.  
\end{Lemma}

\begin{proof}   Suppose $\rho(x)=\<\alpha,\beta>$ and $y\in \dom(x)$ with $\rho(y)=\<\xi,\eta>$.  Clearly $y\in V^{B_\alpha}_\beta$ since $x$ is, so $\xi \leq \alpha$ and $\eta \leq \beta$.  Suppose $\max\{\alpha,\beta\} = \max\{\xi,\eta\}$ and $\alpha = \xi$.  Then both $x$ and $y$ are in $V^{B_\alpha}$ but neither appear in $V^{B_{\alpha^\prime}}$ for any $\alpha^\prime < \alpha$.  However, it follows from the recursive definition of the $V^{B_\alpha}$ that $y \in \dom(x)$ implies that $y$ was created at an earlier stage of the $V^{B_\alpha}$ construction, and so $\beta < \eta$.  Thus $\rho(y) < \rho(x)$. \end{proof}

\begin{Def}\rm\label{D:V^Batomic}For $\B$, $V^\B$ as in Lemma \ref{L:V^Bclasses}, define binary functions ${\bval{\tau \in \sigma}}$, ${\bval{\tau \subset \sigma}},$ and $\bval{\tau = \sigma}$ recursively as follows.
\begin{enumerate}
\parskip=0pt
\item $\bval{\tau\in\sigma}=\sum_{x\in \dom(\sigma)} \bval{\tau = x}\btimes \sigma(x)$.
\item $\bval{\tau \subset \sigma} = \prod_{y\in \dom(\tau)} (\tau(y)\rightarrow\bval{y\in\sigma})$ (recall the Boolean operation $u\rightarrow v$ is  defined $-u+v$).
\item $\bval{\tau=\sigma}=\bval{\tau\subset\sigma}\btimes\bval{\sigma\subset\tau}$.
\end{enumerate}
The recursion is carried out on pairs  $\<\tau,\sigma>$ according to the canonical well ordering of their $\B$-ranks $\<\rho(\tau),\rho(\sigma)>$.
\end{Def}

Notice that all sums and products are taken over sets and thus exist, by set-completeness of $\B$ (Lemma \ref{L:ChainOfCompleteSubalgebras} conclusion \ref{L:ChainOfCompleteSubalgebras:SetComplete}).  It is for exactly this reason that the definition cannot in general be extended to proper class $\P$-names.  However, it is still possible to define the Boolean value of certain additional formulas according to Definition \ref{D:BooleanValuesOfAllFlas}.  For example, the Boolean value of every $\Delta_0$ formula involving only sets will exists in $\B$.  

To see that the recursion in the definition above is well-defined, consider a pair $\<\tau,\sigma>$.  The definition of $\bval{\tau\in\sigma}$ depends on $\bval{\tau= x}$ for $x\in \dom(\sigma)$. Lemma \ref{L:BRankBasicFact} gives $\rho(x) < \rho(\sigma)$, and so $\<\rho(\tau), \rho(x)> < \<\rho(\tau),\rho(\sigma)>$.  Similarly, $\bval{\tau \subset \sigma}$ depends on $\bval{y\in\sigma}$ for $y\in \dom(\tau)$, and  $\<\rho(y),\rho(\sigma)> < \<\rho(\tau),\rho(\sigma)>$.  Finally, $\bval{\tau=\sigma}$ depends on $\bval{\tau\subset\sigma}$ and $\bval{\sigma \subset \tau}$.  The former depends only on earlier stages of the recursion, as argued above, but the latter requires a little more argument.  The value $\bval{\sigma \subset \tau}$ depends on $\bval{y\in \tau}$ for $y\in \dom(\sigma)$, and I must argue that $\<\rho(y),\rho(\tau)> < \<\rho(\tau),\rho(\sigma)>$.  Clearly $\max\{\rho(y),\rho(\tau)\} \leq \max\{\rho(\tau),\rho(\sigma)\}$ since $\rho(y) < \rho(\sigma)$.  If the maximum of the former is strictly less, then I am done.  Otherwise, it must be the case that $\rho(\sigma)\leq\rho(\tau)$.  But $\rho(y) < \rho(\sigma)$, so $\rho(y) < \rho(\tau)$, and so $\<\rho(y),\rho(\tau)> < \<\rho(\tau),\rho(\sigma)>$ according to the canonical well ordering of pairs.  Thus the recursion is well-defined.

What is more, the nature of the recursion is such that $\bval{\tau\in\sigma}$, etc.,  can be calculated correctly in $B_\gamma$, provided $\tau$ and $\sigma$ are in $V^{B_\gamma}$.

\begin{Lemma}\label{L:BVunambiguous}  For all $\gamma$, for all $\tau$ and $\sigma$ in $V^{B_\gamma}$, 
\begin{enumerate}
\parskip=0pt
\item $\bval{\tau\in\sigma}_\B=\bval{\tau\in\sigma}_{B_\gamma}$,
\item $\bval{\tau\subset\sigma}_\B=\bval{\tau\subset\sigma}_{B_\gamma}$,
\item $\bval{\tau=\sigma}_\B=\bval{\tau=\sigma}_{B_\gamma}$,
\end{enumerate}
where, for example, $\bval{\tau\in\sigma}_{B_\gamma}$ is the result of applying Definition \ref{D:V^Batomic} to the Boolean algebra $B_\gamma$ and the universe $V^{B_\gamma}$.
\end{Lemma}
\emph{Sketch of Proof}.  Careful examination of the definition shows that for $\tau$ and $\sigma$ in $V^{B_\gamma}$ the recursion relies only on Boolean values of statements involving other members of $V^{B_\gamma}$.  Furthermore, $B_\gamma \subset_c \B$ shows that the sums and products taken in $B_\gamma$ will agree with those taken in $\B$.  A rigorous proof of the assertion is a straightforward induction on pairs $\<\rho(\tau),\rho(\sigma)>$ for $\tau,\sigma \in V^{B_\gamma}$.

\begin{Lemma}\label{L:V^BisBVM}  The Boolean operations of Definition \ref{D:V^Batomic} satisfy
\begin{enumerate}
\parskip=0pt
\item \label{L:V^BisBVM:tau(x)} $\tau(x) \leq \bval{x\in\tau}$ for every $x\in \dom(\tau)$, 
\item \label{L:V^BisBVM:reflexivity}$\bval{\tau=\tau}=1$ (reflexive),
\item \label{L:V^BisBVM:symmetry}$\bval{\tau=\sigma}=\bval{\sigma=\tau}$ (symmetric),
\item \label{L:V^BisBVM:transitivity} $\bval{\tau=\sigma} \btimes \bval{\sigma=x} \leq \bval{\tau=x}$ (transitive), and
\item \label{L:V^BisBVM:membership}$\bval{\tau\in \sigma}\btimes\bval{\tau=x}\btimes\bval{\sigma=y}\leq \bval{x\in y}$.
\end{enumerate}
\end{Lemma}
In subsequent sections, an ultrafilter on $\B$ will be used to reduce the Boolean-valued statements to standard $2$-valued statements. This lemma ensures that various fundamental properties of equality and membership will hold in the resulting $2$-valued model. Assertion \ref{L:V^BisBVM:tau(x)} describes the basic relationship of $\B$-names to membership, assertions \ref{L:V^BisBVM:reflexivity}-\ref{L:V^BisBVM:transitivity} provide the axioms of equality for sets, and assertion \ref{L:V^BisBVM:membership} ensures that equality and membership interact correctly.

\begin{proof}[Proof of Lemma \ref{L:V^BisBVM}]
Assertions \ref{L:V^BisBVM:tau(x)} and \ref{L:V^BisBVM:reflexivity} are proved simultaneously by induction on $\rho(\tau)$.  For $x\in \dom(\tau)$,
$$\bval{x\in\tau}  =  \sum_{y\in \dom(\tau)} \bval{x = y}\btimes \tau(y)$$
Recalling that $\sum$ is given by the supremum, and observing that one such $y\in \dom(\tau)$ is $x$ itself, I obtain
$$\sum_{y\in \dom(\tau)} \bval{x = y}\btimes \tau(y)  \geq  \bval{x=x}\btimes\tau(x) =  1\btimes \tau(x) =  \tau(x)$$
using the inductive hypothesis to conclude $\bval{x=x}=1$.  For assertion \ref{L:V^BisBVM:reflexivity}, 
\begin{eqnarray*}
	\bval{\tau=\tau} & = &  \bval{\tau\subset\tau}\btimes\bval{\tau\subset\tau} \\
	& = & \bval{\tau\subset\tau} \\
	& = &  \prod_{y\in \dom(\tau)} \({ \tau(y)\rightarrow\bval{y\in\tau}})\\
	& = &  \prod_{y\in \dom(\tau)} \({-\tau(y)+\bval{y\in\tau}})\\
\end{eqnarray*}
The inductive hypothesis gives $\bval{y\in\tau}\geq\tau(y)$, and so $$-\tau(y)+\bval{y\in\tau} \geq -\tau(y)+\tau(y) = 1$$
and thus $\bval{\tau=\tau}=\prod_{y\in \dom(\tau)} 1 = 1$.

For assertion \ref{L:V^BisBVM:symmetry}, 
\begin{eqnarray*}
	\bval{\tau=\sigma} & = &  \bval{\tau\subset\sigma}\btimes\bval{\sigma\subset\tau} \\
	& = & \bval{\sigma\subset\tau}\btimes\bval{\tau\subset\sigma} \\
	& = & \bval{\sigma=\tau}
\end{eqnarray*}

Assertion \ref{L:V^BisBVM:transitivity} requires the most complicated induction.  A number of supporting assertions will be proven simultaneously, each necessary to allow the induction to proceed.  The induction will be performed on triples $\<\tau,\sigma,x>$ according to their $\B$-ranks $\<\rho(\tau),\rho(\sigma),\rho(x)>$, ordered by the canonical well ordering of triples (Definition \ref{D:CanonicalWellOrderingofPairs}).  The inductive proof will be carried out simultaneously for the following assertions:
\begin{eqnarray}
\label{Eq:Bsubset1} \bval{\tau\subset\sigma}\btimes\bval{\sigma=x}\leq\bval{\tau\subset x}\\
\label{Eq:Bsubset2} \bval{\sigma\subset\tau}\btimes\bval{\sigma=x}\leq\bval{x\subset\tau}\\
\label{Eq:Bmembership1} \bval{\sigma\in\tau}\btimes\bval{\sigma=x}\leq\bval{x\in\tau}\\
\label{Eq:Bmembership2} \bval{\sigma\in x}\btimes\bval{\tau=\sigma} \leq \bval{\tau\in x}\\
\label{Eq:Bmembership3} \bval{\tau\in\sigma}\btimes\bval{\sigma=x}\leq\bval{\tau\in x}\\
\label{Eq:Btransitivity} \bval{\tau=\sigma} \btimes \bval{\sigma=x} \leq \bval{\tau=x}
\end{eqnarray}
Note that Lemma \ref{L:V^BisBVM} assertion \ref{L:V^BisBVM:transitivity} is given by (\ref{Eq:Btransitivity}). Suppose the assertions hold for all triples $\<\tau^\prime,\sigma^\prime,x^\prime>$ such that $\<\rho(\tau^\prime),\rho(\sigma^\prime),\rho(x^\prime)> < \<\rho(\tau),\rho(\sigma),\rho(x)>$.  I will show that they must hold for $\<\tau,\sigma,x>$.

Assertion (\ref{Eq:Bsubset1}) can be rewritten as follows, applying the definition of $\bval{\subset}$ and the appropriate rules of Boolean algebras at each step.
\begin{eqnarray*}
	\bval{\tau\subset\sigma}\btimes\bval{\sigma=x} & \leq &\bval{\tau\subset x} \\
	\iff \({\prod_{y\in \dom(\tau)} \tau(y)\rightarrow\bval{y\in\sigma}})\btimes\bval{\sigma=x} & \leq & \prod_{y\in \dom(\tau)} \({\tau(y)\rightarrow\bval{y\in x}})\\
	\iff \prod_{y\in \dom(\tau)} \({\({-\tau(y)+\bval{y\in\sigma}}) \btimes\bval{\sigma=x}}) & \leq & \prod_{y\in \dom(\tau)} \({-\tau(y)+\bval{y\in x}})
\end{eqnarray*}
Thus it suffices to verify for each $y\in \dom(\tau)$ that 
\begin{eqnarray*}
	\({-\tau(y)+\bval{y\in\sigma}}) \btimes\bval{\sigma=x}  &\leq&   -\tau(y)+\bval{y\in x}\\
	\iff -\tau(y)\btimes\bval{\sigma=x}+\bval{y\in\sigma} \btimes\bval{\sigma=x}  &\leq&   -\tau(y)+\bval{y\in x}
\end{eqnarray*}
Since $-\tau(y)\btimes\bval{\sigma=x} \leq -\tau(y)$, this reduces to showing  
$$ \bval{y\in\sigma} \btimes\bval{\sigma=x}  \leq   \bval{y\in x}.$$
This is exactly the inductive hypothesis (\ref{Eq:Bmembership3}) applied to the triple $\<y,\sigma,x>$. Since $y\in \dom(\tau)$ it follows that $\<\rho(y),\rho(\sigma), \rho(x)> < \<\rho(\tau),\rho(\sigma), \rho(x)>$, and so the inductive hypothesis applies.  This proves (\ref{Eq:Bsubset1}).

Assertion (\ref{Eq:Bsubset2}) requires a bit more work.  
$$\bval{\sigma\subset\tau}\btimes\bval{\sigma=x} \leq \bval{x\subset\tau}$$
Rewriting the left side yields
\begin{eqnarray*}
	\lefteqn{\bval{\sigma\subset\tau}\btimes\bval{\sigma=x}}\\
	&=& \bval{\sigma\subset\tau}\btimes\bval{\sigma\subset x}\btimes\bval{x\subset \sigma}\\
	&=& \bval{\sigma\subset\tau}\btimes\bval{\sigma\subset x} \({ \prod_{y\in \dom(x)} -x(y) + \bval{y\in\sigma} }) \\
	&=&   \prod_{y\in \dom(x)} \({\bval{\sigma\subset\tau}\btimes\bval{\sigma\subset x} \btimes \({ -x(y) + \bval{y\in\sigma} }) })
\end{eqnarray*}
And the right side gives 
$$\bval{x\subset\tau} = \prod_{y\in \dom(x)} \({-x(y) + \bval{y\in\tau}})$$

It once again is enough to show inequality of the expressions inside the products for each $y\in \dom(x)$.  The terms involving $-x(y)$ trivially satisfy $\leq$, reducing the problem to showing for each $y\in \dom(x)$
$$\bval{y\in\sigma}\btimes\bval{\sigma\subset\tau}\btimes\bval{\sigma\subset x} \leq \bval{y\in\tau}.$$
In fact, it is enough to show 
$$\bval{y\in\sigma}\btimes\bval{\sigma\subset\tau} \leq \bval{y\in\tau}$$
as the additional factor $\bval{\sigma\subset x}$ will only make the left side smaller.  I now expand the left side using the definitions of $\bval{\sigma\subset\tau}$ and $\bval{y\in\sigma}$ and the basic distributive property of Boolean algebras.
\begin{eqnarray*}
	\lefteqn{\bval{y\in\sigma}\btimes\bval{\sigma\subset\tau} }\\
	 &=& \({\sum_{z\in \dom(\sigma)}  \bval{y=z}\btimes\sigma(z)})\btimes  \bval{\sigma\subset\tau}\\	
	  &=&  \sum_{z\in \dom(\sigma)} \({ \bval{y=z}\btimes\sigma(z)\btimes  \bval{\sigma\subset\tau} })\\
	  &=& \sum_{z\in \dom(\sigma)} \({ \bval{y=z}\btimes\sigma(z)\btimes  \({ \prod_{z^\prime \in \dom(\sigma)} -\sigma(z^\prime)+\bval{z^\prime \in \tau}   }) })\\
	  &=& \sum_{z\in \dom(\sigma)} \prod_{z^\prime \in \dom(\sigma)}   \({     \bval{y=z}\btimes\sigma(z)\btimes\({-\sigma(z^\prime)}) +\bval{y=z}\btimes\sigma(z)\btimes\bval{z^\prime \in \tau} })
\end{eqnarray*}
I must show that this expression is $\leq \bval{y\in\tau}$.  It suffices to see that every member of the sum is $\leq \bval{y\in\tau}$, and since each summand is given by a product it suffices to show that for each summand there is a single member of the product which is $\leq \bval{y\in\tau}$.  Thus I must show that for every $z\in \dom(\sigma)$ there is $z^\prime \in \dom(\sigma)$ such that the expression inside the product above is $\leq \bval{y\in\tau}$.  In fact, I claim that for each $z$, taking $z=z^\prime$ gives the desired result.  For if $z=z^\prime$, I have 
\begin{eqnarray*}
	&&\bval{y=z}\btimes\sigma(z)\btimes\({-\sigma(z^\prime)}) +\bval{y=z}\btimes\sigma(z)\btimes\bval{z^\prime \in \tau}\\
	&=& \bval{y=z}\btimes 0 + \bval{z \in \tau}\btimes\bval{y=z}\btimes\sigma(z)\\
	& \leq & \bval{z \in \tau}\btimes\bval{y=z}
\end{eqnarray*}
Now, $\bval{z \in \tau}\btimes\bval{z=y} \leq \bval{y\in\tau}$ is simply an application of inductive hypothesis (\ref{Eq:Bmembership1}) to the triple $\<\tau,z,y>$.  Since $z\in \dom(\sigma)$ and $y\in \dom(x)$, $\<\rho(\tau),\rho(z), \rho(y)> < \<\rho(\tau),\rho(\sigma), \rho(x)>$ and the inductive hypothesis applies, proving (\ref{Eq:Bsubset2}).

Assertion (\ref{Eq:Bmembership1}) is much more straightforward.
\begin{eqnarray*}
	\bval{\sigma\in\tau}\btimes\bval{\sigma=x}&\leq&\bval{x\in\tau}\\
	\iff \sum_{y\in \dom(\tau)} \bval{\sigma=y}\btimes \tau(y)\btimes\bval{\sigma=x} & \leq & \sum_{y\in \dom(\tau)} \bval{x=y}\btimes \tau(y)
\end{eqnarray*}
It suffices to show, for each $y\in \dom(\tau)$, $$\bval{\sigma=y}\btimes\bval{\sigma=x} \leq \bval{y=x}$$  Applying inductive hypothesis (\ref{Eq:Btransitivity}) to $\<y,\tau,x>$ completes the proof of (\ref{Eq:Bmembership1}).

Assertion (\ref{Eq:Bmembership2}) is proved in exactly the manner of (\ref{Eq:Bmembership1}), applying commutativity of $\bval{=}$ where necessary.

The proof of assertion (\ref{Eq:Bmembership3}) is similar to that of (\ref{Eq:Bsubset2}).
\begin{eqnarray*}
	\bval{\tau\in\sigma}\btimes\bval{\sigma=x}&\leq&\bval{\tau\in x}\\
	\iff \bval{\tau\in\sigma}\btimes\bval{\sigma\subset x}\btimes\bval{x\subset\sigma}&\leq&\bval{\tau\in x}
\end{eqnarray*}
Thus it suffices to show $$ \bval{\tau\in\sigma}\btimes\bval{\sigma\subset x}\leq\bval{\tau\in x}$$
Expanding the left side and applying the distributive laws of Boolean algebras yields
\begin{eqnarray*}
	\lefteqn{ \bval{\tau\in\sigma}\btimes\bval{\sigma\subset x} }\\
	& = &
\sum_{y\in \dom(\sigma)} \prod_{y^\prime \in \dom(\sigma)} \({
\bval{y=\tau}\btimes\sigma(y)\btimes\({-\sigma(y^\prime)}) + \bval{y=\tau}\btimes\sigma(y)\btimes\bval{y^\prime\in x}
})
\end{eqnarray*}
Once again, I must show that for every $y\in \dom(\sigma)$ there is $y^\prime \in \dom(\sigma)$ such that the expression inside the product is $\leq\bval{\tau\in x}$.  Choosing $y^\prime = y$ gives
\begin{eqnarray*}
	&&\bval{y=\tau}\btimes\sigma(y)\btimes\({-\sigma(y)}) + \bval{y=\tau}\btimes\sigma(y)\btimes\bval{y\in x}\\
	&=& 0 + \bval{y=\tau}\btimes\sigma(y)\btimes\bval{y\in x}\\
	&\leq& \bval{y=\tau}\btimes\bval{y\in x}
\end{eqnarray*}
and $\bval{y\in x}\btimes\bval{y=\tau} \leq \bval{\tau\in x}$ is obtained by applying inductive assumption (\ref{Eq:Bmembership2}) to the triple $\<\tau,y,x>$.  This proves (\ref{Eq:Bmembership3}).

Finally, assertion (\ref{Eq:Btransitivity}) is an immediate consequence of (\ref{Eq:Bsubset1}) and (\ref{Eq:Bsubset2}).

This concludes the inductive proof of Lemma \ref{L:V^BisBVM} assertion \ref{L:V^BisBVM:transitivity}.  The final item in the Lemma, assertion \ref{L:V^BisBVM:membership}, is a straightforward consequence of (\ref{Eq:Bmembership2}) and (\ref{Eq:Bmembership3}), keeping in mind that these have been proven for all members of $V^\B$.
\end{proof}

\section{The forcing relation $\forces$}\label{SS:ForcingRelation}

The efficacy of forcing is due in large part to the fact that $V$ has access to information about truth in a generic extension.  This access is in the form of the \emph{forcing relation} $\forces$.  The forcing relation is defined in terms of the forcing language, which consists of the language of set theory together with constants for every \emph{$\P$-name}.  Throughout the remaining development, I will simply take $\P$-names to be the same as $\B$-names, and will use $\V^\P$ and $\V^\B$ interchangeably to refer to the same object (similarly for $V^\P$ and $V^\B$).\footnote{Strictly speaking, $V^\P$ is defined recursively by $\dot{a}\in V^\P \iff \dot{a}$ is a set of pairs $\<\tau,p>$ with $\tau \in V^\P$ and $p\in\P$. Similarly, $\V^\P$ is the collection of all classes $\class{A} \in \V$ such that $\class{A} \subset V^\P \times \P$.  The collection $\<V^\P,\V^\P>$ of \emph{$\P$-names} matches that of $\B$-names given in Definitions \ref{D:V^Bsets} and \ref{D:V^Bclasses}, with the modification that functions are allowed to be set multi-valued.  There is a natural correspondence between $\P$-names and $\B$-names, and in terms of forcing the two objects yield the same forcing extensions.  For the sake of simplicity I will work only with the objects $\V^\B$ and $V^\B$, sometimes denoting them $\V^\P$ and $V^\P$ as convenient.}  
$\P$-names are sometimes denoted by a dotted letter, as in $\dot{a}$ or $\classdot{A}$, and sometimes by $\tau, \sigma, x, y, $ etc. as in the previous section, depending on context and preference.

Throughout this section, I will take $\P=\bigcup P_\alpha$ to be a chain of complete subposets with $\B = \bigcup B_\alpha$ the associated chain of complete subalgebras, and $e:\P\rightarrow \B$ the canonical embedding.  

\begin{Def}\rm\label{D:ForcingRelation}  The \emph{forcing language} is simply the language of Bernays-G\"odel set theory together with constants added for every $\P$-name or member of $\V^\P$.  The \emph{forcing relation} $\forces$ is a subclass of $\P \times \{$sentences in the forcing language, e.g. $\phi(\dot{a}_1,...,\dot{a}_n, \classdot{A}_1,...,\classdot{A}_n)\}$, and is defined recursively on the complexity of $\phi$.   The constants $\dot{a}_1,...,\dot{a}_n, \classdot{A}_1,...,\classdot{A}_n$ are suppressed for readability.
\begin{enumerate}
\parskip=0pt

\item \label{D:ForcingRelation:SetAtomic}If $\phi$ is an atomic formula involving only set constants, or members of $V^\B$, then $p\forces\phi \iff e(p) \leq \bval{\phi}$.

\item \label{D:ForcingRelation:ClassAtomic}If $\phi$ is an atomic formula involving at least one proper class constant, or member of $\V^\B \setminus V^\B$, then $p\forces \phi$ is defined:
	\begin{enumerate}
	\parskip=0pt
	
	\item \label{D:ForcingRelation:ClassAtomic:subset}$p \forces \classdot{A} \subset \classdot{C} \iff \forall \<x,u> \in \classdot{A} \, \{q \leq p \mid \exists \<y,v> \in \classdot{C}$ such that $({e(q) \leq u \rightarrow \bval{x=y}\btimes v})\}$ is dense below $p$ (recall $u\rightarrow v$ is defined $-u+v$).

	\item \label{D:ForcingRelation:ClassAtomic:equality}$p \forces \classdot{A} = \classdot{C} \iff p \forces \classdot{A} \subset \classdot{C}$ and $p \forces \classdot{C} \subset \classdot{A}$.

	\item \label{D:ForcingRelation:ClassAtomic:membership}$p \forces \classdot{A} \in \classdot{C} \iff \{q \mid \exists \<y,v> \in \classdot{C} ({e(q)\leq v \mbox{ and } q\forces \classdot{A}=y  })\}$ is dense below $p$. 

	\end{enumerate}
\item $p\forces \neg\phi \iff \neg \exists q \leq p \(q\forces\phi)$.

\item $p\forces \phi\wedge\psi \iff p\forces \phi$ and $p\forces \psi$.

\item $p\forces \phi\vee\psi \iff p\forces \neg\(\neg\phi \wedge \neg \psi)$.

\item $p\forces \exists x \, \phi(x) \iff \forall q \leq p \,\exists r\leq q \, \exists \classdot{A}\in \V^\P \({r\forces \phi(\classdot{A})})$.

\item $p\forces \forall x \, \phi(x) \iff p\forces \neg \exists x \,\neg \phi(x)$.
\end{enumerate}
\end{Def}

Note that the introduction of the additional relation symbol $\subset$ in case \ref{D:ForcingRelation:ClassAtomic} is for clarity only, as the definition can be modified to avoid it.   It is worth observing that in the case of a true Boolean-valued model, where $\bval{\phi}$ exists for every $\phi$, the  definition used in case \ref{D:ForcingRelation:SetAtomic} above can be extended to all formulas.  In such a case, the remaining parts of the definition are derivable from $e(p)\leq\bval{\phi}$.
Several basic facts about the forcing relation follow.  Note that for proofs by induction on the complexity of formulas it suffices to consider atomic formulas, negations, conjunctions, and existential quantifiers, as the other operators are defined in terms of these.

\begin{Lemma}\label{L:ForcingRelationFacts} Suppose $\forces$ is defined as above.  Then every sentence $\phi$ in the forcing language satisfies the following.
\begin{enumerate}
\parskip=0pt
\item \label{L:ForcingRelationFacts:strongerconditions} If $p\forces \phi$ and $q\leq p$, then $q\forces \phi$.
\item \label{L:ForcingRelationFacts:densetoforce} If $\forall q \leq p \, \exists r \leq q \, \(r\forces \phi)$, then $p\forces \phi$.
\item \label{L:ForcingRelationFacts:doublenegation} $p\forces \phi \iff p\forces \neg\(\neg\phi)$.
\item \label{L:ForcingRelationFacts:disjunction} $p\forces \phi\vee\psi \iff \forall q\leq p \, \exists r\leq q \(r\forces \phi \mbox{ or } r\forces \psi)$.
\item \label{L:ForcingRelationFacts:universal} $p\forces \forall x \, \phi(x) \iff \forall \classdot{A}\in V^\P \({p\forces \phi(\classdot{A})})$.
\item \label{L:ForcingRelationFacts:nocontradiction} If $p\forces \phi$, then $p \nforces \neg \phi$.
\item \label{L:ForcingRelationFacts:everythingisforced} $\forall p\in \P \, \exists q \leq p \(q\forces \phi \vee q\forces \neg\phi)$.
\end{enumerate}
\end{Lemma}

\begin{proof}   Assertion \ref{L:ForcingRelationFacts:strongerconditions} is proved by induction  formulas.  If $\phi$ is an atomic formula involving only sets, then $p\forces \phi \iff e(p)\leq \bval{\phi}$.  If $q\leq p$, then $e(q)\leq e(p) \leq \bval{\phi}$, and so $q\forces \phi$.  For atomic formulas involving proper classes, note that if a class of conditions is dense below $p$ and $q\leq p$, then the class is dense below $q$.  Thus if $p \forces \classdot{A}\subset \classdot{C}$ or $p\forces \classdot{A}\in \classdot{C}$, then $q\leq p$ will also force that statement.  As $p\forces \classdot{A}=\classdot{C}$ is defined in terms of $\subset$, this argument also suffices to prove assertion \ref{L:ForcingRelationFacts:strongerconditions} for equality of classes.  This concludes the atomic case.  Now suppose $p\forces \neg \phi$.  By the definition, this holds only if no $p^\prime\leq p$ forces $\phi$, and clearly this will also be true of any $q\leq p$, so such a $q$ will also force $\neg \phi$.  Next consider the case $p\forces \phi \wedge \psi$.
\begin{eqnarray*}
	p\forces \phi \wedge \psi & \iff & p\forces \phi \mbox{ and } p\forces\psi \mbox{ (definition)}\\
	 & \rightarrow & \forall q \leq p , \(q\forces \phi) \mbox{ and } \forall q \leq p , \(q\forces\psi) \mbox{ (inductive hypothesis)}\\
	 & \iff & \forall q \leq p, \(q\forces\phi\wedge\psi)
\end{eqnarray*}
Finally, suppose $p\forces \exists x \phi(x)$.  By the definition I have $\forall q \leq p \,\exists r\leq q \, \exists \dot{a}\in V^\P \({r\forces \phi(\dot{a})})$.  If I fix $q\leq p$ and consider any $q^\prime \leq q$, then $q^\prime$ will also lie below $p$ and so $\forall q^\prime \leq q \,\exists r\leq q^\prime \, \exists \dot{a}\in V^\P \({r\forces \phi(\dot{a})})$.  Thus $q \forces \exists x \phi(x)$.   Therefore assertion \ref{L:ForcingRelationFacts:strongerconditions} holds for all $\phi$.

Assertion \ref{L:ForcingRelationFacts:densetoforce} states that if it is dense below $p$ that $\phi$ is forced, then $p$ itself forces $\phi$.  Once again the proof is by induction on formulas, but in this case set atomic formulas present the biggest challenge.  Suppose $\phi$ is an atomic formula involving only sets, and $p$ is such that $\forall q \leq p \exists r \leq q \, \(r\forces \phi)$.  Suppose further that $p \nforces \phi$.  Then $e(p) \nleq \bval{\phi}$, and by Definition \ref{D:BAordering} this gives $e(p)\btimes\(-\bval{\phi}) \ne 0$.  Let $e(p)\btimes\(-\bval{\phi}) = u$.  $\P$ embeds densely into $\B \setminus \{0\}$ via the embedding $e$, and so there must be $q\in \P$ such that $e(q) \leq u$.  Clearly $e(q) \leq -\bval{\phi}$.  However, $e(q) \leq e(p)$ implies $q \leq p$, and so by hypothesis $\exists r \leq q \, \(r \forces \phi)$.  Thus $e(r) \leq \bval{\phi}$, but also $e(r)\leq e(q) \leq -\bval{\phi}$.  This is a contradiction, as the only member of $\B$ lying below both $\bval{\phi}$ and $-\bval{\phi}$ is $0$, but $e$ maps into $\B \setminus \{0\}$.  For atomic formulas involving proper classes, assertion \ref{L:ForcingRelationFacts:densetoforce} follows directly from the definition for $p \forces \classdot{A}\subset \classdot{C}$ and $p\forces \classdot{A}\in \classdot{C}$, and for $p \forces \classdot{A}=\classdot{C}$ it follows inductively from $p \forces \classdot{A}\subset \classdot{C}$.

Now suppose $\forall q \leq p \, \exists r \leq q \, \(r\forces \phi\wedge\psi)$.  Then $\forall q \leq p \, \exists r \leq q \, (r\forces \phi$ and $r\forces\psi)$.  The inductive hypothesis gives $p\forces \phi$ and $p\forces \psi$, and so $p\forces \phi \wedge \psi$.

For the existential quantifier case, suppose $\forall q \leq p \, \exists r \leq q \, \({r\forces \exists x \phi(x)})$.  Then $\forall q \leq p \,\exists r\leq q \, \exists \classdot{A}\in \V^\P \({r\forces \phi(\classdot{A})})$, and so $p\forces \exists x \phi(x)$ by the definition.  Thus assertion \ref{L:ForcingRelationFacts:densetoforce} holds.

For assertion \ref{L:ForcingRelationFacts:doublenegation},
\begin{eqnarray*}
	p\forces \neg\(\neg\phi) & \iff & \neg \exists q \leq p \, \(q \forces \neg \phi) \mbox{ (definition)}\\
	& \iff & \neg \exists q \leq p \,  \neg \exists r \leq q \, \( r\forces \phi) \mbox{ (definition)}\\
	& \iff & \forall q \leq p \, \exists r \leq q \, \(r\forces \phi)\\
	& \iff & p \forces \phi \mbox { (assertion \ref{L:ForcingRelationFacts:densetoforce})}
\end{eqnarray*}

Assertions \ref{L:ForcingRelationFacts:disjunction} and \ref{L:ForcingRelationFacts:universal} give direct characterizations of forcing for disjunctions and universal quantifiers.

For assertion \ref{L:ForcingRelationFacts:disjunction}, 
\begin{eqnarray*}
	\lefteqn{p\forces \phi \vee \psi} \\
	& \iff & p\forces \neg\(\neg\phi \wedge \neg\psi) \mbox{ (definition)}\\
	& \iff & \neg \exists q \leq p \({q \forces \(\neg\phi \wedge \neg\psi)}) \mbox{ (definition)}\\
	& \iff & \forall q \leq p \({ q \nforces \(\neg\phi \wedge \neg\psi)}) \mbox{ (logical equivalence)}\\
	& \iff & \forall q \leq p \({ q \nforces \neg\phi \mbox { or } q \nforces \neg\psi}) \mbox{ (definition)}\\
	& \iff & \forall q \leq p \({ \({\exists r \leq q \(r\forces\phi)})     \mbox { or } \({\exists r \leq q \(r\forces\psi)}) }) \mbox{ (definition)} \\
	& \iff & \forall q \leq p \, \exists r \leq q \({r\forces \phi \mbox { or } r \forces \psi })
\end{eqnarray*}
and for assertion \ref{L:ForcingRelationFacts:universal}, 
\begin{eqnarray*}
	\lefteqn{p\forces \forall x \, \phi(x) }\\
	& \iff & p \forces \neg \exists x \, \neg \phi(x) \mbox{ (definition)}\\
	& \iff & \neg \exists q\leq p \, \({q \forces \exists x \, \neg \phi(x) }) \mbox{ (definition)}\\
	& \iff & \forall q\leq p \, \({q \nforces \exists x \, \neg \phi(x) }) \mbox{ (logical equivalence)}\\
	& \iff & \forall q\leq p \, \neg\({\forall r\leq q \, \exists s \leq r \, \exists \classdot{A} \in \V^\P \({s \forces \neg \phi(\classdot{A})}) }) \mbox{ (definition)}\\
	& \iff & \forall q\leq p \, \exists r\leq q \, \forall s \leq r \, \forall \classdot{A} \in \V^\P \, \({s \nforces \neg \phi(\classdot{A})})  \mbox{ (logical equivalence)}\\
	& \iff & \forall q\leq p \, \exists r\leq q \, \forall s \leq r \, \forall \classdot{A} \in \V^\P \, \exists t \leq s \({t \forces \phi(\classdot{A}) }) \mbox{ (definition)}\\
	& \iff & \forall q\leq p \, \exists r\leq q \, \forall \dot{a} \in \V^\P \, \({ r \forces \phi(\dot{a}) }) \mbox{ (assertion \ref{L:ForcingRelationFacts:densetoforce})}\\
	& \iff &  p \forces \phi(\classdot{A}) \mbox{ (assertion \ref{L:ForcingRelationFacts:densetoforce})}\\
\end{eqnarray*}

Assertion \ref{L:ForcingRelationFacts:nocontradiction} can be argued directly for any $\phi$.  Suppose $p\forces \phi$.  Then there is $q\leq p$ such that $q\forces \phi$, since $p$ is such a $q$.  Thus $p\nforces \neg\phi$, by the definition of $p\forces \neg\phi$.

Assertion \ref{L:ForcingRelationFacts:everythingisforced} also follows directly from the definition.  Fix $p$ and $\phi$.  If there is $q\leq p$ such that $q\forces \phi$, then I am done.  If not, then by the definition $p\forces \neg\phi$, and so I am done.  This concludes the proof of Lemma \ref{L:ForcingRelationFacts}. \end{proof}

\section{Existence of generic filters}\label{SS:ConsistencyOfGenericFilters}

In forcing over models of $\BGC$ one encounters the metamathematical issue of the consistency of the existence of $\V$-generic filters.  There are several common approaches to this question in the context of forcing with partially ordered sets over $\ZFC$ models, and in general these approaches are readily adaptable to $\BGC$ models and forcing with proper classes.  For the sake of concreteness, I will take the approach commonly used in introductions to forcing, that of \emph{countable transitive models}.  I will first argue that $\V$-generic filters exist for such models, and I will then show how this argument can be adapted to prove consistency results.  Suppose that $\<V,\V,\in>$ is a countable transitive model of $\BGC$, and consider a  partial order $\P\in\V$.  As there are only countably many dense subclasses of $\P$ in $\V$, I can enumerate them $\<\class{D}_n \mid n\in \omega>$.  I construct a descending sequence $\<p_n \mid n \in  \omega>$ in $\P$ with $p_n \in \class{D}_n$ for each $n$.  Taking $\class{G}$ to the be the upwards closure of $\{p_n \mid n\in \omega\}$ in $\P$, it is easily established that $\class{G}$ is a filter and, as it meets every dense subclass of $\P$ in $\V$, it is $\V$-generic.  This allows me to consider the forcing extension $\<V[\class{G}], \V[\class{G}],\in>$, and forcing arguments can be used to establish that the extension satisfies certain desired properties, e.g. $\BGC$ together with some assertion $\phi$ forced by $\P$.  

This argument presupposes the existence of a countable transitive model of $\BGC$.  As the existence of such a model cannot be proven from $\BGC$ itself, the approach must be modified to achieve the consistency results desired in this paper.  However, a `finite' version of the argument also applies.  If $\Gamma$ is a finite collection of sentences forced by $\P$, then I require only a finite portion of $\BGC$ in the ground model for the forcing argument to be carried out.  That is, for any finite collection $\Gamma$ of sentences forced by $\P$, there is a finite fragment $\Gamma^*$ of $\BGC$ such that any countable transitive model of $\Gamma^*$ has a forcing extension satisfying $\Gamma$.  Now suppose I wish to establish a result of the form ``$\Con(\BGC) \rightarrow \Con(\BGC+\phi)$'' for some forceable assertion $\phi$.  If $\neg\Con(\BGC+\phi)$ holds, then there is a finite fragment $\Gamma$ of $\BGC+\phi$ which is inconsistent.  Taking $\Gamma^*\subset \BGC$ as above, the Reflection Theorem can be used to construct a countable transitive model of $\Gamma^*$.  Forcing is then used to construct a model of $\Gamma$, yielding a contradiction.  Therefore $\Con(\BGC) \rightarrow \Con(\BGC+\phi)$ holds.

The argument above provides the foundation for the applications of forcing found in Chapters \ref{C:TheGroundAxiom} and \ref{C:GAandVnotHOD} as well as for the development of forcing in this appendix.  However, in practice this approach is seldom explicit, and forcing arguments are generally carried out by assuming the existence of a $\V$-generic filter.  While not strictly correct, this is a harmless simplification in light of the discussion above.

\section{The model $\V[\class{G}]$}\label{SS:V[G]}

I now define the generic extension $\<V[\class{G}],\V[\class{G}]>$ where $\class{G}$ is a $\V$-generic filter on a partial order $\P$.  

\begin{Def}\rm\label{D:valuationfunction}Suppose $\P$ is a partial order, $\B$ the associated Boolean algebra and $e:\P\rightarrow\B$ the canonical embedding.  If $\class{G}\subset\P$ is a $\V$-generic filter on $\P$, then the \emph{valuation function} $i_\class{G}$ is defined recursively on $\P$-names by 
$$i_\class{G}(\classdot{A}) = \left\{i_\class{G}(\dot{a}) \mid \exists u \in \B \, \exists p \in \class{G} \, \({ \<\dot{a},u> \in \classdot{A} \mbox{ and } e(p) \leq u }) \right\}.$$
The \emph{generic extension} $\<V[\class{G}],\V[\class{G}],\in>$ consists of $\V[\class{G}]$, the collection of all $i_\class{G}(\classdot{A})$ for $\classdot{A} \in \V^\P$, and $V[\class{G}]$, the collection of all $i_\class{G}(\dot{a})$ for $\dot{a} \in V^\P$, together with the $\in$ relation.  If $\class{A} \in \V[\class{G}]$ and $\classdot{X} \in \V^\P$ such that $i_\class{G}(\classdot{X})=\class{A}$, then $\classdot{X}$ is said to be a \emph{$\P$-name for $\class{A}$}.
\end{Def}

Note that $V^\P \subset \V^\P$ implies $V[\class{G}] \subset \V[\class{G}]$.  The `function' $i_\class{G}$ is not a function in the sense of a class of ordered pairs, as its domain includes many classes.  However, in the same sense that a $\ZFC$ model can have `class functions' defined by first-order formulas, so does $i_\class{G}$ provide a uniform definition associating each $\classdot{A}\in\V^\P$ with its valuation $i_\class{G}(\classdot{A})$.    Definition \ref{D:valuationfunction} assumes the existence of an `associated Boolean algebra $\B$.'  If $\P$ is a set, then $\B = ro(\P$), and  if $\P$ is a chain of complete subposets, then $\B$ is the associated chain of complete subalgebras.  These are the only cases that I will be concerned with here.  However, it is worth observing that both $\V^\P$ and the valuation function can be defined entirely in terms of $\P$ without reference to Boolean algebras.

Although Theorem \ref{T:ForcingTheorem}, the Forcing Theorem, will be required to prove that the extension $\<V[\class{G}],\V[\class{G}]>$ satisfies any significant portion of $\BGC$, a few basic facts follow directly from the definition.  One concerns the definition of the sets $V[\class{G}]$ of $\V[\class{G}]$, and the other two concern transitivity of $\V[\class{G}]$.

\begin{Lemma}\label{T:V[G]BasicFacts}Suppose $\P$ and $\class{G}$ satisfy the hypotheses of Definition \ref{D:valuationfunction}.
\begin{enumerate}
\parskip=0pt
\item \label{T:V[G]BasicFacts:definitionofsets} $V[\class{G}]$ is definable in $\V[\class{G}]$ as $\{\class{X} \mid \exists \class{Y} \, (\class{X} \in \class{Y})\}$.
\item \label{T:V[G]BasicFacts:V[G]transitive} $\V[\class{G}]$ is transitive.  
\item \label{T:V[G]BasicFacts:V[G]transitive2} If $\class{X} \subset V^\P$ and for every $\dot{a}\in \classdot{X}$ the domain $\dom(\dot{a})$ is a subset of $\class{X}$, then $i_\class{G}^{\prime\prime} \class{X}$ is transitive.
\end{enumerate}
\end{Lemma}

\begin{proof}For assertion \ref{T:V[G]BasicFacts:definitionofsets}, it is certainly true that for every name $\dot{a} \in V[\class{G}]$ there will be $\class{Y} \in \V[\class{G}]$ such that $i_\class{G}(\dot{a}) \in \class{Y}$.  Simply consider $\class{Y} = i_\class{G}(\dot{c})$ where $\dot{c}=\{\<\dot{a},1>\}$.  Conversely, suppose there is $\class{Y} \in \V[\class{G}]$ such that $\class{X} \in \class{Y}$.  By definition there is a name $\classdot{Y} \in \V^\P$ such that $i_\class{G}(\classdot{Y}) = \class{Y}$, and there must be $\dot{a}\in \dom(\classdot{Y})$ with $i_\class{G}(\dot{a}) = \class{X}$.  As $\dot{a}\in V^\P$, I have $\class{X} \in i_\class{G}^{\prime\prime} V^\P$.  For transitivity, fix $\class{A}\in \V[\class{G}]$ and $a\in\class{A}$.  There is a name $\classdot{A} \in \V^\P$ such that $i_\class{G}(\classdot{A}) = \class{A}$, and a name $\dot{a}\in \dom(\classdot{Y})$ with $i_\class{G}(\dot{a}) = a$.  Since $\dot{a} \in V^\P$, I have $i_\class{G}(\dot{a}) \in V[\class{G}]$, and so $\V[\class{G}]$ is transitive.  The proof of assertion \ref{T:V[G]BasicFacts:V[G]transitive2} is similar.
\end{proof}

Several $\P$-names deserve special mention.  The first is the collection of \emph{check names} associated with members of $\V$.  These are the canonical names for each class $\class{A}\in \V$.  The check name for $\class{A}$ is denoted $\check{\class{A}}$.

\begin{Def}\rm\label{D:checknames}For $\P$ a partial order with associated Boolean algebra $\B$, define for each $\class{A} \in \V$ the \emph{check name} $\check{\class{A}}$ recursively by
\begin{enumerate}
\parskip=0pt
\item $\check{\emptyset} = \emptyset$, and
\item for $\class{A}\in \V$ nonempty, $\check{\class{A}} = \{ \<\check{b}, 1>  \mid  b \in \class{A}\}$, where $1$ is the top element of $\B$.
\end{enumerate}
\end{Def}

Note that if $a$ is a set, then $\check{a}$ is a also a set.  Check names always have the same valuation, regardless of the generic $\class{G}$.

\begin{Lemma}\label{L:checknames}If $\P$ is a partial order with associated Boolean algebra $\B$ and $\class{G}$  is $\V$-generic for $\P$, then for every $\check{\class{A}}$ the valuation $i_\class{G}(\check{\class{A}})=\class{A}$.
\end{Lemma}

\begin{proof}   Note that for any nonempty filter $\class{G}$, every $p\in \class{G}$ satisfies $e(p) \leq 1$.  Thus $i_\class{G}(\check{\class{A}}) = \{i_\class{G}(\check{b})  \mid  b\in \class{A}\}$.  Proceeding by $\in$-induction, this gives $i_\class{G}(\check{\class{A}})  = \{ b  \mid  b\in \class{A}\} = \class{A}$. \end{proof}

Another important name is the canonical name for the generic class $\class{G}$, denoted $\classdot{G}$.  

\begin{Def}\rm\label{D:Gdot}If $\P$ is a partial order, then 
$\classdot{G} = \{\<\check{p},e(p)> \mid p\in\P\}$.\\
If $\P=\bigcup P_\alpha$ is a chain of complete subposets, then 
$\dot{G}_\alpha = \{\<\check{p},e(p)> \mid p\in\P_\alpha\}$ for each $\alpha$.
\end{Def}

The main fact about $\classdot{G}$ is that, when valuated by $i_\class{H}$ for any $\V$-generic $\class{H}$, it yields $\class{H}$ itself.

\begin{Lemma}\label{L:Gdot}If $\P$ is a partial order and $\class{H}$ is $\V$-generic for $\P$, then $i_\class{H}(\classdot{G}) = \class{H}$.
\end{Lemma}

\begin{proof}  Note that this lemma takes care of both cases of Definition \ref{D:Gdot}, as each $G_\alpha$ is $\V$-generic for the set $P_\alpha$.  Applying the definition of $i_\class{H}$ reveals $i_\class{H}(\classdot{G}) = \left\{i_\class{H}(\check{p})  \mid  \exists q \in \class{H} \, \({ e(q) \leq \classdot{G}(\check{p})=e(p)  }) \right\} = \left\{p  \mid  p\in \class{H}\right\}$.  For the last equality, note that $\({ \exists q \in \class{H} \,  e(q)\leq e(p) }) \iff p \in \class{H}$, as $\class{H}$ is closed upwards. \end{proof}

Many forcing arguments involve functions, and it will be useful to have a shorthand way of referring to names for ordered pairs.  This is given by the function $\op:V^\P \times V^\P \rightarrow V^\P$, which inputs two $\P$-names $\dot{a}$ and $\dot{b}$ and outputs a name for the ordered pair $\<a,b>$, where $a=i_\class{G}(\dot{a})$ and $b=i_\class{G}(\dot{b})$.  Recall that by definition an ordered pair is $\<a,b> = \{a, \{a,b\}\}$.

\begin{Def}\rm\label{D:op} For $\dot{a}$ and $\dot{b} \in V^\P$, define 
$$\op(\dot{a},\dot{b}) = \{ \<\dot{a}, 1>,  \< \{ \<\dot{a}, 1>, \<\dot{b}, 1> \}, 1> \},$$
where $1$ is the top element of the Boolean algebra associated with $\P$.
\end{Def}

Note that for any $\class{G}$ which is $\V$-generic for $\P$, the name $op(\dot{a},\dot{b})$ valuates to $i_\class{G}(op(\dot{a},\dot{b})) = \< i_\class{G}(\dot{a}), i_\class{G}(\dot{b})>$.
A final type of important name is a \emph{nice name for a subset of $\dot{a}$}.

\begin{Def}\rm\label{D:nicenames}If $\P$ is partial order and $\dot{a}\in V^\P$, then $\dot{c} \in V^\P$ is a \emph{nice name for a subset of $\dot{a}$} if and only if $\dot{c} \from \dom(\dot{A}) \rightarrow \B$.
\end{Def}

The point of nice names is that they give representatives of every subset of $\dot{a}$ in a generic extension by $\P$.

\begin{Lemma}\label{L:nicenames}Suppose $\P$ is a partially ordered set or a chain of complete subposets and $\dot{a} \in V^\P$.  If $\class{G}$ is $\V$-generic for $\P$ and $z \in V[\class{G}]$ is a subset of $i_\class{G}(\dot{a})$, then there is $\dot{c}$ a nice name for a subset of $\dot{a}$ such that $i_\class{G}(\dot{c})=z$.  Furthermore, such a $\dot{c}$ exists with $\bval{\dot{c}\subset \dot{a}} = 1$.
\end{Lemma}

\begin{proof}  Fix any name $\dot{z}$ for $z$.  Define $\dot{c} = \{\<\sigma,u_\sigma> \mid \sigma \in \dom(\dot{a}) \}$, where $u_\sigma = \sum_{\tau \in \dom(\dot{z})} \dot{z}(\tau) \btimes \bval{\tau = \sigma} \btimes \bval{\sigma\in\dot{a}}$.  Clearly $\dot{c}$ is a nice name for a subset of $\dot{a}$.  It is left to the reader to verify that $i_\class{G}(\dot{c})=z$ and $\bval{\dot{c} \subset \dot{a}} = 1$. \end{proof}

The recursive nature of the valuation function yields a nice relationship between the rank of $i_\class{G}(\tau)$ and the $\B$-rank $\rho(\tau)$.

\begin{Lemma}\label{L:V[G]rank}For $\P$ and $\class{G}$ as above, if $\tau \in V^\B$, then $\rank(i_\class{G}(\tau)) \leq \rho(\tau)$.
\end{Lemma}

\begin{proof}   By induction on $\rho(\tau)$.  Suppose the Lemma holds for every $\sigma$ with $\rho(\sigma) < \rho(\tau)$.  Then in particular it holds for every $\sigma$ in $\dom(\tau)$.  Clearly $\rank(i_\class{G}(\tau)) \leq \sup\{ \rank(i_\class{G}(\sigma))+1 \mid \sigma \in \dom(\tau)\}$.  Since $\rho(\sigma) < \rho(\sigma)$ for every $\sigma \in \dom(\tau)$, I have $\rank(i_\class{G}(\tau)) \leq \rho(\tau)$. \end{proof}

Finally, the valuation function is well-behaved with respect to complete subposets.

\begin{Lemma}\label{L:CompleteSubposetsandValuation}Suppose $\P \subset_c \Q$ are partial orders with associated Boolean algebras $\B \subset_c \D$, and $\class{H}\subset\Q$ is $\V$-generic.  If $\class{G}=\class{H}\cap\P$, then  $i_\class{G}(\tau) = i_\class{H}(\tau)$ for every $\tau\in V^\P$.
\end{Lemma}

\begin{proof}   By induction on $\rho(\tau)$.  Suppose the Lemma holds for every $\sigma \in \dom(\tau)$.  By definition $i_\class{G}(\tau)$ consists of $i_\class{G}(\sigma)$ for all $\sigma \in \dom(\tau)$ satisfying $\exists p\in \class{G} \, ( e(p)\leq \tau(\sigma) )$.  As $i_\class{H}(\tau)$ is defined similarly, and  $i_\class{H}(\sigma)=i_\class{G}(\sigma)$ for all $\sigma\in \dom(\tau)$ by inductive hypothesis, it suffices to show that for every $u=\tau(\sigma)$ in $\B$, $$\exists p \in \class{G} \, e(p) \leq u \iff \exists q \in \class{H} \, e(q) \leq u.$$ The forward implication is immediate, as $\class{H}\subset \class{G}$.  For the reverse implication, suppose $u=\tau(\sigma)$ and suppose $q\in \class{H}$ is such that $e(q)\leq u$. Consider the following set $A=\{p\in \P \mid e(p)\leq u \mbox{ or } e(p) \perp u\}$.  It is a routine matter to verify that $A$ is dense in $\P$, and so $\class{G}\cap A$ is nonempty.  Suppose $g\in \class{G}\cap A$.  I claim $e(g) \leq u$.  Otherwise $e(g)\perp u$, but $e(q)\leq u$ gives $e(g) \perp e(q)$ in $\D$.  As $e$ preserves compatibility, I conclude $g\perp q$ in $\Q$.  However, $q\in \class{H}$ by hypothesis and $g\in \class{G}  \subset \class{H}$, yielding a contradiction as $\class{H}$ is a filter and so all members are compatible.  \end{proof}

Lemma \ref{L:CompleteSubposetsandValuation} gives a nice characterization of $V[\class{G}]$ in the case where $\P=\bigcup P_\alpha$ is a chain of complete subposets.  If $\class{G}$ is $\V$-generic for $\P$, then $G_\alpha = \class{G} \cap P_\alpha$ is $\V$-generic for $P_\alpha$ for every $\alpha$.  Since every $\P$-name is actually a $P_\alpha$ name for some $\alpha$, Lemma \ref{L:CompleteSubposetsandValuation} shows that $V[\class{G}] = \bigcup V[G_\alpha]$.

\section{The Forcing Theorem}\label{SS:ForcingTheorem}

The Forcing Theorem states that satisfaction in the model $\V[\class{G}]$ is exactly characterized by the forcing relation restricted to $\class{G}$.

\clearpage
\begin{Thm}\label{T:ForcingTheorem}\emph{Forcing Theorem}.  If $\P$ is a poset or a chain of complete subposets and $\class{G}$ is $\V$-generic for $\P$, then for any formula $\phi(x_1,...,x_n, \class{X}_1,...,\class{X}_n)$, for any parameters $a_1,...,a_n \in V[\class{G}]$ and $\class{A}_1,...,\class{A}_n\in \V[\class{G}]$ with respective names $\dot{a}_1,...,\dot{a}_n \in V^\P$ and $\classdot{A}_1,...,\classdot{A}_n\in \V^\P$,  
$$\V[\class{G}] \vDash \phi(a_1,...,a_n, \class{A}_1,...,\class{A}_n) \iff \exists p\in \class{G} \, \({ p \forces \phi(\dot{a}_1,...,\dot{a}_n, \classdot{A}_1,...,\classdot{A}_n) })$$
\end{Thm}

\begin{proof}   The theorem is proved by induction on the complexity of $\phi$.  The atomic case will rely on Definition \ref{D:V^Batomic}, and much of the argument will take place in the Boolean algebra $\B$ associated with $\P$.  With that in mind, I begin by observing that if I take $\overline{\class{G}}$ to be the upwards closure of $e^{\prime\prime}\class{G}$ in $\B$, then by Lemmas \ref{L:ChainOfCompleteSubalgebras} and \ref{L:GenericUltrafilter=GenericFilter}, $\overline{\class{G}}$ is a $\V$-generic ultrafilter on $\B$.  Furthermore, for all $u\in \B$ it is clear that $u\in \overline{\class{G}} \iff \exists p\in \class{G} \, \({ e(p) \leq u })$. 

For atomic formulas $\phi$ involving only sets, $p\forces \phi \iff e(p) \leq \bval{\phi}$.  By the comment above, $\exists p \in \class{G} \, \({ e(p) \leq \bval{\phi}}) \iff \bval{\phi} \in \overline{\class{G}}$. Thus the atomic case of Theorem \ref{T:ForcingTheorem} can be reduced to following statements, which will be proved simultaneously by induction on $\<\rho(\tau),\rho(\sigma)>$ for $\tau,\sigma \in V[\class{G}]$.

\setcounter{equation}{0}
\begin{eqnarray}
	\label{Eq:EvaluationMembership1} i_\class{G}(\tau) \in i_\class{G}(\sigma) &\iff& \bval{\tau\in\sigma} \in \overline{\class{G}}\\
	\label{Eq:EvaluationMembership2} i_\class{G}(\sigma) \in i_\class{G}(\tau)&\iff& \bval{\sigma\in\tau} \in \overline{\class{G}}\\
	\label{Eq:EvaluationSubset1} i_\class{G}(\tau) \subset i_\class{G}(\sigma) &\iff& \bval{\tau\subset\sigma} \in \overline{\class{G}}\\
	\label{Eq:EvaluationSubset2} i_\class{G}(\sigma) \subset i_\class{G}(\tau) &\iff& \bval{\sigma\subset\tau} \in \overline{\class{G}}\\
	\label{Eq:EvaluationEquality} i_\class{G}(\tau) = i_\class{G}(\sigma) &\iff& \bval{\tau=\sigma} \in \overline{\class{G}}
\end{eqnarray}

For statement (\ref{Eq:EvaluationMembership1}), 
\begin{eqnarray*}
	\lefteqn{\bval{\tau\in\sigma} \in \overline{\class{G}}}\\
	& \iff & \({\sum_{x\in \dom(\sigma)} \bval{\tau=x}\btimes\sigma(x)}) \in \overline{\class{G}}\\
	& \iff & -\({\sum_{x\in \dom(\sigma)} \bval{\tau=x}\btimes\sigma(x)}) \notin \overline{\class{G}} \mbox{ ($\overline{\class{G}}$ is ultrafilter)}\\
	& \iff & \({\prod_{x\in \dom(\sigma)} -\({\bval{\tau=x}\btimes\sigma(x)}) }) \notin \overline{\class{G}}\\
	& \iff & \exists x\in \dom(\sigma) \, -\bval{\tau=x}\btimes\sigma(x) \notin \overline{\class{G}} \mbox{ ($\overline{\class{G}}$ is generic)}\\
	& \iff & \exists x\in \dom(\sigma) \, \bval{\tau=x}\btimes\sigma(x) \in \overline{\class{G}} \mbox{ (ultrafilter)}\\
	& \iff & \exists x\in \dom(\sigma) \, i_\class{G}(\tau)=i_\class{G}(x) \wedge \sigma(x)\in \overline{\class{G}} \mbox{ (hypothesis (\ref{Eq:EvaluationEquality}))}\\
	& \iff & i_\class{G}(\tau)\in i_\class{G}(\sigma) \mbox{ (Definition \ref{D:valuationfunction})}
\end{eqnarray*}

Statement (\ref{Eq:EvaluationMembership2}) follows from exactly the same proof with the roles of $\tau$ and $\sigma$ reversed.  Note that the same inductive hypothesis is used in each case.

For statement (\ref{Eq:EvaluationSubset1}), 
\begin{eqnarray*}
	\lefteqn{\bval{\tau\subset\sigma} \in \overline{\class{G}} }\\
	& \iff &  \({ \prod_{x\in \dom(\tau)} -\tau(x)+\bval{x\in\sigma} })  \in \overline{\class{G}}\\
	& \iff & \forall x \in \dom(\tau) \, \({ -\tau(x) + \bval{x\in\sigma} }) \in \overline{\class{G}} \mbox{ (upwards closure)}\\
	& \iff &  \forall x \in \dom(\tau) \, \tau(x) \in \overline{\class{G}} \rightarrow  \bval{x\in\sigma} \in \overline{\class{G}}\\
	& \iff &  \forall x \in \dom(\tau) \, \tau(x) \in \overline{\class{G}} \rightarrow  i_\class{G}(x)\in i_\class{G}(\sigma) \mbox{ (hypothesis (\ref{Eq:EvaluationMembership1}))}\\
	& \iff &  i_\class{G}(\tau) \subset i_\class{G}(\sigma) \mbox{ (Definition \ref{D:valuationfunction})}
\end{eqnarray*}

Statement (\ref{Eq:EvaluationSubset2}) follows from the same proof, reversing the roles of $\tau$ and $\sigma$ and using hypothesis (\ref{Eq:EvaluationMembership2}) in place of (\ref{Eq:EvaluationMembership1}) in the second to last step.

Finally, 
\begin{eqnarray*}
	\bval{\tau=\sigma} \in \overline{\class{G}} & \iff & \bval{\tau\subset\sigma}\btimes\bval{\sigma\subset\tau}  \in \overline{\class{G}}\\
	& \iff & \bval{\tau=\sigma} \in \overline{\class{G}} \mbox{ and } \bval{\sigma=\tau} \in \overline{\class{G}}\\
	& \iff & i_\class{G}(\tau)\subset i_\class{G}(\sigma) \mbox{ and } i_\class{G}(\sigma)\subset i_\class{G}(\tau)\\
	& \iff & i_\class{G}(\tau) = i_\class{G}(\sigma)
\end{eqnarray*}

This concludes the case where $\phi$ is an atomic formula involving only sets.  Next I will consider formulas involving at least one class variable, following the order described in Definition \ref{D:ForcingRelation}, cases \ref{D:ForcingRelation:ClassAtomic:subset}-- \ref{D:ForcingRelation:ClassAtomic:membership}.  For case \ref{D:ForcingRelation:ClassAtomic:subset}, I must show 
$$\V[\class{G}] \vDash i_\class{G}(\classdot{A}) \subset i_\class{G}(\classdot{C}) \iff \exists p \in \class{G} \, \({ p \forces \classdot{A} \subset \classdot{C} }).$$  
For the forward direction, suppose $\V[\class{G}] \vDash i_\class{G}(\classdot{A})$ but there is no $p\in \class{G}$ forcing this statement.  Since it is dense in $\P$ that either $\phi$ or $\neg\phi$ is forced (Lemma \ref{L:ForcingRelationFacts}, assertion \ref{L:ForcingRelationFacts:everythingisforced}), it follows that there is $p\in\class{G}$ such that $p\forces \neg({\classdot{A} \subset \classdot{C}})$.  For each $q\leq p$ and $\<x,u> \in \classdot{A}$, let $\class{D}_{q,\<x,u>} = \{r\leq q \mid  \exists \<y,v> \in \classdot{C} \, \({e(r) \leq (u \rightarrow \bval{x=y}\btimes v}))\}$.
\begin{eqnarray*}
	p\forces \neg({\classdot{A} \subset \classdot{C}}) & \iff & \neg \exists q \leq p \, ({ q \forces \classdot{A} \subset \classdot{C} })\\
	& \iff & \forall q \leq p \, \exists \<x,u> \in \classdot{A} \, \({\class{D}_{q,\<x,u>} \mbox{ is not dense below } q}) \\
	& \iff & \exists q_0 \leq p \, \exists \<\dot{a},u> \in \classdot{A} \, \forall r \leq q_0 \, \({\class{D}_{r,\<\dot{a},u>}=\emptyset })
\end{eqnarray*}
I will show that $i_\class{G}(\dot{a}) \in i_\class{G}(\classdot{A}) \setminus i_\class{G}(\classdot{C})$, by showing that $\classdot{A}(\dot{a}) \in \overline{\class{G}}$ but $i_\class{G}(\dot{a}) \neq c$ for any $c \in i_\class{G}(\classdot{C})$.  Fix $c\in i_\class{G}(\classdot{C})$ and $\<\dot{c},v> \in \classdot{C}$ such that $v \in \overline{\class{G}}$ and $i_\class{G}(\dot{c})=c$.  Untangling the meaning of $\class{D}_{r,\<\dot{a},u>}=\emptyset$ and applying the definition of the Boolean relations $\rightarrow$ and $\leq$, it follows that for each $r\leq q_0$, either $ 0 < e(r) \btimes u \btimes -\bval{\dot{a}=\dot{c}}$ or $0 < e(r) \btimes u \btimes -v$.  Since $\P$ embeds densely into $\B\setminus \{0\}$, this means that there is $r^\prime \in \class{G}$ such that either $e(r^\prime) < e(r) \btimes u \btimes -\bval{\dot{a}=\dot{c}}$ or $e(r^\prime) < e(r) \btimes u \btimes -v$.  The set of such $r^\prime$ is dense below $q_0$, so there exists one such condition $r_0 \in \class{G}$.  It follows that either $e(r_0) < e(q_0) \btimes u \btimes -\bval{\dot{a}=\dot{c}}$ or $e(r_0) < e(q_0) \btimes u \btimes -v$.  Note that either case implies $u \in \overline{\class{G}}$, and so $i_\class{G}(\dot{A}) \in i_\class{G}(\classdot{A})$.  The second case $e(r_0) < e(q_0) \btimes u \btimes -v$ cannot hold, for then $-v \in \overline{\class{G}}$, contradicting our assumption that $v \in \overline{\class{G}}$.  If $e(r_0) < e(q_0) \btimes u \btimes -\bval{\dot{a}=\dot{c}}$ holds, then $-\bval{\dot{a}=\dot{c}} \in \overline{\class{G}}$ and so $i_\class{G}(\dot{a}) \neq i_\class{G}(\dot{c})$ by the Forcing Theorem applied to the formula $\dot{a}=\dot{c}$.  Since the choice of $c\in i_\class{G}(\classdot{C})$ was arbitrary, this shows that $i_\class{G}(\dot{a}) \notin i_\class{G}(\classdot{C})$.

For the reverse direction, suppose $\exists p \in \class{G} \, ({ p \forces \classdot{A} \subset \classdot{C} })$.  Fix $a \in i_\class{G}(\classdot{A})$ and $\<\dot{a},u> \in \classdot{A}$ such that $i_\class{G}(\dot{a})=a$ and $u\in \overline{\class{G}}$.  There is $q \leq p$ in $\class{G}$ and $\<\dot{c},v> \in \classdot{C}$ such that $e(q) \leq u \rightarrow \bval{\dot{a}=\dot{c}}\btimes v$.  Since $u\in \overline{\class{G}}$, I can assume without loss of generality that $e(q)\leq u$.  From the definition of $\rightarrow$, I have $e(q) \leq \({-u + \bval{\dot{a}=\dot{c}}\btimes v})$, and since $e(q)\leq u$ it follows that $e(q) \leq \bval{\dot{a}=\dot{c}}\btimes v$.  Thus both $\bval{\dot{a}=\dot{c}}$ and $v$ are in $\overline{\class{G}}$, and so $a=i_\class{G}(\dot{c})\in i_\class{G}(\dot{C})$.  Since this is true for every $a \in i_\class{G}(\classdot{A})$, I have $\V[\class{G}] \vDash i_\class{G}(\classdot{A}) \subset i_\class{G}(\classdot{C})$.

Continuing with the Forcing Theorem for atomic formulas containing class parameters, I now consider $\classdot{A}=\classdot{C}$  (Definition \ref{D:ForcingRelation} case \ref{D:ForcingRelation:ClassAtomic:equality}).
\begin{eqnarray*}
\lefteqn{\V[\class{G}] \vDash i_\class{G}(\classdot{A}) = i_\class{G}(\classdot{C})}\\
	& \iff & \V[\class{G}] \vDash i_\class{G}(\classdot{A}) \subset i_\class{G}(\classdot{C}) \mbox{ and } \V[\class{G}] \vDash i_\class{G}(\classdot{C}) \subset i_\class{G}(\classdot{A})\\
	& \iff & p\forces \classdot{A} \subset \classdot{C} \mbox{ and } p\forces \classdot{C} \subset \classdot{A}\\
	&  \iff & p\forces \classdot{A} = \classdot{C} 
\end{eqnarray*}

For case \ref{D:ForcingRelation:ClassAtomic:membership}, 
\begin{eqnarray*}
\lefteqn{i_\class{G}(\classdot{A}) \in i_\class{G}(\classdot{C})}\\
	& \iff & \exists \<\dot{c},v> \in \classdot{C} \, \({ v \in \overline{\class{G}} \mbox{ and } i_\class{G}(\classdot{A})=i_\class{G}(\dot{c}) })\\
	& \iff & \exists \<\dot{c},v> \in \classdot{C} \, \exists p \in \class{G} \, \({ e(p) \leq v  \mbox{ and } p \forces \classdot{A}=\dot{c} })\\
	& \rightarrow & q \forces \classdot{A} \in \classdot{C}.
\end{eqnarray*}
Note that the final step is an implication only.  For the reverse direction, it suffices to observe that if $q \forces \classdot{A} \in \classdot{C}$, then there is $q^\prime \leq q$ in $\class{G}$ such that $\exists \<\dot{c},v> \in \classdot{C} \, \({ v \in \overline{\class{G}} \mbox{ and } q^\prime \forces \classdot{A}=\dot{c} })$.  

This concludes the Forcing Theorem for atomic formulas.  I next consider the Boolean connectives.
\begin{eqnarray*}
	\V[\class{G}] \vDash \neg \phi & \iff & \V[\class{G}] \nvDash \phi\\
	& \iff & \neg \exists p \in \class{G} \, \({ p \forces \phi}) \mbox{ (inductive hypothesis)}\\
	& \iff & \exists p \in \class{G} \, \({ p \forces \neg \phi})
\end{eqnarray*}

For conjunctions,
\begin{eqnarray*}
	\lefteqn{\V[\class{G}] \vDash \phi\wedge\psi }\\
	& \iff & \V[\class{G}] \vDash \phi \mbox{ and } \V[\class{G}] \vDash \psi \\
	& \iff & \exists p \in \class{G} \, \({ p \forces \phi}) \mbox{ and } \exists q \in \class{G} \, \({ q \forces \psi})  \mbox{ (inductive hypothesis)}\\
	& \iff & \exists p \in \class{G} \, \({ p \forces \phi}) \mbox{ and } \({p \forces \psi})\mbox{ ($\class{G}$ is directed)}\\
	& \iff & \exists p \in \class{G} \, \({ p \forces \phi\wedge\psi}) 
\end{eqnarray*}

In the existential quantifier case, 
\begin{eqnarray*}
	\lefteqn{\V[\class{G}] \vDash \exists x \phi(x) }\\
	& \iff & \exists \class{A} \in \V[\class{G}] \mbox{ such that } \V[\class{G}]\vDash \phi(\class{A})\\
	& \iff & \exists \classdot{A} \in \V^\B \mbox{ such that } \V[\class{G}]\vDash \phi(i_\class{G}(\classdot{A}))\\
	& \iff & \exists \classdot{A} \in \V^\B \, \exists p \in \class{G} \, \({ p \forces \phi(\classdot{A}) }) \mbox{ (inductive hypothesis)}\\
	& \rightarrow & \exists p\in \class{G} \, \forall q\leq p \, \exists r \leq q \, \exists \classdot{A}\in \V^\B \, \({r\forces \phi(\classdot{A}) })\\
	& \iff & \exists p \in \class{G} \, p \forces \phi(\classdot{A})
\end{eqnarray*}
Note that the penultimate step is only an implication.  For the reverse direction, if $p \in \class{G}$ satisfies $\forall q\leq p \, \exists r \leq q \, \exists \classdot{A}\in \V^\B \, \({r\forces \phi(\classdot{A}) })$,
 then the set $\class{D} = \left\{ r \leq p \mid \exists \classdot{A}\in \V^\B \, \({r\forces \phi(\classdot{A}) }) \right\}$ is dense below $p$.  Fix $r \in \class{D}\cap \class{G}$ and $\classdot{A} \in \V^\B$ such that $r \forces \phi(\classdot{A})$.  Then by inductive assumption $\V[\class{G}] \vDash \phi(i_\class{G}(\classdot{A}))$ and so $\V[\class{G}] \vDash \exists x \phi(x)$.

Thus the Forcing Theorem holds for all $\phi$. \end{proof}

\section{A partial Generic Model Theorem}\label{SS:PartialGenericModelTheorem}

The main goal of this development is the Generic Model Theorem, which states among other things that the model $\<V[\class{G}],\V[\class{G}]>$ satisfies $\BGC$.  Unfortunately the theorem does not hold in general for all partial orders $\P$, but in the case of $P$ a set we have:

\begin{Thm}\label{T:GenericModelTheoremSets}\emph{Generic Model Theorem for set forcing}.  Suppose $\<V,\V>$ is a model of $\BGC$ and $P$ is a poset.  If $G\subset P$ is $\V$-generic for $P$, then
\begin{enumerate}
\parskip=0pt
\item  If $\<W,\W>$ is a transitive model of $\BGC$ such that $\V \subset \W$ and $G \in\W$, then $\V[G]\subset \W$ and $V[G] \subset W$.
\item  $G \in \V[G]$ and $\V \subset \V[G]$.
\item  $\ORD^\V = \ORD^{\V[G]}$.
\item  $\<V[G],\V[G]>$ satisfies $\BGC$.
\end{enumerate}
\end{Thm}

As the main point of this exposition is the corresponding theorem for class $\P$, the above Theorem will not be proved.  Although the proof is instructive, especially as an exercise in working in $\BGC$, it is basically a matter of adapting the usual proof for $\ZFC$.  In the case of $\P$ a proper class, the assumption that $\P=\bigcup P_\alpha$ is a chain of complete subposets will guarantee that the resulting $\<V[\class{G}],\V[\class{G}]>$ satisfies a fragment of $\BGC$ (Theorem \ref{T:PartialGenericModelTheorem}), but even stronger restrictions will be required for the full result.  These stronger restrictions will split into two cases, one dealing with iterations (Section \ref{SS:GMTIterations}) and one with products (Section \ref{SS:GMTProducts}).  These cases are general enough to include a great number of common applications of class forcing, including those presented in this paper.  The specific restrictions will be stated in Definitions \ref{D:IterationRestrictions} and \ref{D:ProductRestrictions}, and it will be shown that a broad class of iterations and products satisfy these restrictions.  Finally, the full Generic Model Theorem will be proved in each of these cases.

\begin{Thm}\label{T:PartialGenericModelTheorem}\emph{Generic Model Theorem (partial)}.  Suppose $\<V,\V>$ is a model of $\BGC$, the partial order $\P=\bigcup P_\alpha$ is a chain of complete subposets, and $\class{G}\subset \P$ is $\V$-generic for $\P$.
\begin{enumerate}
\parskip=0pt
\item \label{T:PartialGenericModelTheorem:Smallest} If $\<W,\W>$ is a transitive model of $\BGC$ such that $\V \subset \W$ and $\class{G}\in\W$, then $\V[\class{G}]\subset \W$ and $V[\class{G}] \subset W$.
\item \label{T:PartialGenericModelTheorem:ContainsV} $\class{G} \in \V[\class{G}]$ and $\V \subset \V[\class{G}]$.
\item \label{T:PartialGenericModelTheorem:SameOrdinals} $\ORD^\V = \ORD^{\V[\class{G}]}$.
\item \label{T:PartialGenericModelTheorem:PartialZFC} $\<V[\class{G}],\V[\class{G}]>$ satisfies all axioms of $\BGC$ except possibly power set and replacement.
\end{enumerate}
\end{Thm}

\begin{proof}   Suppose $\<W,\W>$ satisfies the hypotheses of assertion \ref{T:PartialGenericModelTheorem:Smallest}.  Fix $\class{A} \in \V[\class{G}]$.  There is $\classdot{A} \in \V^\P$ such that $i_G(\classdot{A})=\class{A}$, and by hypothesis $\classdot{A}\in \W$.  As $\W$ is transitive, I also have that $\tc(\{\classdot{A}\})\subset \W$.   The valuation function $i_\class{G}$ depends only on $\class{G}$ and recursively on members of $\tc(\classdot{A})$, so $\W$ will calculate $i_{\class{G}}(\classdot{A})$ correctly.  Thus $\class{A} \in \W$, and so $\V[\class{G}]\subset \W$.  That  $\class{G} \in \V[\class{G}]$ and $\V\subset \V[\class{G}]$ and are given by Lemmas \ref{L:Gdot} and  \ref{L:checknames}.  Comparing the ordinals of $\V$ and $\V[\class{G}]$, it is clear that $\ORD^\V \subset \ORD^{\V[\class{G}]}$.  For the reverse inclusion, note that by  Lemma \ref{L:V[G]rank} every set $\dot{a}\in V^\P$ has $\B$-rank $\rho(\dot{a})$ greater than or equal to the rank of its valuation $i_G(\dot{a})$.  Thus there can be no $\dot{a}$ whose valuation yields an ordinal greater than all the ordinals of $\V$.

It remains to show that $\<V[\class{G}],\V[\class{G}]>$ satisfies the fragment of $\BGC$ given in assertion \ref{T:PartialGenericModelTheorem:PartialZFC}.  In the following, $\B$ will be the Boolean algebra associated with $\P$, and $\overline{\class{G}}$ will be the $\V$-generic ultrafilter on $\B$ associated with $\class{G}$.  That $\V[\class{G}]$ is \emph{extensional} is a consequence of transitivity, Lemma \ref{T:V[G]BasicFacts}.  That every set is a class follows from the observation that $V^\P \subset \V^\P$, and so $V[\class{G}] \subset \V[\class{G}]$.  It was shown in Lemma \ref{T:V[G]BasicFacts} that $\class{X} \in \class{Y}$ implies $\class{X}$ is a set.  The \emph{pairing axiom} asserts
$$\forall x,y \, \exists z \, \forall u \, \({u\in z \iff \({u=x \vee u=y}) }).$$
Fix $a$ and $b$ in $V[G]$ with names $\dot{a}$ and $\dot{b}$.  Setting $\dot{c} = \{\<\dot{a},1>, \<\dot{b},1>\}$, I have $i_G(\dot{c}) = \{i_G(\dot{a}), i_G(\dot{b})\} = \{a,b\}$.  

The \emph{axiom of comprehension} asserts that for any formula $\phi(x, \class{X}_1, ..., \class{X}_n)$ in which only set variables are quantified and for any classes $\class{A}_1, ..., \class{A}_n$, the collection $\class{Y} = \{x \mid \phi(x,\class{A}_1, ..., \class{A}_n) \}$ is a class.  Fix such a formula $\phi$.  The parameters $\class{A}_1, ..., \class{A}_n$ will be suppressed, as the proof is essentially the same.  Working in $\V$, let $\classdot{Y} = \{\<\dot{a},e(p)> \mid  p \forces \phi(\dot{a}) \}$.  I claim that $i_\class{G}(\classdot{Y}) = \{ x \mid \V[\class{G}] \vDash \phi(x) \}$.  For the forward inclusion, fix $a \in i_\class{G}(\classdot{Y})$.  There must be $\<\dot{a},u> \in \classdot{Y}$ such that $u\in \overline{G}$ and $i_\class{G}(\dot{a})=a$.  From the definition of $\classdot{Y}$ it must be the case that $u=e(p)$ and $p \forces \phi(a)$.  Since $e(p)\in \overline{\class{G}}$ it follows that $p\in \class{G}$, and so by the Forcing Theorem $\V[\class{G}] \vDash \phi(i_\class{G}(\dot{a}))$, or $\V[\class{G}] \vDash \phi(a)$.  For the reverse inclusion, suppose $\V[\class{G}] \vDash \phi(a)$ for some $a \in V[\class{G}]$.  Fix a name $\dot{a}$ for $a$.  There must be $p\in \class{G}$ such that $p \forces \phi(\dot{a})$, and so $\<\dot{a},e(p)> \in \classdot{Y}$.  Thus $i_\class{G}(\dot{a}) \in i_\class{G}(\classdot{Y})$.

The \emph{axiom of infinity}, asserting the existence of an inductive set, is witnessed in $\V[\class{G}]$ by any inductive set from $V$, as $V \subset \V[\class{G}]$.  

The \emph{union axiom} will require a little more work.  I would like to show
$$\forall x \, \exists y \, \forall z \, \({z\in y \iff \exists u \in x \, \({z\in u})}).$$
Fix $a\in V[\class{G}]$ with name $\dot{a}$, and define $\dot{b} = \{\<\sigma,u_\sigma>  \mid  \exists \tau \in \dom(\dot{a}), \,  \sigma \in \dom(\tau) \}$, where $u_\sigma  = \sum_{\tau\in \dom(\dot{a})} \{\dot{a}(\tau)\btimes\tau(\sigma) \mid  \sigma \in \dom(\tau) \}$.  I claim that $i_G(\dot{b}) = \bigcup i_G(\dot{a})$.  Fix $x\in i_G(\dot{b})$.  There is $\sigma \in \dom(\dot{b})$ such that $\dot{b}(\sigma) \in \overline{\class{G}}$ and $i_\class{G}(\sigma)=x$.
\begin{eqnarray*}
	 \dot{b}(\sigma)\in \overline{\class{G}} & \iff & u_\sigma \in \overline{\class{G}}\\
	& \iff & \sum_{\tau\in \dom(\dot{a})} \{\dot{a}(\tau)\btimes\tau(\sigma)  \mid   \sigma \in \dom(\tau) \} \in \overline{\class{G}}\\
	& \iff & \exists \tau \in \dom(\dot{a}), \, \sigma \in \dom(\tau) \, \({\dot{a}(\tau)\btimes\tau(\sigma) \in \overline{\class{G}} }) \mbox{ (genericity)}\\
	& \rightarrow & \dot{a}(\tau) \in \overline{\class{G}} \mbox{ and } \tau(\sigma) \in \overline{\class{G}}\\
	& \rightarrow & i_\class{G}(\sigma)\in i_\class{G}(\tau) \in i_\class{G}(\dot{a}) \mbox{ (Forcing Theorem)}
\end{eqnarray*}
The third bi-implication follows from the fact that if a sum is in a generic ultrafilter, one of the summands must be as well.  This is a straightforward consequence of the definition of generic ultrafilter on a Boolean algebra, Definition \ref{D:BAfilters}, applying the identity $\sum\class{X}=-\prod-\class{X}$.
This shows that $i_\class{G}(\dot{b}) \subset \bigcup a$.  For the reverse inclusion, suppose $y \in \bigcup a$.  Then there is $\tau\in \dom(\dot{a})$ and $\sigma \in \dom(\tau)$ such that $\dot{a}(\tau), \tau(\sigma) \in \overline{\class{G}}$ and $i_\class{G}(\sigma)=y$.
\begin{eqnarray*}
	\dot{a}(\tau), \tau(\sigma) \in \overline{\class{G}} & \iff & \dot{a}(\tau)\btimes\tau(\sigma) \in \overline{\class{G}}\\
	& \rightarrow & u_\sigma \in \overline{\class{G}}\\
	& \rightarrow & i_\class{G}(\sigma) \in i_\class{G}(\dot{b})\\
	& \rightarrow & y \in i_\class{G}(\dot{b})
\end{eqnarray*}
Thus $i_\class{G}(\dot{b}) = \bigcup a$, and so $\V[\class{G}]$ satisfies the union axiom.

The \emph{axiom of regularity} states that every nonempty class has an $\in$-minimal element.
$$\forall \class{X} \, \({ \class{X} \ne \emptyset \rightarrow \exists y\in \class{X} \, \(y\cap \class{X} = \emptyset) }).$$
Lemma \ref{L:V[G]rank} shows that every set in $V[\class{G}]$ has rank.  Fix $\class{X}\neq \emptyset$ and let $\alpha$ be least such that there is $y\in \class{X}$ of rank $\alpha$.  Then $y\cap \class{X} = \emptyset$, since $z\in y \rightarrow \rank(z) < \rank(y)$, and so $z \in \class{X}$ would contradict `$\alpha$ is least.'  Thus $\V[\class{G}]$ satisfies regularity.  

For the \emph{axiom of choice}, asserting 
$$\exists \class{F} \mbox{ a function } \forall x \neq \emptyset \, (\class{F}(x)\in x),$$
I will first show that there is a \emph{definable} choice function in $\V[\class{G}]$.  Note that the choice function $\class{F} \in \V$ can be used to well-order any set $A\in V$ in the usual fashion, by first defining a sequence $A_0 = A$, at successors $A_{\alpha+1} = A_\alpha \setminus \{\class{F}(A_\alpha)\}$, and at limits $A_\lambda = \cap_{\alpha<\lambda} A_\alpha$, and defining for all $a,b \in A$ the ordering $a<b$ if and only if there is an $\alpha$ for which $b\in A_\alpha$ but $a \notin A_\alpha$.  From this one can define a well-ordering $\prec$ of all sets $V$, ordering first according to rank and then using $\class{F}$ to order the sets of rank $\alpha$ for each $\alpha$.  To define a choice function $\phi(x,y)$ on sets in $\V[\class{G}]$, first note that both $V^\P$ and $\prec$ are classes in $\V[\class{G}]$.  For each set $A \in V[\class{G}]$ let $\dot{A}$ be the $\prec$-least member of $V^\P$ such that $i_\class{G}(\dot{A}) = A$, and let $\dot{a}$ be the $\prec$-least member of $\dom(\dot{A})$ such that $\dot{A}(\dot{a})\in \overline{\class{G}}$.  Then setting $x_A = i_\class{G}(\dot{a})$, the function $\phi(x,y)$ selecting $x_A$ from $A$, i.e. $\V[\class{G}] \vDash \forall A, x, \, (\phi(A,x) \iff x=x_A)$, is definable in $\V[\class{G}]$ from parameters $\class{F}$ and $V^\P$ together with the generic $\class{G}$.  As the axiom of comprehension holds in $\V[\class{G}]$, it follows that the function defined by $\phi$ is a class and so $\V[\class{G}]$ satisfies the axiom of choice.
\end{proof}

\section[The Generic Model Theorem for iterations]{The Generic Model Theorem \\ for iterations}\label{SS:GMTIterations}

The general theory of iterations will be developed in section \ref{SS:Iterations}, and readers may wish to review the definitions in that section before proceeding.  The following condition is sufficient to allow the Generic Model Theorem to go through.  It will be shown in Section \ref{SS:Iterations} that the iterations used in Chapters \ref{C:TheGroundAxiom} and \ref{C:GAandVnotHOD} satisfy this condition.

\clearpage
\begin{Def}\rm\label{D:IterationRestrictions}$\P= \bigcup P_\alpha$ is a \emph{progressively closed iteration} if and only if $\P$ is a chain of complete subposets and for arbitrarily large regular cardinals $\delta$ there are arbitrarily large $\alpha$ such that there is a $P_\alpha$ name $\classdot{P}_{[\alpha,\infty)} = \{ \< \op(\check{\beta}, \dot{P}_{[\alpha,\beta)} ), 1> \mid \beta > \alpha \}$ satisfying 
\begin{enumerate}
\parskip=0pt
\item for every $\beta > \alpha$ the poset $P_\beta$ is isomorphic to the two-stage iteration $P_\beta \cong P_\alpha * \dot{P}_{[\alpha,\beta)}$,
\item $P_\alpha \forces \check{\delta} \mbox{ is a regular cardinal and }\dot{P}_{[\alpha,\beta)} \mbox{ is } <\!\check{\delta} \mbox{-closed}$,
\item for $\beta^\prime > \beta > \alpha$ the isomorphisms at $\beta$ and $\beta^\prime$ yield complete subposets $P_\alpha * \dot{P}_{[\alpha,\beta)} \subset_c P_\alpha * \dot{P}_{[\alpha,\beta^\prime)}$ such that the following diagram commutes:
 \[
 \begin{CD}
P_\beta @>\subset_c>> P_{\beta^\prime} \\
@V{\cong}VV @V{\cong}VV\\
P_\alpha * \dot{P}_{[\alpha,\beta)} @> \subset_c>> P_\alpha * \dot{P}_{[\alpha,\beta^\prime)} \\
 \end{CD}
 \]
\item $P_\alpha \forces \classdot{P}_{[\alpha,\infty)}$ is a chain of complete subposets, that is, $P_\alpha \forces \dot{P}_{[\alpha,\beta)} \subset_c  \dot{P}_{[\alpha,\beta^\prime)}$ for $\beta^\prime > \beta > \alpha$ .
\end{enumerate}
In such a case, I will say \emph{$\P$ factors at $\alpha$ with closure $\delta$.}

\end{Def}

Note that the definition of $\classdot{P}_{[\alpha,\infty)}$ relies on the `ordered pair' function, Definition \ref{D:op}. The essence of Definition \ref{D:IterationRestrictions} is that as you progress through the iteration, the remaining forcing becomes more and more closed.  The following lemma shows that if we force with an initial part $P_\alpha$, then in that partial extension the `remaining part of $\P$' still forms a progressively closed iteration.  This fact will be important as it will allow us to carry out forcing arguments within some $\V[G_\alpha]$ while preserving hypotheses on $\P$.

\begin{Lemma}\label{L:FactoringLemma}\emph{Factoring Lemma.} Suppose $\P = \bigcup P_\alpha$ is a progressively closed iteration and $\class{G}$ is $\V$-generic for $\P$. Then for every ordinal $\alpha$ and regular cardinal $\delta$ such that $\P$ factors at $\alpha$ with closure $\delta$, setting  $G_\alpha = \class{G} \cap P_\alpha$,
\begin{enumerate}
\parskip=0pt

\item \label{L:FactoringLemma:PCI} $\P_{[\alpha,\infty)} = i_{G_\alpha}(\classdot{P}_{[\alpha,\infty)})$ is a progressively closed iteration in $\V[G_\alpha]$.

\item \label{L:FactoringLemma:closure} $V[G_\alpha] \vDash \P_{[\alpha,\infty)}$ is $< \delta$-closed.

\setcounter{saveenum}{\value{enumi}}
\end{enumerate}
Suppose further that for each $\beta > \alpha$, the set $G_{[\alpha,\beta)}$ is the $V[G_\alpha]$-generic subset of  $i_{G_\alpha}\({\dot{P}_{[\alpha,\beta)} })$ such that $G_\beta = G_\alpha *  G_{[\alpha,\beta)}$ (as in Lemma \ref{L:P*Q}), and $\class{G}_{[\alpha,\infty)} = \bigcup G_{[\alpha,\beta)}$.  Then 
\begin{enumerate}
\parskip=0pt
\setcounter{enumi}{\value{saveenum}}

\item \label{L:FactoringLemma:TailGeneric}$\class{G}_{[\alpha,\infty)}$ is $V[G_\alpha]$-generic for $\P_{[\alpha, \infty)}$. 

\item \label{L:FactoringLemma:FactorV[G]}$\V[\class{G}]=\V[G_\alpha][\class{G}_{[\alpha,\infty)}]$.

\item \label{L:FactoringLemma:ForcingExtension}Thus, $\V[\class{G}]$ is a forcing extension of $\V[G_\alpha]$ by $\P_{[\alpha, \infty)}$.

\end{enumerate}
\end{Lemma}

\begin{proof}  Note the abuse of notation in which the sequence given by $\P_{[\alpha,\infty)} = \{\<\beta,i_{G_\alpha}(\dot{P}_{[\alpha,\beta)})> \mid \beta>\alpha\}$ is confused with the actual partial order formed by $\bigcup \ran(\P_{[\alpha,\infty)})$.  To see that $\P_{[\alpha, \infty)}$ is a chain of complete subposets, observe that for $\beta < \beta^\prime$ I have $P_\alpha \forces \dot{P}_{[\alpha,\beta)} \subset_c  \dot{P}_{[\alpha,\beta^\prime)}$ and so $V[G_\alpha] \vDash i_{G_\alpha}({\dot{P}_{[\alpha,\beta)} }) \subset_c  i_{G_\alpha}({ \dot{P}_{[\alpha,\beta^\prime)} })$.  For the progressive closure condition, fix $\kappa$ and fix $\alpha^\prime > 
\alpha$ such that $\P$ factors at $\alpha^\prime$ with closure $\kappa$.  Since I can factor $P_{\alpha^\prime} \cong P_\alpha * \dot{P}_{[\alpha, \alpha^\prime)}$, and for $\beta > \alpha^\prime$ I can also factor $P_\beta \cong P_{\alpha^\prime} * \dot{P}_{[\alpha^\prime,\beta)}$, it follows that $P_\beta \cong \({ P_\alpha * \dot{P}_{[\alpha, \alpha^\prime)} }) * \dot{P}_{[\alpha^\prime,\beta)}^{\star}$.  The final factor $\dot{P}_{[\alpha^\prime,\beta)}^{\star}$ is not strictly equal to $\dot{P}_{[\alpha^\prime,\beta)}$, as the former is a $P_\alpha * \dot{P}_{[\alpha, \alpha^\prime)}$-name and the latter is a $P_{\alpha^\prime}$-name.  In fact, the isomorphism $P_{\alpha^\prime} \cong P_\alpha * \dot{P}_{[\alpha, \alpha^\prime)}$ induces an isomorphism on names, so to each $P_{\alpha^\prime}$ name $\tau$ there is associated an isomorphic $P_\alpha * \dot{P}_{[\alpha, \alpha^\prime)}$ name $\tau^\star$, and it is this name that is used as the third factor.  This and related properties of iterations will be explored in Section \ref{SS:Iterations}.  It is a  matter of inductively redefining names to re-associate the three factors as $P_\beta \cong  P_\alpha * \({ \dot{P}_{[\alpha, \alpha^\prime)}  * \dot{P}_{[\alpha^\prime,\beta)}^{\star\star} })$ in such a way that $P_\alpha \forces \({ \dot{P}_{[\alpha, \alpha^\prime)}  * \dot{P}_{[\alpha^\prime,\beta)}^{\star\star} }) \cong \dot{P}_{[\alpha,\beta)}$.  The precise nature of this transformation is given in Definition \ref{D:ReassociateP*Q}.  Once again, $\dot{P}_{[\alpha^\prime,\beta)}^{\star\star}$is not exactly $\dot{P}_{[\alpha^\prime,\beta)}^\star$, but is in fact a $P_\alpha$ name for a $\dot{P}_{[\alpha, \alpha^\prime)}$ name for $\dot{P}_{[\alpha^\prime,\beta)}^\star$.  Thus $\V[G_\alpha] \vDash i_{G_\alpha}\({\dot{P}_{[\alpha,\beta)} }) \cong i_{G_\alpha}\({ \dot{P}_{[\alpha, \alpha^\prime)}  * \dot{P}_{[\alpha^\prime,\beta)}^{\star\star} }) \cong i_{G_\alpha}\({ \dot{P}_{[\alpha, \alpha^\prime)} })  * i_{G_\alpha}\({ \dot{P}_{[\alpha^\prime,\beta)}^{\star\star} })$, and $\V[G_\alpha] \vDash i_{G_\alpha}\({ \dot{P}_{[\alpha, \alpha^\prime)} }) \forces i_{G_\alpha}\({ \dot{P}_{[\alpha^\prime,\beta)}^{\star\star} })$ is $<\!\check{\kappa}$-closed.  A consequence of Definition \ref{D:ReassociateP*Q} and Lemma \ref{L:ReassociateP*Q} is that in $\V[G_\alpha]$ the isomorphisms $i_{G_\alpha}\({\dot{P}_{[\alpha,\beta)} }) \cong i_{G_\alpha}\({ \dot{P}_{[\alpha, \alpha^\prime)} })  * i_{G_\alpha}\({ \dot{P}_{[\alpha^\prime,\beta)}^{\star\star} })$ will commute with the embeddings $i_{G_\alpha}\({\dot{P}_{[\alpha,\beta)} }) \subset_c i_{G_\alpha}\({\dot{P}_{[\alpha,\beta^\prime)} })$.
Another consequence is that for $\alpha < \alpha^\prime < \beta < \beta^\prime$,
$$V[G_\alpha] \vDash i_{G_\alpha}\({ \dot{P}_{[\alpha,\alpha^\prime)} }) \forces i_{G_\alpha}\({ \dot{P}_{[\alpha^\prime,\beta)}^{\star\star} }) \subset_c i_{G_\alpha}\({ \dot{P}_{[\alpha^\prime,\beta^\prime)}^{\star\star} }).$$ 
It follows that in $\V[G_\alpha]$ the assertion that $\P_{[\alpha,\infty)}$ factors at $\alpha^\prime$ with closure $\kappa$ is witnessed by the $i_{G_\alpha}\({ \dot{P}_{[\alpha, \alpha^\prime)} })$-name $\Q=\left\{\op\(\check{\beta},i_{G_\alpha}({ \dot{P}_{[\alpha^\prime,\beta)}^{\star\star} ) }) \mid \beta > \alpha^\prime \right\}$.

The closure of $\P_{[\alpha, \infty)}$ in $V[G_\alpha]$ follows immediately from the closure of each $i_{G_\alpha}\({\dot{P}_{[\alpha,\beta)} })$.  The final three points are consequences of the partial Generic Model Theorem (Theorem \ref{T:PartialGenericModelTheorem}) and facts about iterations (Section \ref{SS:Iterations}, Lemmas \ref{L:P*Q} and \ref{L:ReassociateP*Q}).
\end{proof}

In forcing arguments it is common to fix an $\alpha$ and consider $\P$ to be the forcing $P_\alpha$, in which $G_\alpha$ is adjoined to the universe, followed by the progressively closed iteration $\bigcup i_{G_\alpha}\({ \dot{P}_{[\alpha,\beta)} })$ in $V[G_\alpha]$.  In such a case, the latter forcing is refered to as the `tail forcing' or the `tail of the iteration'.  Thus the usual application of progressive closure is to work in $V[G_\alpha]$ for $\alpha$ sufficiently large that the tail forcing is highly closed, and so adds no short sequences over the ground model.  This consequence of closure is presented in any exposition of forcing and will be stated here but not proved.

\begin{Lemma}\label{L:ClosedForcing}Suppose $\kappa$ is a regular cardinal and $\P$ is a $<\!\kappa$-closed partial order.  If $\class{G}$ is $\V$-generic for $\P$, then for every $\delta < \kappa$, for every class $\class{J}\in V[\class{G}]$ such that $\class{J}:\delta \rightarrow V$, the function $\class{J}$ is a set in the ground model $\V$.
\end{Lemma}

The full Generic Model Theorem holds for progressively closed iterations.

\begin{Thm}\label{T:GenericModelTheoremIterations}\emph{Generic Model Theorem (iterations)}.  Suppose $\P = \bigcup P_\alpha$ is a progressively closed iteration and $\class{G}\subset \P$ is $\V$-generic for $\P$.  Then
$\V[\class{G}] \vDash \BGC$.
\end{Thm}

\begin{proof}Note that the partial Generic Model Theorem (Theorem \ref{T:PartialGenericModelTheorem}) applies and so $\V[\class{G}]$ immediately satisfies a fragment of $\BGC$.  It remains to show $\V[\class{G}]$ satisfies the axioms of power set and replacement.

The \emph{power set} axiom asserts
$$\forall x \, \exists y \, \forall \class{Z} \, \({ \class{Z} \subset x \iff \class{Z} \in y }).$$
Fix $a\in \V[\class{G}]$.  Fix $\alpha$ sufficiently large that $a\in V[G_\alpha]$ and $\P$ factors at $\alpha$ with closure $\({ |a|^+ })^{\V[G_\alpha]}$.  Notice that every subclass of $a$ in $\V[\class{G}]$ is already in $\V[G_\alpha]$ by the closure of the tail forcing (Lemma \ref{L:ClosedForcing}), so it suffices to show that $\mathcal{P}(a)$ exists in $\V[G_\alpha]$.  But this is true, since $\V[G_\alpha]\vDash \BGC$ by the Generic Model Theorem for sets (Theorem \ref{T:GenericModelTheoremSets}).  Thus $\V[\class{G}]$ satisfies the power set axiom.

The \emph{axiom of replacement} asserts
$$\forall \class{J} \, \forall a \, \({  \class{J} \mbox{ is a function} \rightarrow \exists z \, \({ z = \{ \class{J}(x) \mid x \in a \}  })  }).$$
Note that since the axiom of choice gives a bijection of every set with an ordinal, it is sufficient to consider domains $a = \beta$ an ordinal.  Fix a function $\class{J}\in \V[\class{G}]$ and an ordinal $\beta$.  Let $\alpha$ be sufficiently large that $\P$ factors at $\alpha$ with closure $(|\beta|^+)^{\V[\class{G}]}$.  The remainder of the proof will be carried out in $\V[G_\alpha]$.  Fix $\classdot{J}$ a $\P_{[\alpha,\infty)}$-name for $\class{J}$.

\emph{Claim.} For any $p \in \P_{[\alpha,\infty)}$ such that $p\forces \class{J}$ is a function, there is a condition $q \leq p$ and a $\beta$-sequence $\<\dot{a}_\eta \mid \eta < \beta>$ of $\P_{[\alpha,\infty)}$-names such that $q \forces \classdot{J}(\check{\eta})=\dot{a}_\eta$  for every $\eta < \beta$.  

If I can establish this claim, then it follows directly that such $q$ are dense below such $p$.  Since one such $p$ lies in $\class{G}_{[\alpha,\infty)}$ by hypothesis, there must be such a $q \in G_{[\alpha,\infty)}$.  For this $q$, consider the associated sequence $\<\dot{a}_\eta \mid \eta < \beta>$ and define a name $\dot{b} = \{\<\dot{a}_\eta, 1>  \mid  \eta< \beta\}$.  Clearly $\dot{b} \in \V[G_\alpha]$, and $i_{\class{G}_{[\alpha,\infty)}} (\dot{b})$ is exactly equal to the image of $\beta$ under $\classdot{J}$ in $\V[\class{G}]$, since $q \forces \classdot{J}(\check{\eta})=\dot{a}_\eta$ for every $\eta < \beta$.  Thus $\V[\class{G}]$ will satisfy Replacement.

It remains to prove the claim. Fix such a $p$.  For each $\eta < \beta$, note that $p \forces \exists x \, \classdot{J}(\check{\eta})=x$.  Consider the class $D_\eta = \{q \leq p  \mid  \exists \dot{a} \, p\forces \classdot{J}(\check{\eta})=\dot{a}\}$.  Each $D_\eta$ is dense in $\P_{[\alpha,\infty)}$ below $p$. Now recursively construct two sequences, a descending sequence $\<q_\eta \mid \eta<\beta>$ in $\P_{[\alpha,\infty)}$ and a sequence $\<\dot{a}_\eta>$ of $\P_{[\alpha,\infty)}$-names, satisfying  $q_\eta \forces \classdot{J}(\check{\eta}) = \dot{a}_\eta$ for each $\eta<\beta$.  Suppose the sequences have been defined up to $\eta$.  Closure of $\P_{[\alpha,\infty)}$ allows me to find a $q$ lying below $q_\delta$ for all $\delta < \eta$.  Density of $D_\eta$ allows me to choose $q_\eta \in D_\eta$ lying below $q$.  Since $q_\eta \in D_\eta$ there is a corresponding $\dot{a}_\eta$ such that $q_\eta \forces \classdot{J}(\check{\eta}) = \dot{a}_\eta$.  Thus the recursion continues for all $\eta < \beta$.  Finally, closure allows me to fix $q$ lying below all $q_\eta$.  This completes the proof of the claim, and shows that $\V[\class{G}]$ satisfies the axiom of replacement.  Thus $\V[\class{G}] \vDash \BGC$.
\end{proof}

\section[The Generic Model Theorem for products]{The Generic Model Theorem \\ for products}\label{SS:GMTProducts}

For general remarks on product forcing, see Section \ref{SS:Products}.  For the purposes of establishing $\BGC$ in the extension, the following condition is sufficient.

\begin{Def}\rm\label{D:ProductRestrictions}A partial order $\P=\bigcup P_\alpha$ is a \emph{progressively closed product} if and only if $\P$ is a chain of complete subposets and for arbitrarily large regular $\delta$, $\P$ can be written as a product $$\P \cong \P_1 \times P_2$$
such that $P_2$ is a set, $P_2$ has the $\delta^+$-c.c., $\P_1$ is a chain of complete subposets, and $\P_1$ is $\leq\!\delta$-closed.  For such $\delta$, I will say \emph{$\P$ factors at $\delta$}.
\end{Def}

The important fact about products such as $\P_1 \times P_2$ is given in the following lemma.  The techniques used in the proof of the lemma will be employed in the proof of the Generic Model Theorem for products, Theorem \ref{T:GenericModelTheoremProducts}.

\begin{Lemma}\label{L:ClosedProducts} Suppose $\delta$ is a regular cardinal and $\class{G}\times H$ is $\V$-generic for $\P_1\times P_2$ such that $P_2$ is a set, $P_2$ has the $\delta^+$-c.c., and $\P_1$ is $\leq\!\delta$-closed.  Then any class $\class{J} \in \V[\class{G}\times H]$  which is a function $\class{J}:\delta \rightarrow V$ is a set in $\V[H]$.
\end{Lemma}

\begin{proof} The fact that $V[\class{G}\times H] = V[\class{G}][H] = V[H][\class{G}]$ will be discussed in Section \ref{SS:Products}, Lemma \ref{L:PxQclasses}.

Fix $\class{J} \in V[\class{G}\times H]$ as in the lemma and a $\P_1\times P_2$-name $\classdot{J}$ for $\class{J}$.  Consider any condition $\<g,h> \in \P_1\times P_2$ such that $\<g,h> \forces \classdot{J}:\check{\delta}\rightarrow \check{V}$, and fix $\alpha < \delta$.  Clearly I can extend $\<g,h>$ to $\<p,q>$ such that there is $a\in V$ and $\<p,q> \forces \classdot{J}(\check{\alpha}) = \check{a}$. However, I claim that even more is true.

\emph{Claim}.  If $\<g,h> \forces \classdot{J}:\check{\delta}\rightarrow \check{V}$ and $\alpha < \delta$, then there is $p\leq g$ and $A \subset P_2$ such that $A$ is a maximal antichain below $h$ and for every $q\in A$ there is $a_q \in V$ such that $\<p,q> \forces \classdot{J}(\check{\alpha}) = \check{a}_q$.

To prove the claim, recursively construct a descending sequence $\<p_\beta>$ below $g$, a sequence of pairwise incompatible elements $\<q_\beta>$ below $h$, and a sequence of sets $\<a_\beta>$ such that $\<p_\beta,q_\beta> \forces \classdot{J}(\check{\alpha}) = \check{a}_\beta$ for every $\beta$.  Suppose the sequences have been constructed below $\beta$.  If $\{q_\eta \mid \eta < \beta \}$ forms a maximal antichain below $h$, then the construction halts.  Otherwise, there is $q < h$ that is pairwise incompatible with $q_\eta$ for every $\eta < \beta$.  Closure of $\P_1$ allows me to find $p$ lying below all $p_\eta$ for $\eta < \beta$ provided $\beta < \delta^+$.  Since $\<p,q> \leq \<g,h>$ I can extend $\<p,q>$ to a condition $\<p_\beta,q_\beta>$ such that there is $a_\beta \in V$ and $\<p_\beta,q_\beta> \forces \classdot{J}(\check{\alpha}) = \check{a}_\beta$.  Note that the construction must halt at some stage below $\delta^+$, as the $q_\beta$ form an antichain and $P_2$ has the $\delta^+$-c.c..  After the construction halts, I can once again use closure in $\P_1$ to obtain $q$ lying below all $q_\beta$.  This $q$, together with $A = \{p_\beta\}$ and the associated $a_\beta$, witness the properties of the claim.

This argument shows that for $\<g,h>$ as in the claim and for every $\alpha < \delta$, the following class $D_\alpha$ is open and dense below $g$.
\begin{eqnarray*}
	D_\alpha = \left\{   p \mid \exists A^{\alpha}_p \subset P_2 \mbox{ max antichain below }h \mbox{ such that }  \vphantom{\classdot{J}}\right.\\ 
	\left. \forall q \in A^{\alpha}_p \, \exists a^{\alpha}_{p,q} \, \({ \<p,q> \forces \classdot{J}(\check{\alpha}) = \check{a}^\alpha_{p,q} }) \right\} 
\end{eqnarray*}
In fact, the $\delta^+$-closure of $\P_1$ ensures that the intersection of $\delta$ many open dense sets is dense, and so $\cap_{\alpha<\delta} D_\alpha$ is dense below $g$.  To prove this fact, for any $p < g$ construct a descending sequence $\<p_\alpha>$ below $p$, meeting $D_\alpha$ at stage $\alpha$ and using closure of $\P_1$ to get through limit stages.  

Since it is true by assumption that $\class{J}:\delta \rightarrow V$ in $\V[\class{G}\times H]$, there must be $\<g,h>\in \class{G}\times H$ as in the claim.  It follows from the argument above that there is $\<g^\prime, h^\prime> \in \class{G}\times H$ such that $\<g^\prime, h^\prime> \forces \classdot{J} $ is a function from $\check{\delta}$ to $\check{V}$ and $g^\prime \in D_\alpha$ for every $\alpha < \delta$. Fix such a $\<g^\prime, h^\prime>$.  I now work in $V[H]$.  Note that $H$ is $\V$-generic for $P_2$ and $h^\prime \in H$, so $H$ must contain a unique element of each maximal antichain $A^{\alpha}_{g^\prime}$ below $h^\prime$.  As $\V[H]$ is a generic extension of $\V$ by set forcing I have $\V \subset \V[H]$, so the forcing relation $\forces_{\P_1 \times P_2}$, the collection $\{ A^{\alpha}_{g^\prime} \mid \alpha < \delta\}$, and the corresponding $\{a^{\alpha}_{g^\prime, h} \mid \alpha < \delta \mbox{ and } h \in  A^{\alpha}_{g^\prime} \}$ are all in $\V[H]$.  Therefore, $\V[H]$ can define $\class{J}$ by $\class{J}(\alpha) = $ the unique $a^{\alpha}_{g^\prime, h} \in A^{\alpha}_{g^\prime}$ such that $h \in A^{\alpha}_{g^\prime}$.  Thus $\class{J} \in \V[H]$.  
\end{proof}

\begin{Thm}\label{T:GenericModelTheoremProducts}\emph{Generic Model Theorem (products)}.  Suppose $\P = \bigcup P_\alpha$ is a progressively closed product, and $\class{G}\subset \P$ is $\V$-generic for $\P$.  Then
$\V[\class{G}] \vDash \BGC$.
\end{Thm}

\begin{proof}Once again the partial Generic Model Theorem, Theorem \ref{T:PartialGenericModelTheorem}, applies and so it suffices to show that $\V[\class{G}]$ satisfies the axioms of power set and replacement.  For precise statements of these axioms, see the proof of Theorem \ref{T:GenericModelTheoremIterations}.  For the power set axiom, fix $a \in V[\class{G}]$.  Fix a regular cardinal $\delta > (|\tc(a)|)^{\V[\class{G}]}$ such that $\P$ factors at $\delta$ as $\P \cong \P_1 \times P_2$.   The isomorphism $\P \cong \P_1 \times P_2$ induces an isomorphism $\class{G}=\class{G}_1 \times G_2$.  As $a$ is coded as set of ordinals $\corners{a}\subset |\tc(a)|$ in $\V[\class{G}]$,  Lemma \ref{L:ClosedProducts} shows that $a \in \V[G_2]$, as is every subclass of $a$ from $\V[\class{G}]$.  Therefore it suffices to show that the power set of $a$ exists in $\V[G_2]$.  However, $\V[G_2]$ is a set-forcing extension of $\V$ and thus satisfies $\BGC$ by the Generic Model Theorem for sets, so the power set of $a$ exists in $\V[G_2]$.  Thus $\V[\class{G}]$ satisfies the power set axiom.

For the axiom of replacement, fix a function $\class{J} \in \V[\class{G}]$ and $\classdot{J}$ a $\P$-name for $\class{J}$.  By the axiom of choice in $\V[\class{G}]$ it suffices to consider domains of the form $\alpha$ an ordinal.  Fix $\alpha$, and choose $\delta>\alpha$ a regular cardinal such that $\P$ factors at $\delta$ as $\P \cong \P_1\times P_2$.  Let $\class{G}=\class{G}_1 \times G_2$ be the corresponding factorization of $\class{G}$.  I will construct a $\P$ name $\sigma$ such that $i_\class{G}(\sigma)$ is exactly $\{ \class{J}(\beta) \mid \beta < \alpha \}$.  The construction mirrors that in the proof of Lemma \ref{L:ClosedProducts}, except that the range of the function is no longer assumed to be member of $\V$.  Consider any condition $\<g,h>$ forcing that $\classdot{J}$ is a function.

\emph{Claim}.  If $\<g,h>$ forces $\classdot{J}$ defines a function on $\alpha$, and $\beta < \alpha$, then there is $p\leq g$ and $A \subset P_2$ such that $A$ is a maximal antichain below $h$ and for every $q\in A$ there is $\dot{a}_q \in V^\P$ such that $\<p,q> \forces \classdot{J}(\check{\beta}) = \dot{a}_q$.

The proof of the claim proceeds just as in the proof of Lemma \ref{L:ClosedProducts}, recursively constructing a descending sequence $\<p_\eta>$ below $g$, an antichain $\<q_\eta>$ below $h$, and a sequence of names $\<\dot{a}_\eta>$.  The claim shows that for every $\beta < \alpha$, the following class $D_\beta$ is open and dense below $p$.
\begin{eqnarray*}
	D_\beta = \left\{p  \mid  \exists A^{\beta}_p \subset P_2 \mbox{ max antichain below }h \mbox{ such that }\right.\\ 
	\left.\forall q \in A^{\beta}_p \, \exists \dot{a}^{\beta}_{p,q} \, \({ \<p,q> \forces \classdot{J}(\check{\beta}) = \dot{a}^{\beta}_{p,q} }) \right\} 
\end{eqnarray*}
Once again, closure of $\P_1$ shows that $\bigcap_{\beta<\alpha} D_\beta$ is dense.  Thus there is $\<g^\prime, h^\prime> \in \class{G}_1\times G_2$ such that $\<g^\prime, h^\prime> \forces \classdot{J}$ is a function, and $g^\prime \in D_\beta$ for every $\beta < \alpha$.

Working in $V[G_2]$, I define the $\P$-name $\sigma = \{ \<a^{\beta}_{g^\prime, h}, e(\<g^\prime,h>)> \mid h \in A^{\beta}_{g^\prime} \}$.  Note that although $\sigma$ may not be a ground model $\P$-name, $V[\class{G}]$ can still perform the valuation $i_\class{G}(\sigma)$.  Furthermore, $i_\class{G}(\sigma) = \{ i_\class{G}(a^{\beta}_{g^\prime, h}) \mid h \in G_2 \}$.  For each $\beta$ there will be a unique $h_\beta \in A^\beta_{g^\prime}$ such that $\<g^\prime,h_\beta> \in \class{G}_1\times G_2$, and it follows that $i_\class{G}(\classdot{J})(\beta)=i_\class{G}(a^{\beta}_{g^\prime, h_\beta})$.  Thus $i_\class{G}(\sigma)$ is exactly $\class{J}^{\prime\prime} \beta$ in $\V[\class{G}]$, and so $\V[\class{G}]$ satisfies replacement.  Thus $\V[\class{G}] \vDash \BGC$.  \end{proof}

The remaining two sections are general discussions of class iterations and class products.  In particular, they will show that a great many commonly used iterations and products satisfy the conditions given in Definitions \ref{D:IterationRestrictions} and \ref{D:ProductRestrictions}.

I begin with the theory of iterations.

\section{Proper class iterations}\label{SS:Iterations}

I will give the basic details of the theory of iterations below.  However, many elementary facts will be stated without proof.  For more details, I recommend to the reader the excellent exposition given in Kunen \cite{Kunen:Independence}, pp268--276.

Although the machinery of iterated forcing can be complicated, the basic idea is simple.  I begin by forcing with a partially ordered set $P\in \V$, obtaining a generic extension $\V[G]$.  I then force over $\V[G]$ with another partial order $Q \in \V[G]$, obtaining a further extension $\V[G][H]$, a two-stage forcing iteration of $\V$.  The idea behind iterated forcing is to show that $\V[G][H]$ is in fact a forcing extension of $\V$ by a single partial order $P*\dot{Q}$ in $\V$.  The definition of $P*\dot{Q}$ will be given below.  Note that iterations will be defined for sets only, as our proper class iterations will be built out of set forcing at each stage.

It will be useful to make some assumptions about the partial orders used in iterations.  First, I assume that every partial order $P$ has a largest element, denoted $1_P$.  Furthermore, while $P$ can be a partial order on any underlying set, it will be useful to have a uniform representative of the largest element.  Thus, I will consider only $P$ with $\emptyset \in P$ and $\emptyset = 1_P$.  Henceforth I will assume that `being a partial order' includes these restrictions.

Now suppose $\V[G]$ is a forcing extension by a partially ordered set $P$, and $Q\in \V[G]$ is a partial order.  I would like to fix a name $\dot{Q}$ for $Q$ in $V^P$, but in particular I would like a name such that $1_P \forces \dot{Q}$ is a partial order.  This is possible, using the theory of \emph{mixed names}, as $P$ is a set.  For more on mixed names and the mixing lemma see Kunen \cite{Kunen:Independence} p226.  If $1_P \forces \dot{Q}$ is a partial order, then I will call $\dot{Q}$ a \emph{$P$-name for a partial order}.  Note that $1_P \forces \check{\emptyset}$ is the top element of $\dot{Q}$.  The following definition is slightly nonstandard, with a modification allowing $P$ to be a subset of $P*\dot{Q}$.

\clearpage
\begin{Def}\rm\label{D:P*Q}Suppose $P$ is a partially ordered set and $\dot{Q}$ is a $P$-name for a partial order.  Then the partial order $P*\dot{Q}$ consists of $P$ together with all pairs $\<p,\tau>$ such that $p\in P$ and $\tau\in \dom(\dot{Q})$, with ordering
\begin{enumerate}
\parskip=0pt
	\item $p \leq \<q,\sigma> \mbox{ if and only if } p\leq q \mbox{ and } p \forces \check{\emptyset} \leq \sigma$,
	\item $\<p, \tau> \leq q \mbox{ if and only if } p\leq q$, and
	\item $\<p,\tau> \leq \<q,\sigma> \mbox{ if and only if } p\leq q \wedge p\forces \tau \leq \sigma$.
\end{enumerate}
The top element of $P*\dot{Q}$ is $1_{P*\dot{Q}} = 1_P = \emptyset$.\footnote{In the standard definition, $P*\dot{Q}$ consists only of pairs $\<p,\tau>$, and instead of $P\subset P*\dot{Q}$ one considers the embedding of $P$ into $P*\dot{Q}$ given by pairing each $p\in P$ with the top element of $\dot{Q}$,  $p \mapsto \<p,\check{\emptyset}>$.} 
\end{Def}

Note that $P*\dot{Q}$ may not be a partial order in the strict sense, e.g. there may be more than one name in $\dom(\dot{Q})$ for the element $1_Q$.  The main fact about $P*\dot{Q}$ is that forcing with it gives the same result as forcing with $P$ to obtain $\V[G]$ and then forcing with $i_G(\dot{Q})$ to obtain $\V[G][H]$.

\begin{Def}\rm\label{D:G*H}If $G$ is $\V$-generic for $P$ and $H\subset i_G(\dot{Q})$, then 
$$G*H = G \cup \{\<p,\tau> \mid p \in G \mbox{ and } i_G(\tau) \in H \}.$$
\end{Def}

\clearpage
\begin{Lemma}\label{L:P*Q}Suppose $P$ is a poset and $\dot{Q}$ is a $P$-name for a poset.  Then  $$P \subset_c P*\dot{Q}.$$ 
Furthermore, if $K$ is $\V$-generic for $P*\dot{Q}$, then setting $G=K \cap P$ and $H = \{i_G(\tau) \mid \exists p \, \<p,\tau> \in K \}$, it follows that
\begin{enumerate}
\parskip=0pt
\item $G$ is  $\V$-generic for $P$,
\item $H$ is $\V[G]$-generic for $i_G(\dot{Q})$,
\item $K= G*H$, and
\item $\V[K]=\V[G][H]$.
\end{enumerate}
\end{Lemma}

For the proof see Kunen \cite{Kunen:Independence} page 270, Theorem 5.5.  Now suppose I have a $P*\dot{Q}$-name $\tau$.  Let us take $\tau$ to be a name for a familiar object, such as an ordinal, (or a set of ordinals, or a partial order, etc.).  Using the Lemma above, it will often be the case that I wish to work in the partial extension $\V[G]$.  What happens to $\tau$ in $\V[G]$?  If $\tau$ is a $P*\dot{Q}$-name for an ordinal, then it seems that applying $i_G$ should yield, in $\V[G]$, a $Q$-name for an ordinal, where $Q = i_G(\dot{Q})$.  Strictly speaking, I cannot apply $i_G$ to $\tau$, for $\tau$ is not a $P$ name. However, there is a natural correspondence between `$P*\dot{Q}$-names' and `$P$-names for $Q$-names.'

\begin{Def}\rm\label{D:ReassociateP*Q} Recursively for $\classdot{A} \in \V^{P*\dot{Q}}$, let
\begin{eqnarray*}
	j(\classdot{A}) & = & \{ \<{ \<j(\sigma), p> , \<{ \{\<j(\sigma), p>, \<\dot{q}, p>\}, p}> , p}> \mid \<{ \sigma, \<p,\dot{q}> }> \in \classdot{A} \}\\
	& \cup & \{ \<{ \<j(\sigma), p> , \<{ \{\<j(\sigma), p>, \<\check{\emptyset}, p>\}, p}> , p}> \mid \<{ \sigma, p }> \in \classdot{A} \}
\end{eqnarray*}
\end{Def}

The recursion is performed first on sets in $V^{P*\dot{Q}}$ according to $P*\dot{Q}$-rank $\rho(\classdot{A})$ and then on proper class names.  The second part of the definition takes into account the fact that $P \subset P*\dot{Q}$.

Note that each $j(\classdot{A})$ is in $\V^P$, and for $G$ a $\V$-generic subset of $P$, the valuation $i_G(\classdot{A}) = \{\<i_G(j(\sigma), i_G(\dot{q})> \mid \exists p \in G \, \<{ \sigma, \<p,\dot{q}> }> \in \classdot{A} \} \cup \{\<i_G(j(\sigma), 1_Q> \mid \exists p \in G \, \<{ \sigma, p}> \in \classdot{A} \}$.  Verifying this fact is simply a matter of untangling the definition, paying careful attention to the application of the valuation function $i_G$ and recalling the definition of an ordered pair $\<a,b> = \{ a, \{a, b\} \}$.

\begin{Lemma}\label{L:ReassociateP*Q}Given $P * \dot{Q}$ and $K$ a $\V$-generic filter on $P*\dot{Q}$, suppose $G$ and $H$ are defined as in Lemma \ref{L:P*Q} above.  Then for every $\classdot{A} \in \V^{P*\dot{Q}}$, 
$$i_K(\classdot{A})=i_H(i_G(j(\classdot{A}))).$$
Furthermore, for any formula $\phi(x)$,
$$\({ \({ \<p,\dot{q}> \forces_{P*\dot{Q}} \phi(\classdot{A}) })^\V \mbox{ and } p \in G }) \rightarrow \({ i_G(\dot{q}) \forces_Q \phi(i_G(j(\classdot{A}))) })^{\V[G]}$$
where $Q = i_G(\dot{Q})$.
\end{Lemma}

I now turn my attention to iterations of arbitrary length.  Note that a $2$-stage iteration $P*\dot{Q}$ contains all pairs in $P \times \dom(\dot{Q})$.  An  $\alpha$-stage iteration will therefore consist of $\alpha$-sequences rather than pairs.  An important point to consider is whether one takes all $\alpha$-sequences, or only some subset.  This leads to the notion of \emph{support} of an iteration, which in turn relies on \emph{ideals}.

\begin{Def}\rm\label{D:Ideal}For any non-empty set $A$, a set $\I \subset \mathcal{P}(A)$ is an \emph{ideal on $A$} if and only if
\begin{enumerate}
\parskip=0pt
\item $\emptyset \in \I$,
\item $\forall X, Y \in \I \, \({ X\cup Y \in \I })$, and
\item \label{D:Ideal:downwardsclosed}$\forall X \in \I, \, \forall Y \, \({ Y \subset X \rightarrow Y \in \I })$.
\setcounter{saveenum}{\value{enumi}}
\end{enumerate}
I will also require that ideals contain all singletons (and thus all finite subsets),
\begin{enumerate}
\parskip=0pt
\setcounter{enumi}{\value{saveenum}}
\item $\forall a \in A, \, \{a\} \in \I$.
\end{enumerate}
\end{Def}

A technical problem involving the definition above will arise in forcing iterations.  Suppose $\I$ is an ideal in $\V$ and $\V[G]$ is a forcing extension of $\V$.  It may be the case that $\I$ is not longer an ideal in $\V[G]$, for if $X\in \I$ and the forcing adds a subset to $X$, then $\I$ will no longer satisfy condition \ref{D:Ideal:downwardsclosed} in $\V[G]$.  For this reason I will rely on a weaker notion which is preserved by forcing, that of a \emph{sub-ideal} on an ordinal.

\begin{Def}\rm\label{D:Subideal}For an ordinal $\alpha$, a set $\I \subset \mathcal{P}(\alpha)$ is a \emph{sub-ideal} if and only if
\begin{enumerate}
\parskip=0pt
\item $\emptyset \in \I$,
\item $\forall X, Y \in \I \, \({ X\cup Y \in \I })$,
\item $\forall X \in \I, \, \forall \beta<\alpha \, \({ X\cap\beta \in \I })$, and
\item $\forall a \in A, \, \{a\} \in \I$.
\end{enumerate}
\end{Def}

\begin{Def}\rm\label{D:IdealCompleteness}For $\I$ an ideal or sub-ideal and  $\kappa$ a regular cardinal
\begin{enumerate}
\parskip=0pt
\item $\I$ is \emph{$\kappa$-complete} if and only if for all $Z\subset \I$ such that $|Z|<\kappa$, the union $\bigcup Z \in \I$.
\item $\I$ is \emph{$\kappa$-subcomplete} if and only if for all $Z\subset \I$ such that $|Z|<\kappa$, there is $Y\in \I$ such that $\bigcup Z \subset Y$.
\end{enumerate}
\end{Def}

Note that if we take $\mathcal{P}(A)$ to be ordered by the $\subset$ relation, then an ideal on $A$ is the dual notion of a filter  on the ordering $\<A,\subset>$ (Definition \ref{D:POClassFilters}).  Also note that it is usually assumed $A\notin \I$.  However, certain commonly used ideals in forcing iterations fail to satisfy this property when restricted to certain stages of the iteration, and so it is omitted.

\clearpage
\begin{Def}\rm\label{D:alphaStageIteration}For $\alpha$ an ordinal and $\I$ a sub-ideal on $\alpha + 1$, an \emph{$\alpha$-stage forcing iteration with supports in $\I$} consists of two sequences
$$\< P_\beta \mid \beta \leq \alpha> \mbox{ and } \<\pi_\beta \mid \beta < \alpha>$$
satisfying:
\begin{enumerate}
\parskip=0pt
\item Each $P_\beta$ is a partial order consisting of functions $p$ with domain a subset of $\beta$.
\item If $p\in P_\beta$ and $\xi < \beta$, then $p\upharpoonright \xi \in P_\xi$.
\item Each $\pi_\beta$ is a $P_\beta$-name for a partial order.
\item If $p \in P_\beta$ and $\<\xi, \tau>\in p$, then $\tau \in \dom(\pi_\xi)$.
\item For each $\beta$, the top element $1_\beta$ is the empty sequence $\emptyset$.
\item If the \emph{support} of $p \in P_\beta$ is defined $\supt(p)=\dom(p)$, then for every $p$, $\supt(p)\in\I$.
\item $P_0 = \{\emptyset\}$.
\setcounter{saveenum}{\value{enumi}}
\end{enumerate}
At successor stages $\beta+1$, the iteration satisfies:
\begin{enumerate}
\parskip=0pt
\setcounter{enumi}{\value{saveenum}}
\item A function $p$ is in $P_{\beta+1}$ if and only if $p$ is a function with domain $\subset \beta + 1$, the restriction $p\upharpoonright \beta \in P_\beta$, and $\beta\in \dom(p) \rightarrow p(\beta)\in \dom(\pi_\beta)$.
\item For $p$ and $q$ in $P_{\beta+1}$, the ordering is given by $p\leq_{P_{\beta+1}} q$ if and only if $p\upharpoonright\beta \leq_{P_\beta} q\upharpoonright\beta$ and one of the following holds.
	\begin{enumerate}
	\parskip=0pt
	\item $\beta \in \dom(p) \cap \dom(q)$ and $p\upharpoonright\beta \forces_{P_\beta} p(\beta) \leq_{P_{\beta+1}} q(\beta)$.
	\item $\beta \in \dom(p) \setminus \dom(q)$.
	\item $\beta \in \dom(q) \setminus \dom(p)$  and $p\upharpoonright\beta \forces_{P_\beta} \check{\emptyset} \leq_{P_{\beta+1}} q(\beta)$.
	\end{enumerate}
\setcounter{saveenum}{\value{enumi}}
\end{enumerate}
At limit stages $\beta$, the iteration satisfies:
\begin{enumerate}
\parskip=0pt
\setcounter{enumi}{\value{saveenum}}
\item A function $p$ is in $P_\beta$ if and only if $p$ is a partial function with domain $\subset \beta$, the support of $p$ is in $\I$, and  $p\upharpoonright \xi \in P_\xi$ for every $\xi < \beta$.
\item For $p$ and $q$ in $P_\beta$, the ordering is given by $p\leq_{P_\beta} q$ if and only if $p\upharpoonright\xi \leq_{P_\xi} q\upharpoonright\xi$ for every $\xi < \beta$.
\end{enumerate}
\end{Def}

Note that the definition of $P_{\beta+1}$ is almost exactly given by $P_\beta * \pi_\beta$.  The difference is technical, as $P_{\beta+1}$ consists of partial $\beta+1$-sequences, but $P_\beta * \pi_\beta$ consists of pairs whose first member is a partial $\beta$-sequence.  However, there is a canonical isomorphism $P_{\beta+1} \cong P_\beta * \pi_\beta$ obtained by simply rearranging each member of $P_\beta * \pi_\beta$ into a partial $\beta+1$ sequence.  This idea will be generalized in Definition \ref{D:ReassociateIteration} below.  Before proceeding, however, I note that every initial part of an iteration embeds completely in to the later stages.

\begin{Lemma}\label{L:alphaStageIteration}Suppose $\< P_\beta \mid \beta \leq \alpha> \mbox{ and } \<\pi_\beta \mid \beta < \alpha>$ form an $\alpha$-stage forcing iteration with supports in $\I$.  Then $P_\xi \subset_c P_\beta$ for all $\xi<\beta\leq \alpha$.
\end{Lemma}

This means, in particular, that a generic for $P_\alpha$ gives generics for all of the $P_\beta$.  It also shows that $P_\alpha = \bigcup P_\beta$ satisfies the conditions of being a chain of complete subposets (Definition \ref{D:CCP}), except that the union is taken only over $\alpha$ and not over all of the ordinals.  The generalization to $\ORD$-length iterations appears below.

Forcing arguments using iterations often involve working in some partial extension obtained by forcing over $P_\beta$ for some $\beta$.  I would like to have a sensible way of describing what happens to the remainder of the iteration in such a case.  In fact, for any $\beta<\alpha$ there is a $P_\beta$-name for a partial order $\dot{P}_{[\beta,\alpha)}$ such that $P_\alpha \cong P_\beta * \dot{P}_{[\beta,\alpha)}$ by a natural isomorphism in the spirit of Definition \ref{D:ReassociateP*Q}.  The general definition is given next.

\begin{Def}\rm\label{D:ReassociateIteration}Suppose $\< P_\beta \mid \beta \leq \alpha> \mbox{ and } \<\pi_\beta \mid \beta < \alpha>$ form an $\alpha$-stage forcing iteration with supports in $\I$.  Then for $\xi < \beta \leq \alpha$, define $j_{\xi,\beta}:P_\beta \rightarrow V^{P_\xi}$ as follows.
For $p\in P_\beta$, let $p_0 = p \upharpoonright \xi$ and $p_1 = p \upharpoonright [\xi,\beta)$.  If $p_1$ is empty, then $j_{\xi,\beta}(p)=\emptyset$.  Otherwise,
\begin{eqnarray*}
	j_{\xi,\beta}(p)&= &\{ \<{  \op(\check{\gamma},j_{\xi,\beta}(\sigma)), p_0 }> \mid \<{ \gamma, \sigma }> \in \dom(p_1) \}
\end{eqnarray*}
Let $\dot{P}_{[\xi,\beta)} = \{j_{\xi,\beta}(p) \mid p\in P_\beta \mbox{ and } p \notin P_\xi\}$.
\end{Def}

This map is well-behaved with respect to the inclusion $P_\xi \subset_c P_\beta$.

\clearpage
\begin{Lemma}\label{L:ReassociateIterationsMapCommutes}For $\P$ and $j_{\xi,\beta}$ defined as above, for $\xi < \beta \leq \alpha$ the map $p \mapsto \<p_0, j_{\xi,\beta}(p_1)>$ defines  an isomorphism $P_\beta \cong P_\xi * \dot{P}_{[\xi,\beta)}$.  Furthermore, for $\xi < \beta < \beta^\prime \leq \alpha$ the isomorphisms commute with the inclusion map,
 \[
 \begin{CD}
P_\beta @>\subset_c>> P_{\beta^\prime} \\
@V{\cong}VV @V{\cong}VV\\
P_\alpha * \dot{P}_{[\alpha,\beta)} @> \subset_c>> P_\alpha * \dot{P}_{[\alpha,\beta^\prime)} \\
 \end{CD}
 \]
\end{Lemma}

Note that $j_{\xi,\beta}$ induces a map between $P_\beta$ names and `$P_\xi$ names for $\dot{P}_{[\xi,\beta)}$ names.'  Abusing notation, I will also denote this map $j_{\xi,\beta}$.  The next lemma essentially says that these maps are well-behaved.

\begin{Lemma}\label{L:ReassociateIteration}Suppose $\< P_\beta \mid \beta \leq \alpha> \mbox{ and } \<\pi_\beta \mid \beta < \alpha>$ form an $\alpha$-stage forcing iteration with supports in $\I$ and the maps $j_{\xi,\beta}$ for $\xi<\beta\leq\alpha$ are defined as above.  In addition, if $G$ is $\V$-generic for $P_\alpha$, then let $G_\beta$ = $G \cap P_\beta$ for each $\beta<\alpha$.  For all $\xi < \beta < \alpha$, let $G_{[\xi,\beta)} = \{ i_{G_\xi}(j_{\xi,\beta}(p)) \mid p \in G_\beta \setminus G_\xi \}$.  Then
\begin{enumerate}
\parskip=0pt
\item $G_\xi$ is $\V$-generic for $P_\xi$,
\item $G_{[\xi,\beta)}$ is $\V[G_\xi]$-generic for $i_{G_\xi}(\dot{P}_{[\xi,\beta)})$,
\item For all $\tau \in \V^{P_\beta}$, $i_{G_\beta}(\tau) = i_{G_{[\xi,\beta)}}(i_{G_\xi}(j_{\xi,\beta}(\tau)$, and 
\item $\V[G_\beta] = \V[G_\xi][G_{[\xi,\beta)}]$.
\end{enumerate}
Furthermore, for $p\in P_\beta$, $\tau \in V^{P_\beta}$ and any formula $\phi(x)$,
$$\({ \({ p \forces_{P_\beta} \phi(\tau) })^\V \mbox{ and } (p\upharpoonright \xi) \in G_\xi }) \rightarrow \({ i_{G_\xi}(j_{\xi,\beta}(p)) \forces_Q \phi(i_{G_\xi}(j_{\xi,\beta}(\tau))) })^{\V[G_\xi]}$$
where $Q = {i_{G_\xi}(\dot{P}_{[\xi,\beta)})}$. 
Finally, in $\V[G_\xi]$, the sequences $\< i_{G_\xi}(\dot{P}_{[\xi,\beta)}) \mid \beta \in [\xi,\alpha]>$ and $\<i_{G_\xi}(j_{\xi,\beta}(\pi_\beta) \mid \beta \in [\xi,\alpha)>$ form an $\alpha$-stage forcing construction with supports in $\I$.
\end{Lemma}

Finally, the closure property of an iteration depends on the closure of each factor together with the subcompleteness of the ideal $\I$.

\begin{Lemma}\label{L:ClosedIterations}Suppose $\< P_\beta \mid \beta \leq \alpha> \mbox{ and } \<\pi_\beta \mid \beta < \alpha>$ form an $\alpha$-stage forcing iteration with supports in $\I$.  If $\kappa$ is a regular cardinal, $\I$ is $\kappa$-subcomplete, and  $P_\beta \forces ``\pi_\beta$ is $<\!\kappa$-closed'' for all $\beta<\alpha$, then $P_\beta$ is $<\!\kappa$-closed for all $\beta \leq \alpha$.
\end{Lemma}

For proper class iterations, we generalize the notions above to allow $\ORD$-many stages.

\begin{Def}\rm\label{D:ORDIdeal}A \emph{sub-ideal on $\ORD$} is a class $\I$ consisting of sets of ordinals and satisfying the four conditions of Definition \ref{D:Subideal}.  For a regular cardinal $\kappa$, the definitions of $\kappa$-complete and $\kappa$-subcomplete remain the same.
\end{Def}

\clearpage
\begin{Def}\rm\label{L:ORDStageIteration}For $\I$ a sub-ideal on $\ORD$, an \emph{$\ORD$-stage forcing iteration with supports in $\I$} consists of two classes $\< P_\alpha \mid \alpha \in \ORD>$ and $\<\pi_\alpha \mid \alpha \in \ORD>$ such that for every ordinal $\alpha$, the sequences $\< P_\beta \mid \beta \leq \alpha>$ and $\<\pi_\beta \mid \beta < \alpha>$ form an $\alpha$-stage forcing iteration with supports in $\I$.
\end{Def}

Lemma \ref{L:alphaStageIteration} claims that, under this definition, such an iteration is in fact a chain of complete subposets.

\begin{Lemma}\label{L:ORDStageIterationIsChainOfCompleteSubposets}If $\< P_\alpha \mid \alpha \in \ORD>$ and $\<\pi_\alpha \mid \alpha \in \ORD>$ is an $\ORD$-stage forcing iteration with supports in $\I$, then $\P = \bigcup_{\alpha\in\ORD} P_\alpha$ is a chain of complete subposets.
\end{Lemma}

\begin{proof}The Lemma follows directly from Lemma \ref{L:alphaStageIteration} applied at each stage $\alpha$.
\end{proof}

Finally, I would like to give conditions under which an $\ORD$-stage forcing iteration is a progressively closed iteration (Definition \ref{D:IterationRestrictions}).  Showing that the tails of the iteration remain closed requires attention both to closure in the stages of forcing and to the completeness of the support.  The most common support is $\emph{Easton support}$.

\clearpage
\begin{Def}\rm\label{D:EastonSupport}An iteration is said to have \emph{Easton support} if it uses the ideal $\I$ on $\ORD$ consisting of all sets of ordinals $A$ such that  $|A\cap \kappa| < \kappa$ for every inaccessible cardinal $\kappa$.  Intuitively, Easton support uses bounded support at inaccessibles and full support at all other limits.
\end{Def}

\begin{Lemma}\label{L:ProgressivelyClosedIterations}If $\< P_\alpha \mid \alpha \in \ORD>$ and $\<\pi_\alpha \mid \alpha \in \ORD>$ is an $\ORD$-stage forcing iteration with Easton support, and for arbitrarily large regular cardinals $\delta$ there exists $\alpha \geq \delta$ such that $P_\beta \forces ``\pi_\beta$ is $<\check{\delta}$-closed'' for every $\beta \geq \alpha$, then $\P = \bigcup P_\alpha$ is a progressively closed iteration.
\end{Lemma}

\begin{proof}   That $\P$ is chain of complete subposets was given in Lemma \ref{L:ORDStageIterationIsChainOfCompleteSubposets}.  For the closure condition, fix $\delta$ and fix $\alpha$ as in the lemma.  For every $\beta > \alpha$, Lemma \ref{L:ReassociateIteration} shows that $$P_\alpha \forces \({ \dot{P}_{[\alpha,\beta)} \forces j_{\alpha,\beta}(\pi_\beta) \mbox{ is } <\!\check{\delta} \mbox{-closed} }).$$
Taking $\alpha\geq\delta$ ensures that Easton support is $\delta$-complete when restricted to sets `lying above $\alpha$,' that is, sets $A$ such that $A\cap \alpha = \emptyset$.  Thus Lemma \ref{L:ClosedIterations} applies and so 
$$P_\alpha \forces \({ \dot{P}_{[\alpha,\beta)} \mbox{ is } <\!\check{\delta} \mbox{-closed} }).$$
That the commutativity condition is satisfied follows from Lemma \ref{L:ReassociateIterationsMapCommutes}.  Finally,
$$\alpha < \lambda < \lambda^\prime \rightarrow P_\alpha \forces \dot{P}_{[\alpha,\lambda)} \subset_c  \dot{P}_{[\alpha,\lambda^\prime)}$$
is a consequence the final assertion of Lemma \ref{L:ReassociateIteration}, stating that in the extension by $P_\alpha$, the tail forcing is itself an iteration. \end{proof}

Thus, for example, the Easton support iteration adding a single subset to each regular $\kappa$ is a progressively closed iteration.  So is any variation on this theme, adding subsets only to certain $\kappa$ or adding many subsets to $\kappa$, in any definable combination.  Thus the $\GCH$ coding used to force the Ground Axiom, as well as the Kurepa coding used to force $\GCH + \GA$, are progressively closed iterations and so forcing with them preserves $\BGC$.

\section{Proper class products}\label{SS:Products}

A product is simply a special case of an iteration in which all $\pi_\beta$ are in fact partial orders in $\V$.  For example, if $P$ and $Q$ are posets, then the product $P\times Q$ can be viewed simply as the iteration $P*\check{Q}$, and the material of the previous section applies.  However, there are two important facts which contribute to making the theory of products more than a special case of the theory of iterations.  First, although both $P$ and $Q$ may be $<\!\kappa$-closed in $V$, it does not follow that $P \forces \check{Q}$ is $<\!\kappa$-closed.  In fact, this is often not true.  This will prevent us from viewing a product of closed posets as a progressively closed iteration. Second, since both $P$ and $Q$ are in $\V$ there is no clear notion of which forcing should take place first.  Intuitively, forcing with  $P*\check{Q}$ should give the same result as forcing with $Q*\check{P}$.  This commutative property, $P\times Q \cong Q\times P$, lies at the heart of many forcing argument involving products. I will give a general definition of $P\times Q$ and state the basic facts about such a product in the case where both are sets.  I will then generalize to allow one of them to be a proper class, as in the case of progressively closed products, Definition \ref{D:ProductRestrictions}.  Once again I assume that all partial orders have a top element $\emptyset$.

\begin{Def}\rm\label{D:PxQ}Suppose $\P$ and $\Q$ are partial orders.  Then the partial order $\P\times \Q$ consists of all pairs $\<p,q>$ such that $p\in P$ and $q\in Q$, with ordering $\<p_0, q_0> \leq \<p_1, q_1>$ if and only if $p_0 \leq p_1$ and $q_0 \leq q_1$.
\end{Def}

\begin{Def}\rm\label{D:GxH}If $\P$ and $\Q$ are partial orders, $\class{G}\subset \P$ and $\class{H}\subset \Q$, then 
$$\class{G}\times \class{H} = \{\<p,q>\in \P\times \Q \mid p\in \class{G} \wedge q\in \class{H} \}.$$
\end{Def}

It is a matter of checking definitions to see that $\class{G}\times \class{H}$ is a filter on $\P\times\Q$ if and only if $\class{G}$ is a filter on $\P$ and $\class{H}$ is a filter on $\Q$.  I would like to replace `filter' with `generic filter' in this statement, but this will require additional restrictions on $\P$ and $\Q$ (Lemma \ref{L:PxQclasses}).

In the case where both $P$ and $Q$ are sets, the following facts are standard.

\clearpage
\begin{Lemma}\label{L:PxQ}Suppose $P$ and $Q$ are posets. Then the map $p \mapsto \<p,1_Q>$ is a complete embedding of $P$ into $P\times Q$.  Furthermore, suppose $G\subset P$ and $H\subset Q$.  Then the following are equivalent:
\begin{enumerate}
\parskip=0pt
\item $G\times H$ is $\V$-generic for $P\times Q$.
\item $H\times G$ is $\V$-generic for $Q\times P$.
\item $G$ is  $\V$-generic for $P$ and $H$ is $V[G]$-generic for $Q$.
\item $H$ is $\V$-generic for $P$ and $G$ is $V[H]$-generic for $Q$.
\end{enumerate}
Furthermore, if any of the above hold, then $\V[G\times H]=\V[G][H]=\V[H][G]=\V[H\times G]$.
\end{Lemma}

The next goal is the generalization of this lemma to apply in the case of progressively closed products (Definition \ref{D:ProductRestrictions}).  

\clearpage
\begin{Lemma}\label{L:PxQclasses}Suppose $\delta$ is a regular cardinal, $Q$ is a poset with the $\delta^+$-c.c., $\P$ is a chain of complete subposets and $\P$ is $\leq\!\delta$-closed.  Fix $H \subset Q$ and $\class{G}\subset \P$. Then the following are equivalent:
\begin{enumerate}
\parskip=0pt
\item \label{L:PxQclasses:V[HxG]} $H\times \class{G}$ is $\V$-generic for $Q\times \P$.
\item \label{L:PxQclasses:V[GxH]} $\class{G}\times H$ is $\V$-generic for $\P\times Q$.
\item \label{L:PxQclasses:V[H][G]} $H$ is $\V$-generic for $Q$ and $\class{G}$ is $V[H]$-generic for $\P$.
\item \label{L:PxQclasses:V[G][H]} $\class{G}$ is  $\V$-generic for $\P$ and $H$ is $V[\class{G}]$-generic for $Q$.
\end{enumerate}
Furthermore, if any of the above hold, then $\V[\class{G}\times H]=\V[\class{G}][H]=\V[H][\class{G}]=\V[H\times \class{G}]$.
\end{Lemma}

\begin{proof}That \ref{L:PxQclasses:V[HxG]}$\iff$\ref{L:PxQclasses:V[GxH]} follows from the  obvious isomorphism $Q\times \P \cong \P\times Q$.  I will show \ref{L:PxQclasses:V[HxG]}$\iff$\ref{L:PxQclasses:V[H][G]} and \ref{L:PxQclasses:V[GxH]}$\iff$\ref{L:PxQclasses:V[G][H]}.

For \ref{L:PxQclasses:V[HxG]}$\rightarrow$\ref{L:PxQclasses:V[H][G]}, I will begin by showing $H$ is $\V$-generic for $Q$. Fix $D_Q$ a dense subset of $Q$.  Consider $\class{D} = D_Q \times \P$.  Clearly $\class{D}$ is dense in $Q\times \P$, and so $H\times \class{G}$ intersects it.  Thus $H$ intersects $D_Q$, and so $H$ is $\V$-generic for $Q$.  Therefore it makes sense to consider the forcing extension $\V[H]$.  I will next show that $\class{G}$ is $\V[H]$-generic for $\P$.  Since $\V \subset \V[H]$, it follows that $\P$ is a class in $\V[H]$.  Now suppose that $\class{D} \in \V[H]$ is a dense subclass of $\P$, and fix $\classdot{D}\in \V^Q$ a name for $\class{D}$.  There is $q\in H$ forcing that $\classdot{D}$ is dense in $\P$.    I claim that $\class{B} = \{\<q^\prime,p^\prime> \in Q\times \P \mid q^\prime \forces \check{p}^\prime\in\classdot{D} \}$ is dense below $\<q,1_\P>$ in $Q\times \P$.  To see this, fix any $\<q_0, p_0> < \<q,1_\P>$ and note that $q_0 \forces \classdot{D}$ is dense in $\P$.  Thus $q_0 \forces \exists p \leq p_0 \, (p \in \classdot{D})$, and so there is $q_1 \leq q_0$ and $p_1 \leq p_0$ such that $q_0 \forces \check{p}_1 \in \classdot{D}$.  This shows that $\class{B}$ is dense below $\<q,1_\P>$ in $Q\times \P$.  Since $H\times \class{G}$ is $\V$-generic, and since it contains $q$, it must contain some member $\<q^\prime, p^\prime>$ of $\class{B}$.  Thus $p^\prime \in \class{G} \cap \class{D}$, and it follows that $\class{G}$ is $V[H]$-generic for $\P$.

For the reverse direction, \ref{L:PxQclasses:V[H][G]}$\rightarrow$\ref{L:PxQclasses:V[HxG]}, consider $\class{D}$ a dense subclass of $Q\times \P$.  Without loss of generality I can assume that $\class{D}$ is open (Lemma \ref{D:POGenericFilters}).  In $\V[H]$, let $\class{D}_\P = \{p\in\P \mid \exists q\in H \, \<q,p> \in \class{D}\}$.  I claim that $\class{D}_\P$ is dense in $\P$.  Fix $p_0\in\P$, and working in $\V$ define $D_Q = \{q \in Q \mid \exists p\leq p_0 \, \<q,p>\in \class{D}\}$.  It is clear that $D_Q$ is dense in $Q$, and so there must be $q_1 \in D_Q \cap H$ and $p_1\leq p_0$ such that $\<q_1,p_1> \in \class{D}$.  Thus $p_1\in \class{D}_\P$, and so $\class{D}_\P$ is dense in $\P$ in $\V[H]$.  Therefore there is $p\in \class{G}\cap \class{D}_\P$ with an associated $q\in H$ such that $\<q,p> \in \class{D}$.  This shows that $\class{D} \cap H\times \class{G}$ is nonempty, and it follows that $H\times \class{G}$ is $\V$-generic for $Q\times \P$.

The equivalence \ref{L:PxQclasses:V[GxH]}$\iff$\ref{L:PxQclasses:V[G][H]} is proved in a similar fashion, with the strategies adapted slightly to accomodate the reversal of roles of the set and class forcing.  For the forward direction, $\V$-genericity of $\class{G}$ follows as before.  Note also that $H$ is $\V$-generic for $Q$.  To see that $H$ is $\V[\class{G}]$-generic for $Q$, first observe that $\P$ is $\leq\!\delta$-closed and so adds no maximal antichains to $Q$ as $Q$ has the $\delta^+$-c.c.. Thus $H$ is also $\V[\class{G}]$-generic for $Q$.  The reverse direction follows as in the proof of  \ref{L:PxQclasses:V[H][G]}$\rightarrow$\ref{L:PxQclasses:V[HxG]}. \end{proof}

This concludes the main facts used in the proof of the Generic Model Theorem for products (Theorem \ref{T:GenericModelTheoremProducts}).  It remains to show that the products used throughout this paper are, in fact, progressively closed products.  For this, a general definition of class products is necessary.

\begin{Def}\rm\label{D:GeneralProducts}Suppose $\I$ is a sub-ideal on $\ORD$ and $\<Q_\alpha \mid \alpha\in\ORD>$ is a class such that each $Q_\alpha$ is a poset.  Then the \emph{product $\P$ of $\<Q_\alpha>$ with supports in $\I$}, denoted $\prod_{\alpha\in\ORD} Q_\alpha$, consists of all functions $p$ satisfying
\begin{enumerate}
\parskip=0pt
\item $\dom(p) \in \I$, and 
\item $p(\alpha)\in Q_\alpha$ for every $\alpha\in \dom(p)$,
\end{enumerate}
with ordering
$$p\leq q \iff \forall \alpha \in \({ \dom(p)\cap \dom(q) }), \, p(\alpha) \leq q(\alpha).$$
\end{Def}

\begin{Lemma}\label{L:GeneralProductIsChainofCompleteSubposets}Suppose $\P=\<Q_\alpha \mid \alpha\in\ORD>$ is a product with supports in $\I$.  For each $\alpha$ let $P_\alpha = \{p\in \P \mid \dom(p)\subset \alpha\}$, with ordering inherited from $\P$.  Then $\P = \bigcup_{\alpha\in\ORD} P_\alpha$ is a chain of complete subposets.
\end{Lemma}

\begin{proof}The top element of each $P_\alpha$ is the empty sequence $\emptyset$, conforming to the restriction given in the introduction to Section \ref{SS:Iterations}.  The definition shows that the ordering on $P_\alpha$ agrees with the ordering on $\P$, so  $p \leq_{P_\alpha} q \rightarrow p \leq_\P q$.  Furthermore $p \parallel_{P_\alpha} q \iff p \parallel_\P q$ for all $p$ and $q$ in $P_\alpha$, for if there was a condition $r\in\P$ lying below both $p$ and $q$, then $r\upharpoonright \alpha$ would also lie below both $p$ and $q$.  This shows that the notions of ordering $\leq$, compatibility $\parallel$, and incompatibility $\perp$ are unambiguous among all $P_\alpha$.  Fix $\beta < \alpha$.  Clearly $P_\beta$ is a subposet of $P_\alpha$.  To show completeness, fix $A\subset P_\beta$ a maximal antichain and suppose there is $r\in P_\alpha$ such that $r \perp p$ for every $p \in A$. Fix $p\in A$.  If $r \perp p$, then  $r\upharpoonright \beta \perp_{P_\alpha} p$, because compatibility and incompatibility are completely determined by the parts of $p$ and $r$ lying in their common domain.  Thus $(r\upharpoonright \beta) \perp p$ for every $p\in A$, contradicting that $A$ was maximal in $P_\beta$.  Thus $P_\beta \subset_c P_\alpha$, and so $\P=\bigcup P_\alpha$ is a chain of complete subposets.
\end{proof}

Finally, an additional fact relating to closure is required to achieve a progressively closed product.

\begin{Lemma}\label{L:ClosedProducts2}Suppose $\P=\bigcup P_\alpha$ is the product of $\<Q_\alpha \mid \alpha \in \ORD>$ with supports in $\I$.  Suppose $\delta$ is a regular cardinal, $\I$ is $\delta$-subcomplete, and  $Q_\alpha$ is $<\!\delta$-closed for all $\alpha \in \ORD$.  Then $P_\alpha$ is $<\!\kappa$-closed for all $\alpha$, as is $\P$. 
\end{Lemma}

\clearpage
\begin{Lemma}\label{L:ProgressivelyClosedProducts}Suppose $\P=\bigcup P_\alpha$ is the product of $\<Q_\alpha \mid \alpha \in \ORD>$ with Easton support, and for arbitrarily large regular $\delta$ there exists $\alpha\geq\delta$ such that $P_\alpha$ has the $\delta^+$-c.c. and  $Q_\beta$ is $\leq\!\delta$-closed for every $\beta \geq \alpha$, then $\P = \bigcup P_\alpha$ is a progressively closed product.
\end{Lemma}

\begin{proof}   That $\P=\bigcup P_\alpha$ is a chain of complete subposets was shown above.  Fix $\delta$ and $\alpha$ as in the Lemma.  To see that $\P$ factors at $\delta$, for each $p\in\P$ let $p_2 = p\upharpoonright \alpha$ and $p_1 = p \setminus p_2$.  The  map $p \mapsto \<p_1,p_2>$ is an isomorphism $\P \cong \P_1 \times P_2$.  Note that $P_2 = P_\alpha$, which has the $\delta^+$-c.c. by hypothesis.  It is a matter of checking definitions to see that $\P_1$ is exactly the product of $\<Q_\beta \mid \beta \geq \alpha>$ with Easton support, and so $\P_1$ is a chain of complete subposets.  Furthermore, Easton support restricted to sets `above $\alpha$' will be $\leq\!\alpha$-subcomplete and hence $\leq\!\delta$-subcomplete, and so Lemma \ref{L:ClosedProducts2} shows that $\P_1$ is $\leq\!\delta$-closed.  Thus $\P$ factors at $\alpha$, and since this holds for arbitrarily large $\alpha$, it follows that $\P$ is a progressively closed product. \end{proof}

Any product that adds Cohen subsets to $\kappa$ at stage $\kappa$, for example the forcing used by Easton to control the continuum function, will be progressively closed.  Indeed, all product forcing used in this paper, including the forcing that produces models of $\neg\BA+\VHOD$, is progressively closed.

\end{appendices}

\backmatter

\addcontentsline{toc}{chapter}{Bibliography}
\bibliographystyle{alpha}
\bibliography{MathBiblio,HamkinsBiblio}

\end{document}